\documentclass[10pt]{article} 
\usepackage{hyperref}
\usepackage{amssymb, amsmath}
\usepackage{amsmath}
\usepackage{latexsym}
\usepackage{amssymb}

\usepackage{amsthm}

\font\sevengoth=eufm7 at 6pt
\newfam\gothfam
\scriptfont\gothfam=\sevengoth

\newcommand{\mlabel}[1]{\marginpar{#1}\label{#1}}

\newcommand{\fB}{{\mathfrak B}}

\newcommand{\fS}{{\mathfrak S}}

\newcommand{\g}{{\mathfrak g}}

\newcommand{\fg}{{\mathfrak g}}
\newcommand{\fh}{{\mathfrak h}}

\newcommand{\fo}{{\mathfrak o}}
\newcommand{\fq}{{\mathfrak q}}

\renewcommand{\:}{\colon}
\newcommand{\1}{\mathbf{1}}

\newcommand{\cD}{\mathcal{D}}
\newcommand{\cE}{\mathcal{E}}
\newcommand{\cF}{\mathcal{F}}
\newcommand{\cG}{\mathcal{G}}
\newcommand{\cH}{\mathcal{H}}
\newcommand{\cK}{\mathcal{K}}
\newcommand{\cL}{\mathcal{L}}
\newcommand{\cM}{\mathcal{M}}
\newcommand{\cN}{\mathcal{N}}
\newcommand{\cO}{\mathcal{O}}

\newcommand{\cS}{\mathcal{S}}

\newcommand{\eset}{\emptyset}

\renewcommand{\phi}{\varphi}
\newcommand{\dd}{{\tt d}}

\newcommand{\subeq}{\subseteq}

\newcommand{\into}{\hookrightarrow}

\newcommand{\shalf}{{\textstyle{\frac{1}{2}}}}

\newcommand{\N}{{\mathbb N}}
\newcommand{\Z}{{\mathbb Z}}
\newcommand{\R}{{\mathbb R}}
\newcommand{\C}{{\mathbb C}}

\newcommand{\K}{{\mathbb K}}

\newcommand{\T}{{\mathbb T}}

\newcommand{\bS}{{\mathbb S}}

\renewcommand{\hat}{\widehat}

\renewcommand{\tilde}{\widetilde}

\renewcommand{\L}{\mathop{\bf L{}}\nolimits}

\newcommand{\GL}{\mathop{{\rm GL}}\nolimits}

\newcommand{\SO}{\mathop{{\rm SO}}\nolimits}

\newcommand{\OO}{\mathop{\rm O{}}\nolimits}

\newcommand{\U}{\mathop{\rm U{}}\nolimits}

\renewcommand{\Re}{\mathop{{\rm Re}}\nolimits}

\newcommand{\pr}{\mathop{{\rm pr}}\nolimits}

\newcommand{\Herm}{\mathop{{\rm Herm}}\nolimits}

\newcommand{\Heis}{\mathop{{\rm Heis}}\nolimits}

\newcommand{\Aut}{\mathop{{\rm Aut}}\nolimits}

\newcommand{\id}{\mathop{{\rm id}}\nolimits}

\renewcommand{\dim}{\mathop{{\rm dim}}\nolimits}

\newcommand{\supp}{\mathop{{\rm supp}}\nolimits}

\newcommand{\Spann}{\mathop{{\rm span}}\nolimits}
\newcommand{\ev}{\mathop{{\rm ev}}\nolimits}

\newcommand{\Rarrow}{\Rightarrow}
 
\newcommand{\oline}{\overline}

\newcommand{\la}{\langle}
\newcommand{\ra}{\rangle}

\newcommand{\Mot}{{\rm Mot}}
\newcommand{\up}{\mathop{\uparrow}}

\newcommand{\res}{\vert}

\newcommand{\Spec}{{\rm Spec}}

\newcommand{\ssssarr}{\hbox to 15pt{\rightarrowfill}}
\newcommand{\sssarr}{\hbox to 20pt{\rightarrowfill}}
\newcommand{\ssarr}{\hbox to 30pt{\rightarrowfill}}
\newcommand{\sarr}{\hbox to 40pt{\rightarrowfill}}
\newcommand{\arr}{\hbox to 60pt{\rightarrowfill}}
\newcommand{\larr}{\hbox to 60pt{\leftarrowfill}}
\newcommand{\Arr}{\hbox to 80pt{\rightarrowfill}}

\def\theoremname{Theorem}
\def\propositionname{Proposition}
\def\corollaryname{Corollary}
\def\lemmaname{Lemma}
\def\remarkname{Remark}
\def\conjecturename{Conjecture} 

\def\definitionname{Definition}
\def\exercisename{Exercise}
\def\examplename{Example}
\def\examplesname{Examples}
\def\problemname{Problem}
\def\problemsname{Problems}
\def\proofname{Proof}

\def\satzname{Satz} 
\def\koroname{Korollar}
\def\folgname{Folgerung}
\def\bemerkname{Bemerkung}
\def\aufgname{Aufgabe}

\def\beisname{Beispiel}
\def\beissname{Beispiele}
\def\bewname{Beweis}

\def\@thmcounter#1{\noexpand\arabic{#1}}
\def\@thmcountersep{}
\def\@begintheorem#1#2{\it \trivlist \item[\hskip 
\labelsep{\bf #1\ #2.\quad}]}
\def\@opargbegintheorem#1#2#3{\it \trivlist
      \item[\hskip \labelsep{\bf #1\ #2.\quad{\rm #3}}]}
\makeatother
\newtheorem{theor}{\theoremname}[section]
\newtheorem{propo}[theor]{\propositionname}
\newtheorem{coro}[theor]{\corollaryname}
\newtheorem{lemm}[theor]{\lemmaname}

\newenvironment{thm}{\begin{theor}\it}{\end{theor}}

\newenvironment{prop}{\begin{propo}\it}{\end{propo}}

\newenvironment{cor}{\begin{coro}\it}{\end{coro}}

\newenvironment{lem}{\begin{lemm}\it}{\end{lemm}}

\newtheorem{rema}[theor]{\remarkname}

\newenvironment{rem}{\begin{rema}\rm}{\end{rema}}

\newtheorem{stepnow}[theor]{}

\newtheorem{defin}[theor]{\definitionname} %% write

\newenvironment{defn}{\begin{defin}\rm}{\end{defin}}

\newtheorem{exerc}{\exercisename}[section]

\newtheorem{exa}[theor]{\examplename}

\newenvironment{ex}{\begin{exa}\rm}{\end{exa}}

\newtheorem{exas}[theor]{\examplesname}

\newenvironment{exs}{\begin{exas}\rm}{\end{exas}}

\newtheorem{conj}[theor]{\conjecturename}

\newtheorem{pro}[theor]{\problemname}

\newtheorem{prs}[theor]{\problemsname}

\newenvironment{Proof*}{\begin{trivlist}\item[\hskip%
\labelsep{\bf\proofname.\quad}]}% 
{\end{trivlist}}

\newenvironment{prf}{\begin{proof}}{\end{proof}}
  
{\hfill\qed\end{trivlist}}

\newenvironment{beweis*}{\begin{trivlist}\item[\hskip%
\labelsep{\bf\bewname.\quad}]}% 
{\end{trivlist}}

\newtheorem{satzn}[theor]{\satzname}

\newtheorem{koro}[theor]{\koroname}

\newtheorem{folg}[theor]{\folgname}

\newtheorem{bem}[theor]{\bemerkname}

\newtheorem{aufg}[theor]{\aufgname}

\newtheorem{aufgn}[theor]{\aufgname}

\newtheorem{beis}[theor]{\beisname}

\newtheorem{beiss}[theor]{\beissname}

\def\phi{\varphi}
\renewcommand{\mlabel}{\label}
\newcommand{\bx}{\mathbf{x}}
\newcommand{\by}{\mathbf{y}}
\newcommand{\bp}{\mathbf{p}}
\newcommand{\bt}{\mathbf{t}}
\newcommand{\bn}{\mathbf{n}}
\date{}
\begin{document} 

\title{Reflection positive one-parameter groups and dilations} 
\author{Karl-Hermann Neeb,
\begin{footnote}{
Department  Mathematik, FAU Erlangen-N\"urnberg, Cauerstrasse 11, 
91058-Erlangen, Germany; neeb@math.fau.de}
\end{footnote}
\begin{footnote}{Supported by DFG-grant NE 413/7-2, Schwerpunktprogramm 
``Darstellungstheorie''.} 
\end{footnote}
Gestur \'Olafsson
\begin{footnote}{Department of mathematics, Louisiana State University, 
Baton Rouge, LA 70803, USA; olafsson@math.lsu.edu}
\end{footnote}
\begin{footnote}
{The research of G. \'Olafsson was supported by NSF grants 
DMS-0801010, DMS-1101337 and the Emerging Fields Project 
``Quantum Geometry'' of the University of Erlangen.} 
\end{footnote}
}

\maketitle

\begin{abstract}The concept of reflection positivity has its origins in the work of Osterwalder--Schrader 
on constructive quantum field theory. It is a fundamental tool to
construct a relativistic quantum field theory as a unitary representation
of the Poincar\'e group from a non-relativistic field theory as a representation
of the euclidean motion group. This is the second article in a series on the mathematical
foundations of reflection positivity. We develop the theory of 
reflection positive one-parameter 
groups and the dual theory of dilations of contractive hermitian semigroups. 
In particular, we 
connect reflection positivity with the outgoing realization of unitary one-parameter
groups  by Lax and Phillips. We further show that our results provide effective tools 
to construct reflection positive representations of general symmetric 
Lie groups, including the $ax+b$-group, the Heisenberg group, the euclidean motion group and 
the euclidean conformal group. 
\end{abstract}

\section*{Introduction} 

This is the second article in a series of contributions to 
the mathematical foundations of \textit{reflection positivity}, a
basic concept in constructive quantum field theory 
(\cite{GJ81, JOl98, JOl00, JR08,JR07a, JR07, JP13}). 
In our first article \cite{NO12} we developed the the theory of reflection 
positive distributions and reflection positive distribution vectors. This approach 
is natural to obtain  classification results, especially 
in the abelian case.  In particular, we obtained a generalization of the Bochner--Schwartz Theorem 
to positive definite distributions on open convex cones and discussed the important example
of   reflection positive distribution vectors 
for complementary series representations 
of the conformal group $\OO_{1,n+1}^+(\R)$. 

The main objective of the present article is to develop the theory of reflection positive 
one-parameter groups and hermitian contractive 
semigroups as one key to understand reflection positivity for general 
symmetric Lie groups and their representations. One of our central methods is an 
$L^2$ realization of reflection positive one-parameter semigroups which is in some sense dual 
to the outgoing realization of one-parameter semigroups of isometries developed 
by Lax-Phillips \cite{LP64,LP67}. A more general
integration theory for reflection positive representations (of Lie algebras) 
will be developed in the sequel \cite{MNO14}.

The underlying fundamental concept is that of a 
{\it reflection positive Hilbert spaces}, introduced in Section \ref{sec:1}. 
This is a triple $(\cE,\cE_+,\theta)$, 
where $\cE$ is a Hilbert space, $\theta : \cE \to \cE$ is a unitary involution
and $\cE_+$ is a closed subspace of $\cE$ which is $\theta$-positive in the sense that 
the hermitian form $\langle u,v\rangle _\theta =\langle \theta u,v\rangle$ is
positive semidefinite on $\cE_+$. Let $\hat\cE$ denote the Hilbert space determined by 
this form and $q \: \cE_+ \to \hat\cE, v \mapsto \hat v$ be the canonical map. Then 
$\cE_0 = \{ v \in \cE_+ \: \theta v = v\}$ is the maximal subspace of $\cE_+$ on which 
$q$ is isometric. In Section~\ref{sec:1} we provide a collection of natural examples illustrating 
this concept. 

In Section \ref{sec:2} we introduce the notion of Osterwalder--Schrader quantization, as a 
passage from (densely defined) operators $T$ on $\hat\cE$ to operators $\hat T$ on $\hat\cE$ such 
that $\hat T \hat v = \hat{Tv}$, whenever this is well-defined. 
In the present context, we mostly consider operators obtained by 
restriction of unitary operators $U$ preserving $\cE_+$. If 
$U$ commutes with $\theta$, this leads to a unitary operator $\hat U$, and if 
$\theta U \theta = U^{-1}$, then $\hat U$ is a hermitian contraction. 
In particular, any unitary one-parameter group $(U_t)_{t \in \R}$ which is
{\it reflection positive} in the sense that 
$U_t(\cE_+)\subseteq \cE_+$ for $t \geq 0$ and $\theta U_t=U_{-t}\theta$ 
is quantized to a semigroup $(\hat U_t)_{t \geq 0}$ of hermitian contractions 
on $\hat\cE$. 

Section \ref{sec:3} is dedicated to a thorough inspection of reflection positive 
unitary one-parameter groups. Under the natural assumption 
that $\cE_+$ is cyclic, we derive some interesting consequences. 
One of the main results of Section~\ref{sec:3} is 
 that the subspace $\cE_{\rm fix}$ of $U$-fixed points in $\cE$ is contained 
in $\cE_0$ and maps onto the subspace $\hat\cE_{\rm fix}$ 
of fixed points of the corresponding semigroup 
$(\hat U_t)_{t > 0}$ of contractions on $\hat\cE$ (Proposition~\ref{prop:e.5}). 
We think of this fact as an instance of ``quantization commutes with reduction''. 
We also connect our theory with the scattering theory of
Lax-Phillips/Sinai \cite{LP64,LP67,Sin61}. In particular we show in 
Proposition \ref{prop:4.11}  that, if $\cE_+$ is $U$-cyclic in
$\cE$ and $\cE_{\rm fix}=\{0\}$, then
$\cE_+$ is an outgoing subspace of $\cE$ and hence $(U,\cE)$ is equivalent to the 
right translation group on $L^2(\R, \cM)$, where $\cM$ is a Hilbert space 
(representing the multiplicity). Unfortunately this realization does not 
provide a good picture for the involution $\theta$. The Lax--Phillips 
realization is unique up to unitary isomorphism of $\cM$. 
As a byproduct, it implies that 
the spectrum of $U$, resp., its infinitesimal generator, must be all of $\R$ if 
$\cE_+$ is cyclic. In view of the examples arising in physics, the 
case where the even smaller subspace $\cE_0 \subeq \cE_+$ is $U$-cyclic, is of crucial 
importance and sometimes even the much stronger condition 
$q(\cE_0) = \hat\cE$ is satisfied. A characterization of this case is provided 
in Proposition \ref{prop:3.9}. It implies that 
$\cE_+$ is generated by $U_t\cE_0$, $t > 0$, and that $\cE_+$ is maximal 
$\theta$-positive. A particularly interesting class of examples where 
$\cE_0$ is cyclic arises from the work of A.~Klein \cite{Kl77,Kl78} 
on Osterwalder--Schrader positive processes (cf.~Example~\ref{ex:klein}) 
and for the Hardy space of the real line (Subsection~\ref{ex:3.13}). 

In Section \ref{sec:dil} we change our perspective and start with a 
hermitian contraction semigroup $(C_t)_{t \geq 0}$ on a Hilbert space $\cH$. 
We call a reflection positive unitary one-parameter group $U$ on 
$(\cE,\cE_+,\theta)$ a {\it euclidean realization} of $(C,\cH)$ if 
$(\hat U, \hat\cE)$ is equivalent to $(C,\cH)$. Our starting point is the 
observation that euclidean realizations for which $\cE_0$ is cyclic 
always exist. This is derived from the positive  definiteness of the operator-valued 
function $\phi(t) = C_{|t|}$ on $\R$, which is a classical result in dilation theory of operators  
(see~\cite{SzNBK10}). To make the corresponding Hilbert space $\cE$ more concrete, 
we use the representation theory of positive definite functions by positive operator-valued 
measures developed in \cite{Ne98} to obtain a realization as a vector-valued 
$L^2$-space. An interesting case that is discussed in detail 
is the case where we have a \textit{time-zero} realization, i.e., where 
$q \:  \cE_0 \to \hat\cE\cong \cH$ is a unitary isomorphism. 
We also give a spectral characterization of the subspace $\cE_+$ as those elements in $\cE$ whose
inverse Fourier transform is supported in $[0,\infty [$. 
Here elements of $\cE_0$ correspond to functions whose Fourier transform is supported 
in $\{0\}$, hence to polynomial functions. 
In this picture, the basic building blocks are the spaces 
$\cE =L^2(\R\times \R_+,d\zeta)$, where 
$d\zeta = \frac{1}{\pi} \frac{\lambda}{x^2+\lambda^2}\, dx\, d\rho (\lambda)$ 
(which correspond to a cyclic contraction semigroup). 
We conclude Section~\ref{sec:dil} with an explicit outgoing realization of the 
$L^2$-space $\cE$ in the sense of Lax--Phillips theory (Proposition \ref{prop:4.9}). 

In Section \ref{sec:6} we apply the theory of reflection 
positive one-parameter groups to construct 
reflection positive representations of more general symmetric Lie groups, i.e., pairs 
$(G,\tau)$ consisting of a Lie group $G$ and an involutive automorphism $\tau$. This leads 
to a natural concept of a ``euclidean realization'' of a unitary representation of the 
dual symmetric Lie group $G^c$ whose Lie algebra is $\g^c = \fh + i \fq$,
where $\fh = \{ x \in \g \: \tau(x) = x\}$ and $\fq = \{ x \in \g \: \tau(x) = -x\}$. 
In all these examples, the realization theory based on dilations developed in 
Section~\ref{sec:dil} turns out to be amazingly effective to obtain the 
correct Hilbert spaces $\cE$, so that the main point is to implement the representation 
of the group $G$. We show that all representations of the 
$ax+b$-group and the Heisenberg group 
which satisfy the obvious necessary spectral condition do indeed have euclidean realizations. 
Furthermore, the dilation process immediately provides representations 
of the euclidean motion group of $\R^d$ associated to generalized free fields 
(cf.\ \cite{KL82}) and we show that the conformally invariant 
among these representations form the complementary series 
representations of the conformal group 
$\OO_{1,d+1}(\R)_+$ of $\R^d$, resp., its conformal completion 
$\bS^d$ (see \cite{LM75} for a discussion of 
conformal covariance in quantum field theories). 
This observation builds a bridge to the prequel \cite{NO12} 
where these representations were studied in some detail. 

In the appendices we collect and prove several
technical results, some of which are contained in one or the other form in 
the physical literature (\cite{GJ81}, \cite{RS75}). 
For the convenience of the reader, we include them with precise statements and proofs.

\tableofcontents 

\subsection*{Notation and terminology} 

Throughout this note, the space $C^{-\infty}(M)$ of distributions on a manifold 
$M$ consists of antilinear functionals on the space $C^\infty_c(M)$ of smooth 
compactly supported functions. 

We write elements $x \in \R^d$ as pairs $x = (x_0, \bx)$ with 
$x_0 \in \R$ and $\bx \in \R^{d-1}$. For the euclidean scalar product we 
accordingly write 
\[ \la x,y \ra = xy = x_0 y_0  +  \bx \by, \qquad 
x^2 = \|x\|^2 = x_0^2 + \bx^2\]
where $\bx \by=\la \bx , \by \ra$, 
and for the Lorentzian scalar product 
\[ [x,y] =  x_0 y_0  -  \bx \by, \qquad [x,x] = x_0^2 - \bx^2.\] 
The subset 
\[ V_+ := \{ x \in \R^d \: x_0 > 0, [x,x] > 0\} \] 
is called the open {\it forward light cone}. It is invariant under the action 
of the {\it orthochronous Lorentz group}
\[ L^\uparrow := \{ g \in \OO_{1,n+1}(\R) \: g_{00} > 0\}. \]

For the {\it Fourier transform of a measure $\mu$} on the
dual $V^*$ of a finite-dimensional real vector space $V$, we write
\[ \hat\mu(x) := \int_{V^*}e^{-i\alpha(x)}\, d\mu(\alpha).\]
The {\it Fourier transform of an $L^1$-function $f$} on $\R^d$ is defined by
\begin{equation}
  \label{eq:ftrafo}
\hat f(\xi) := (2\pi)^{-d/2} \int_{\R^d} f(x) e^{-i\la \xi,x\ra}\,  dx 
\end{equation}
and we likewise define convolution on $\R^d$ in terms of the Haar measure 
$\frac{dx}{(2\pi)^{d/2}}$. 

For tempered distributions $D \in \cS'(\R^d)$, we define the Fourier transform by 
\begin{equation}
  \label{eq:2}
\hat D(\phi) := D(\tilde \phi), \quad \mbox{ where} 
\quad \tilde \phi(\xi) = \hat\phi(-\xi) 
= (2\pi)^{-d/2} \int_{\R^d} \phi(x) e^{i\la \xi,x\ra}\,  dx.
\end{equation}
For $D(\phi) = \int_{\R^d} \oline{\phi(x)} f(x)\, dx$, we then obtain 
$\hat D(\phi) = \int_{\R^d} \oline{\phi(x)} \hat f(x)\, dx,$ 
so that $\hat D$ is represented by the function $\hat f$.

If $G$ is a locally compact topological group,
then $\Delta =\Delta_G: G\to [0,\infty[$ denotes the modular function, $\varphi^\vee(g)=\overline{\varphi (g^{-1})}$ and $\varphi^* 
=\varphi^\vee \cdot\Delta^{-1}$. If $\tau : G\to G$ is an involution, allowed to be the identity map, and
 $S\subset G$ is a semigroup invariant under $s\mapsto s^\ast=\tau(s)^{-1}$ then 
$\varphi^\sharp=\varphi^*\circ \tau$, or
 $\varphi^\sharp (g)=\overline{\varphi (\tau (g)^{-1})}\Delta (g^{-1})$.

\section{Reflection positive Hilbert spaces} 
\mlabel{sec:1} 
In this section we introduce the notion of a 
\textit{reflection positive Hilbert space}, a
concept that will play a fundamental role in the rest of this article.  
We also discuss various typical classes of examples. 

\begin{defn} \mlabel{def:x.1} Let $\cE$ be a Hilbert space and 
$\theta \in \U(\cE)$ be an involution. 
We call a closed subspace $\cE_+ \subeq \cE$ {\it $\theta$-positive} 
if $\la \theta v,v \ra \geq 0$ for $v \in \cE_+$. 
We then say that  the triple $(\cE,\cE_+,\theta)$ is a {\it reflection positive 
Hilbert space}. In this case we write 
\[\cN 
:= \{ v \in \cE_+ \: \la \theta v, v \ra = 0\} 
= \{ v \in \cE_+ \: (\forall w \in \cE_+)\ \la \theta w, v \ra = 0\} 
= \cE_+ \cap \theta(\cE_+)^\bot, \] 
$q \: \cE_+ \to \cE_+/\cN, v \mapsto \hat v = q(v)$ for the quotient map 
and $\hat\cE$ for the Hilbert completion of $\cE_+/\cN$ with respect to 
the norm $\|\hat v\|_{\hat{\cE}}:=\|\hat v\| := \sqrt{\la \theta v, v \ra}$. 
% not needed? For $v, w \in \cE$, we also write $\la v,w\ra_\theta = \la \theta v ,w\ra$. 

We write $\cE_0$ for the maximal subspace of $\cE_+$ 
on which the contraction $q \: \cE_+ \to \hat\cE$ is 
isometric. Writing $v \in \cE_+$ as $v = v_+ + v_-$ with 
$\theta v_\pm = \mp v_\pm$, we obtain 
\[ \la \theta v, v \ra = \|v_+\|^2 - \|v_-\|^2 
\leq \|v_+\|^2 + \|v_-\|^2 = \|v\|^2,\] 
and equality is equivalent to $v_- = 0$. This implies that 
\begin{equation}
  \label{eq:e_0}
 \cE_0 = (\cE_+)^\theta  = \{ v \in \cE_+ \: \theta(v) = v\} 
= \cE_+ \cap \theta(\cE_+),
\end{equation}
where the last equality follows from the fact that a $\theta$-invariant 
$\theta$-positive subspace of $\cE_+$ is contained in~$\cE^\theta$. 
\end{defn}

\begin{ex} \mlabel{ex:1.3} (a) Let $X$ be a set, $K \: X \times X \to \C$ be a positive definite 
kernel and $\cE = \cH_K \subeq \C^X$ the corresponding reproducing kernel 
Hilbert space. This is the unique Hilbert subspace of $\C^X$ on which all 
point evaluations $f \mapsto f(x)$ are continuous linear maps given by 
$f(x) = \la f, K_x\ra$ for $K_x(y) = K(y,x)$ (cf.\ Definition~\ref{def:1.5}). 
Suppose that $\tau \: X \to X$ is an involution 
leaving $K$ invariant, and $X_+ \subeq X$ a subset with the property that the 
kernel 
\[ K_+\: X_+ \times X_+ \to \C, \quad K_+(x,y) := K(\tau x, y) \] 
is positive definite. Then the closed subspace $\cE_+ \subeq \cE$ generated 
by the elements $K_x$, $x \in X_+$, is $\theta$-positive for 
$(\theta f)(x) := f(\tau x)$. We call $K$ {\it reflection positive} 
with respect to $(X,X_+, \tau)$. 

(b) Typical examples arise if $\tau$ is an involution on a group $G$ and $S \subeq G$ a subsemigroup 
invariant under $s \mapsto s^\sharp := \tau(s)^{-1}$. Then a function 
$\phi \: G \to \C$ is called {\it reflection positive} if the kernel 
$K(x,y) := \phi(xy^{-1})$ is reflection positive with respect to $(G,S,\tau)$ 
in the sense of (a) (see \cite{NO12} for a discussion of this concept). 
Prototypical examples are the functions 
$\phi(t) = e^{-\lambda|t|}$, $\lambda \geq 0$, for 
$(\R,\R_+, -\id_\R)$. 

(c) There are other important examples not related to subsemigroups. 
For  $\beta > 0$, the corresponding circle group 
$G := \R/\beta \Z$ and the domain 
$G_+ := \big[0,\frac{\beta}{2}\big] + \beta \Z\subeq G$, the functions 
$\phi \: G \to \C$ correspond to $\beta$-periodic functions on $\R$. 
In \cite{KL81} such a function is called {\it (OS)-positive} if the kernel 
$K(x,y) := \phi(x-y)$ is reflection positive for $(G,G_+, \tau)$ and 
$\tau(g) = g^{-1}$ in the sense of (a). 
In particular, integral representations of such functions 
are obtained. Typical examples are the $\beta$-periodic functions whose 
restriction to $[0,\beta]$ is given by 
$f_a(t) := e^{-ta} + e^{-(\beta - t)a}$ for $a \geq 0$. 
\end{ex}

\begin{ex} \mlabel{ex:1.4} Let $M$ be a smooth manifold and 
$D \in C^{-\infty}(M \times M)$ be a positive definite distribution. 
Suppose further that $\tau \: M \to M$ is an involutive diffeomorphism and 
that $M_+ \subeq M$ is an open subset such that the distribution $D_+$ 
on $M_+ \times M_+$ defined by 
\[ D_+(\phi) := D(\phi \circ (\tau \times \id_M)) \] 
is positive definite. We then say that $D$ is {\it reflection positive} with respect to 
$(M,M_+, \tau)$. 
Let $\cE$ denote the Hilbert space completion of $C^\infty_c(M)$, endowed with the 
scalar product 
\[ \la \phi, \psi \ra := D(\oline\phi \otimes \psi).\] 
Then the closed subspace $\cE_+$ generated by $C^\infty_c(M_+)$ is 
$\theta$-positive with respect to $\theta \phi := \phi \circ \tau$. 
\end{ex}

Before we turn to the next class of examples, we recall the concept of a 
distribution vector of a unitary representations. 
\begin{defn} (Distribution vectors)
Let $(\pi, {\cal H})$ be a continuous unitary
representation of the Lie group $G$ on the
Hilbert space ${\cal H}$. We write $\cH^\infty$ 
%(or $\cH^\infty(G)$ if several groups act on $\cH$) 
for the linear subspace of 
smooth vectors, i.e., of all elements $v \in \cH$ for which the orbit map
$\pi^v \: G \to \cH, g \mapsto \pi(g)v$ is smooth. Identifying
$\cH^\infty$ with the closed subspace of equivariant maps
in the Fr\'echet space $C^\infty(G,\cH)$, we obtain a natural Fr\'echet
space structure on $\cH^\infty$ for which the
$G$-action on this space is smooth
and the inclusion $\cH^\infty \to \cH$ (corresponding to evaluation in $\1 \in
G$) is a continuous linear map  (cf.\ \cite{Mag92}, \cite{Ne10}).

We write ${\cal H}^{-\infty}$ for the space of continuous antilinear
functionals on ${\cal H}^\infty$, the space of {\it distribution vectors},
and note that we have a natural linear
embedding $\cH \into \cH^{-\infty}, v \mapsto \la v, \cdot \ra$.
Accordingly, we also write $\la \alpha, v \ra 
= \oline{\la v, \alpha\ra}$ for $\alpha(v)$ for 
$\alpha \in \cH^{-\infty}$ and $v \in \cH^\infty$.
The group $G$ acts naturally on $\cH^{-\infty}$ by
\[ (\pi^{-\infty}(g)\alpha)(v) := \alpha(\pi(g)^{-1}v),\]
so that we obtain a $G$-equivariant chain of continuous inclusions
\begin{equation}\label{eq:rig}
  {\cal H}^\infty \subeq {\cal H} \subeq {\cal H}^{-\infty} 
\end{equation}
(cf.\ \cite[Sect.~8.2]{vD09}). It is $\cD(G)$-equivariant,
if we define the representation of $\cD(G)$ on $\cH^{-\infty}$ by
\begin{align*}
(\pi^{-\infty}(\phi)\alpha)(v)
&:= \int_G \phi(g) \alpha(\pi(g)^{-1}v)\, d\mu_G(g)
= \alpha(\pi(\phi^*)v).
\end{align*}
\end{defn}

\begin{ex} \mlabel{ex:1.2} If $\tau$ is an involution on a Lie group $G$ and $S \subeq G$ an open subsemigroup 
invariant under $s \mapsto s^\sharp := \tau(s)^{-1}$, then 
$C^\infty_c(S)$ is a $*$-algebra with respect to the convolution product 
and the $*$-operation $\phi \mapsto \phi^\sharp := \phi^* \circ \tau$. Accordingly, 
we call a distribution $D \in C^{-\infty}(S)$ {\it positive definite} if 
\[ D(\phi^\sharp * \phi) \geq 0 \quad \mbox{ for } \quad 
\phi \in C^\infty_c(S).\] 

We call a distribution 
 $D \in {\cal D}'(G)$  {\it reflection positive} with respect to 
$(G,\tau,S)$ if the following conditions are satisfied: 
\begin{description}
\item[\rm(RP1)] $D$ is positive definite, i.e., 
$D(\phi^* * \phi) \geq 0$ for $\phi \in C^\infty_c(G)$. 
\item[\rm(RP2)] $\tau D = D$, i.e., $D(\phi \circ \tau) = D(\phi)$ for 
$\phi \in C^\infty_c(G)$, and 
\item[\rm(RP3)] $D\res_S$ is positive definite 
as a distribution on the involutive semigroup $(S,\sharp)$, i.e., 
$D(\phi^\sharp * \phi) \geq 0$ for $\phi \in C^\infty_c(S)$. 
\end{description}

This is a special case of the situation in Example~\ref{ex:1.4}, 
where $G = M$, $S = M_+$, 
and the distribution $D^\sharp \in C^{-\infty}(G \times G)$ is defined by 
$D^\sharp(\phi \otimes \psi) := D(\phi^\vee * \psi)$, 
where $\phi^\vee(g) = \phi(g^{-1})\Delta_G(g)^{-1}$. 
The corresponding reproducing kernel Hilbert space 
$\cE = \cH_D \subeq C^{-\infty}(G)$ carries the unitary representation $\pi_D$ of $G$ whose 
integrated form is given by $\pi_D(\phi)E = \phi* E$ for $\phi \in C^\infty_c(G)$. 
Then $D \in \cE^{-\infty}$ is a distribution vector and 
the closed subspace $\cE_+$ generated by $\pi^{-\infty}(C^\infty_c(S))D = C^\infty_c(S) * D 
\subeq \cE$ is $\theta$-positive for $(\theta E)(\phi) := E(\phi \circ \tau)$. 
\end{ex}

\begin{ex} \mlabel{ex:1.5} 
(a) Let $M$ be a smooth manifold and 
$\mu$ a measure on $C^{-\infty}(M)$ with respect to the smallest $\sigma$-algebra 
for which all evaluation maps 
$\phi^*(D) := D(\phi)$, $\phi \in C^\infty_c(M)$, are measurable  
(cf.\ \cite{GV64}). Let $\cE := L^2(C^{-\infty}(M),\mu)$ be the corresponding 
$L^2$-space. 

For an open subset $M_+ \subeq M$, we consider the corresponding 
closed subspace $\cE_+ \subeq \cE$ generated by the functions of the form 
$e_\phi(D) := e^{iD(\phi)}$, $\phi \in C^\infty_c(M_+)$. Further, let $\tau \: M \to M$ 
be an involutive isomorphism whose action on $C^{-\infty}(M)$ preserves $\mu$, 
and write $(\theta F)(D) := F(\tau D)$ for the corresponding unitary 
involution on $\cE$. 

We then say that $\mu$ is {\it reflection positive with respect to 
$(M,M_+, \tau)$} if $(\cE, \cE_+, \theta)$ is reflection positive. 
This structure naturally occurs in euclidean QFT, where 
$M = \R^d$, $M_+ = \R^d_+$ and $\tau(x) = (-x_0, x_1,\ldots, x_{d-1})$ 
(cf.\ \cite{GJ81} and also \cite{JR08,JR07a, JR07, JP13}).

(b) To interpret all that in the context of Example~\ref{ex:1.3}, we assume that 
$\mu$ is finite. Then its Fourier transform 
\[ S(\phi) := \hat\mu(\phi) := \int_{C^{-\infty}(M)} e^{-iD(\phi)}\, d\mu(D)\] 
is a sequentially continuous positive definite function. 
Then $K(\phi, \psi) := S(\phi - \psi)$ is a positive definite kernel 
on $X = C^\infty_c(M)$ which is reflection positive with respect to 
$X_+ := C^\infty_c(M_+)$ and the involution $\tau \phi := \phi \circ \tau$. 
Here $K$ is non-linear, whereas 
it is sesquilinear for the kernels corresponding to distributions 
$D \in C^{-\infty}(M \times M)$ as in Example~\ref{ex:1.4}. 
\end{ex}

\begin{exs} \mlabel{ex:1.6} 
(a) Let  $(\cE,\cE_+, \theta)$ be a reflection positive Hilbert 
space. We write $\oline\cE$ for the complex Hilbert space obtained 
by endowing $\cE$ with the opposite complex structure and the conjugate 
scalar product. For 
$X := \oline\cE$, $X_+ := \oline{\cE_+}$, $\tau := \theta$, 
and the positive definite kernel 
$K(x,y) := \la y,x\ra$, we then obtain 
$\cH_K \cong \cE$ and $\cH_{K_+} \cong \cE_+$ (cf.\ \cite[Ex.~I.1.10]{Ne00}). 
Therefore every reflection 
positive Hilbert space can be obtained in the context of Example~\ref{ex:1.3}. 

In Example~\ref{ex:1.4} we can put 
\[ X := C^\infty_c(M), \quad X_+ := C^\infty_c(M_+) \quad \mbox{ and } \quad 
K(\psi, \phi) = D(\oline\phi \otimes \psi) \]  
to put the construction of the associated reflection positive Hilbert 
space into the context of Example~\ref{ex:1.3}. 

(b) If $D\in C^{-\infty}(M \times M)$ is reflection positive with respect to 
$(M,M_+,\tau)$, then 
\[S(\phi) := e^{-\shalf D(\oline \phi \otimes \phi)}\]
is a 
positive definite function on the nuclear space $C^\infty_c(M)$  and the corresponding 
reproducing kernel Hilbert space is isomorphic to the 
Fock space $\cF(\cH_D)$ (cf.\ Remark~\ref{rem:1.7} below). 
In view of the Bochner--Minlos Theorem (\cite{GV64}), there exists a unique Gaussian measure 
$\mu$ on $C^{-\infty}(M)$ with $\hat \mu = S$, which implies that 
$\cF(\cH_D) \cong L^2(C^{-\infty}(M), \mu)$. It is easy to see that 
the reflection positivity of $D$ with respect to $(M,M_+, \tau)$ is equivalent 
to the reflection positivity of $\mu$ with respect to $(M,M_+, \tau)$ 
(cf.~Remark~\ref{rem:1.7} and Lemma~\ref{lem:exp-crit} below). 
\end{exs}

\begin{rem} \mlabel{rem:1.7} (Fock spaces and reflection positivity) 
We have seen in Example~\ref{ex:1.6}(a), that we can consider a Hilbert space 
$\cE$ as a reproducing kernel space of antilinear functions on $\cE$, defined by the kernel 
$(v,w) \mapsto \la w,v\ra$. Accordingly, we may consider the Fock space 
$\cF(\cE)$ as a reproducing kernel space on 
$\cE$ with kernel $K(v,w) := e^{\la w, v\ra}$ (this is a realization by antiholomorphic functions). 
Then each operator $U \in \U(\cE)$ induces a unitary operator 
$\hat U$ in $\cF(\cE)$ by $(\hat U F)(v) := F(U^{-1}v)$, which leads 
to the relation $\hat U K_w  = K_{Uw}$ for $w \in \cE.$

Suppose that $\cE_+ \subeq \cE$ is a $\theta$-positive subspace. 
Let $\cF(\cE)_+ \subeq \cF(\cE)$ denote the closed subspace 
generated by the functions $K_w$, $w \in \cE_+$. 
Then the kernel 
\[ \la \hat\theta K_v, K_w \ra 
=  \la K_{\theta v}, K_{w} \ra = K(w, \theta v) = e^{\la \theta v, w \ra} \] 
is positive definite on $\cE_+$. This implies that 
$\cF(\cE)_+$ is $\hat\theta$-positive and that the corresponding 
Hilbert space can be identified with the Fock space 
$\cF(\hat\cE)$, where $\hat\cE$ is the ``quantum Hilbert space'' associated
 to $(\cE, \cE_+, \theta)$. 
If $U \in \U(\cE)$ preserves the subspace $\cE_+$, 
then $\hat U K_w = K_{U w}$ implies that $\cF(\cE)_+$ is 
invariant under $\hat U$. 
%This discussion shows that the passage 
%to the Fock space respects reflection positive subspaces. 
%for unitary one-parameter groups and, first of all, 
%that the ``functor $\cF$'' commutes with OS-quantization, i.e., 
%the passage from $(\cE, \cE_+, \theta)$ to $\hat\cE$. 
\end{rem}

\section{OS-quantization}
\mlabel{sec:2} 
In this section we introduce Osterwalder--Schrader quantization as a 
method to pass from operators on a reflection positive Hilbert space 
$(\cE,\cE_+,\theta)$ to operators on $\hat\cE$.  
Our terminology follows \cite[\S VII.7]{Ja08}. 
If $T$ is a linear operator from a Hilbert space $\cH$ to a Hilbert space 
$\cK$ which is not everywhere defined, we write 
$\cD(T) \subeq \cH$ for its domain of definition. 
%If $T$ is continuous then $\cD (T)=\cH$.

\begin{defn} \mlabel{def:e.1-os} (OS-quantization)
Let $(\cE,\cE_+,\theta)$ be a reflection positive Hilbert space. 
Suppose that $T\: \cD(T) \to \cE_+$ is an operator on $\cE_+$, possibly unbounded,
with  $T(\cD(T) \cap \cN) \subeq \cN.$ 
Then $T$ induces a linear operator 
\[ \hat T \: \cD(\hat T) := \{ \hat v \: v \in \cD(T) \} \to \hat\cE, \quad 
\hat T \hat v := \hat{Tv}.\] 
The passage from $T$ to $\hat T$ is called {\it Osterwalder--Schrader 
quantization} (or $OS$-quantization for short).
\end{defn}

\begin{lem}
  \mlabel{lem:d.1} 
Let $(\cE,\cE_+,\theta)$ be a reflection positive Hilbert space 
and $\cD \subeq \cE_+$ be a subspace whose image $\hat\cD$ in $\hat\cE$ is dense. 
Suppose that $T, U \: \cD \to \cE_+$ are linear maps satisfying 
\[ \la \theta Tv, w \ra = \la \theta v, U w\ra \quad \mbox{ for } \quad v,w \in \cD.\] 
Then the following assertions hold: 
\begin{description}
\item[\rm(i)]  $\hat T \hat v := \hat{Tv}$ and $\hat U\hat v := \hat{Uv}$ for 
$v \in \cE_+$ define linear operators with domain $\hat\cD \subeq \hat\cE$. In 
particular, both map $\cN \cap \cD$ into $\cN$. 
\item[\rm(ii)] $\hat U \subeq \hat T^*$. 
\item[\rm(iii)] If $U = T$, then $\hat T$ is symmetric. 
\item[\rm(iv)] If $T= U$ is bounded and $\cD = \cE_+$, 
then  $\hat T$ is bounded with $\|\hat T\| \leq \|T\|$. 
\end{description}
\end{lem}

\begin{prf} (i) For $v,w \in \cD$ we have 
\begin{equation}
  \label{eq:symrel}
\la \hat{Tv}, \hat w \ra = \la \theta T v, w \ra 
= \la \theta v, U w \ra = \la \hat v, \hat{Uw} \ra.
\end{equation}
It follows in particular, that $\hat v = 0$ implies that 
$\hat{Tv} \in \hat \cD^\bot = \{0\}$. Therefore $T$ induces by 
$\hat T \hat v := \hat{Tv}$ a well-defined operator 
on $\hat\cD$, and the same argument applies to $U$. 

(ii) and (iii) are immediate consequences of \eqref{eq:symrel}. 

We now prove (v) (cf.\ \cite{JOl00}). Suppose that $T$ is bounded and $\cD = \cE_+$. 
Then $\hat T\hat\cD \subeq \hat\cD$ implies that all powers $\hat T^n$ define 
symmetric operators on the pre-Hilbert space $\hat\cD$. First, 
\[ \|\hat T \hat v\|^2 \leq \|\hat T^2 \hat v\| \|\hat v\| 
\quad \mbox{ for } \quad v \in \cE_+\] 
follows from the Cauchy--Schwarz inequality. Iterating this argument leads to 
\[ \|\hat T \hat v\|^{2^n} 
\leq \|\hat T^2 \hat v\|^{2^{n-1}} \|\hat v\|^{2^{n-1}} 
\leq \|\hat T^4 \hat v\|^{2^{n-2}} \|\hat v\|^{2^{n-1} + 2^{n-2}} 
\leq \cdots \leq \|\hat T^{2^n} \hat v\| \|\hat v\|^{2^n-1}.\] 
We also have 
$\|\hat T^m \hat v\|^2 \leq \|T\|^{2m} \|v\|^2,$ 
which leads to 
\[ \|\hat T \hat v\|^{2^n} 
\leq \|T\|^{2^n} \|v\| \|\hat v\|^{2^n-1}.\] 
We conclude that 
\[ \|\hat T \hat v\| \leq \|T\| \limsup_n  \|v\|^{2^{-n}} \|\hat v\|^{1 - 2^{-n}} 
= \|T\| \|\hat v\|.\] 
Therefore $\hat T$ is bounded with 
$\|\hat T\| \leq \|T\|$. 
\end{prf}

\begin{ex} \mlabel{ex:d.2} Typical examples to which the preceding lemma applies are 
unitary operators $T$ with $\theta T \theta = T^{-1}$ and $T\cE_+ \subeq \cE_+$, and 
hermitian operators $T$ commuting with $\theta$. In both cases
$\hat T : \hat E\to \hat E$ is continuous and $\hat T^*=\hat T$.  
\end{ex}

Here is a simple observation about unitary operators mapping $\cE_+$ onto itself. 
\begin{lem} \mlabel{lem:2.5} 
Let $(\cE,\cE_+,\theta)$ be reflection positive and $U :\cE\to \cE$ be unitary 
with $U\cE_+ = \cE_+$. Then the following assertions hold: 
\begin{description}
\item[\rm(a)] If $\theta U \theta = U$, then $\hat U :\hat \cE\to \hat\cE$ is unitary. 
\item[\rm(b)] If $\theta U \theta = U^{-1}$, then $\hat U^2 = \id_{\hat\cE}$. 
Moreover, $\cE$ is a direct sum of reflection positive Hilbert subspaces 
$(\cF,\cF \cap \cE_+, \theta\res_{\cF})$ and $(\cG,\cG \cap \cE_+, \theta\res_{\cG})$, 
invariant under $U$ and $U^{-1}$, such that 
$\hat\cG = \{0\}$ and $(U\res_{\cF})^2 = \1$. 
\end{description}
\end{lem}

\begin{prf} (a) is immediate from the definitions. 

(b) Lemma~\ref{lem:d.1}, applied to $U$, implies that 
$\hat U$ is a symmetric contraction. Applying the same lemma to 
$V := U^{-1}$ leads to another symmetric contraction $\hat V$. 
Now $\hat U \hat V \hat v = \hat{UV}\hat v = \hat v$ for every $v \in \cE_+$ 
implies that $\hat U \hat V = \id_{\hat\cE}$. We likewise 
get $\hat V \hat U = \id_{\hat\cE}$, so that $\hat U^{-1} = \hat V$, and 
thus $\Spec(\hat U) \subeq \{-1,1\}$ leads to $\hat U^2 = \id_{\hat\cE}$.

Since the closed subspace $\cN = \cE_+ \cap \theta(\cE_+)^\bot$ of $\cE_+$ is also 
invariant under $U$ and $V= U^{-1}$ (Lemma~\ref{lem:d.1}(i)), the subspace 
$\cG := \oline{\cN \oplus \theta(\cN)^\bot}$ is invariant under $U, U^{-1}$ and $\theta$. Its orthogonal 
complement $\cF$ contains $\cE_+ \cap \cN^\bot$, so that $\hat \cG = \{0\}$ and 
$\hat\cF = \hat\cE$. In particular, $q\res_{\cF_+}$ is injective. 
Now $q \circ U\res_{\cE_+} = \hat U \circ q$ implies that $U_+ := U\res_{\cF_+}$ 
also satisfies the relation $U_+^2 = \1$. Since $\cF_+ + \theta(\cF_+)$ is invariant 
under $U,U^{-1}$ and $\theta$, it is dense in $\cF$, and this leads to 
$(U\res_{\cF})^2 = \1$.
\end{prf}

As a consequence we get:

\begin{prop} \mlabel{prop:osquant-commute} 
Let $(\cE,\cE_+,\theta)$ be a reflection positive Hilbert space. 
Suppose that $(U,\cE)$ is a strongly continuous unitary representation 
of a topological group $G$ on $\cE$ such that 
\[ U_g \cE_+ \subeq \cE_+ \quad \mbox{ and } \quad U_g \theta = \theta U_g 
\quad \mbox{ for } \quad g \in G.\] 
Then OS-quantization defines a continuous 
unitary representation $(\hat U, \hat\cE)$ of $G$ on $\hat\cE$. 
\end{prop}

\begin{rem} Suppose that $\cE$ is finite-dimensional and that $U \in \U(\cE)$ 
satisfies $U\cE_+ \subeq \cE_+$ and 
$\theta U \theta = U^{-1}$. Then the finite dimension 
implies that $U \cE_+ = \cE_+$, so that Lemma~\ref{lem:2.5}(b)  
show that $\hat U^2 = \1$. 
\end{rem}

\section{Reflection positive unitary one-parameter groups} 
\mlabel{sec:3} 

In this section we define the concept of a reflection positive 
one-parameter group $(U_t)_{t \in \R}$ on a reflection positive 
Hilbert space $(\cE,\cE_+,\theta)$ by the requirements that 
$U_t \cE_+ \subeq \cE_+$ for $t>0$ and 
$\theta U_t \theta = U_{-t}$ for $t\in \R$. 
Under the natural assumption 
that $\cE_+$ is cyclic, we derive some interesting consequences. 
The main results of this section are
 that the subspace $\cE_{\rm fix}$ of $U$-fixed points is contained 
in $\cE_0$ and maps onto the space $\hat\cE_{\rm fix}$ 
of fixed points of the corresponding one-parameter semigroup 
$(\hat U_t)_{t > 0}$ of contractions on $\hat\cE$ 
(Proposition~\ref{prop:e.5}). If $\cE_{\rm fix}$ is trivial, 
we show that $\cE_+$ is outgoing in the sense of 
Lax--Phillips scattering theory (Proposition~\ref{prop:4.11}). 
In Subsection~\ref{subsec:3.3}, we take a closer look at the case 
where $\cE_0$ is cyclic. Here the case 
$q(\cE_0) = \hat\cE$ is of particular interest, 
and we show in Proposition~\ref{prop:3.9} that in this case 
$\cE_+$ is generated by $U_t\cE_0$, $t > 0$, and that it is maximal 
$\theta$-positive. A particularly interesting class of examples where 
$\cE_0$ is cyclic arises from the work of A. Klein \cite{Kl77,Kl78} 
on Osterwalder--Schrader positive processes 
(cf.~Example \ref{ex:klein}).

\subsection{The associated contraction semigroup} 
\mlabel{subsec:3.1}

\begin{defn} Let $(\cE, \cE_+,\theta)$ be a reflection positive Hilbert space. 
A {\it reflection positive} unitary one-parameter group on 
$(\cE,\cE_+,\theta)$ is a strongly continuous unitary one-parameter group 
on $\cE$ for which $\cE_+$ is invariant under $U_t$ for $t > 0$ and 
$\theta U_t \theta  = U_{-t}$ for $t \in \R$. 
\end{defn}

\begin{prop} \mlabel{prop:os-onepar} If $(U_t)_{t \in \R}$ is a unitary one-parameter group 
which is reflection positive with respect to 
$(\cE,\cE_+, \theta)$, 
then each $U_t$, $t > 0$, induces a hermitian contraction 
$\hat U_t$ on $\hat\cE$. If $U$ is strongly continuous, we thus obtain a 
strongly continuous one-parameter semigroup of contractions. 
\end{prop}

\begin{prf} From Example~\ref{ex:d.2} we obtain the hermitian contractions 
$\hat U_t$, $t \geq 0$. The defining relation $\hat U_t \hat v = \hat{U_t v}$ 
easily implies that the operators $(\hat U_t)_{t \geq 0}$ form a one-parameter semigroup. 
Suppose that $U$ is continuous. 
Then the weak continuity of 
$\hat U$ follows from the fact that, for $v \in \cE_+$,  
the function 
$t \mapsto \la \hat U_t \hat v, \hat v \ra 
= \la \theta U_tv, v  \ra$  
is continuous. Now \cite[Cor.~IV.1.18]{Ne00} implies that $\hat U$ is strongly continuous.
\end{prf}

\begin{defn} \mlabel{def:3.3} In the context of the preceding proposition, we call 
the quadruple $(\cE,\cE_+, \theta, U)$ a {\it euclidean realization} 
of the contraction semigroup $(\hat U, \hat\cE)$.  
\end{defn}

We shall see in Proposition~\ref{prop:exist-real} below 
that every strongly continuous contraction semigroup on a Hilbert space has a 
euclidean realization, but in general there are many non-equivalent 
realizations, as the following example shows. 

\begin{ex} (see Section ~\ref{ex:4.5} and example \ref{ex:4.11} for details). 
For a finite measure $\rho$ on $\R_+$, consider 
on $\cH := L^2(\R_+,\rho)$ the contraction semigroup 
$(C_t f)(x) = e^{-tx} f(x)$. Then a euclidean realization 
of infinite multiplicity is given by 
\[ \cE = L^2(\R \times \R_+, \zeta), \quad 
d\zeta(x,\lambda ) = \frac{1}{\pi} \frac{\lambda }{\lambda^2 + x^2}\, dx\, d\rho(\lambda ), \quad 
(U_t f)(x,\lambda ) = e^{-itx} f(x,\lambda )\] 
(see Section~\ref{sec:dil}). 
But there are also multiplicity free realizations, such as 
\[ \cE = L^2(\R,\nu), \quad 
d\nu(x) = \frac{1}{\pi}\int_0^\infty \frac{\lambda }{\lambda^2 + x^2}\, d\rho(\lambda )\, dx, \quad 
(U_t f)(x) = e^{-itx} f(x).\] 

In both cases $\cE_0$ is cyclic, but in the first case $q(\cE_0) = \hat\cE \cong \cH$, 
whereas in the second case $\cE_0 = \C 1$ is one-dimensional. 
\end{ex}

\begin{lem}\mlabel{lem:eplus-cyc}
Let $(U_t)_{t \in \R}$ be a unitary one-parameter group on 
$\cE$ and $\cE_+$ be a subspace invariant under $U_t$ for $t > 0$. 
Then the following  conditions are equivalent: 
\begin{description}
  \item[\rm(i)] The subspaces $U_t \cE_+$, $t < 0$, span a dense subspace of $\cE$. 
  \item[\rm(ii)] There exists a dense subspace $\cD \subeq \cE$ such that 
$U_t v \in \cE_+$ for $v \in \cD$ and $t$ sufficiently large. 
\item[\rm(iii)] $\cE_+$ is $U$-cyclic in $\cE$. 
\end{description}
\end{lem}

\begin{prf} Since $U_t \cE_+ \subeq U_s \cE_+$ for $t > s$, the subset 
$\cE_{-\infty} := \bigcup_{t \in \R} U_t \cE_+= \bigcup_{t < 0} U_t \cE_+$ 
is a linear subspace. Assertions (i) and (iii) mean 
that $\cE_{-\infty}$ is dense and (ii) means that $\cE_{-\infty}$ contains a dense subspace. 
Therefore (i), (ii) and (iii) are equivalent.
\end{prf}

The following lemma provides a criterion for the 
density of a subspace of $\hat\cE$. We shall use it 
to verify that certain operators on $\hat \cE$ are densely defined. 
\begin{lem} \mlabel{lem:e.6} 
Let $(U_t)_{t \in \R}$ be a reflection positive unitary one-parameter group 
on $(\cE,\cE_+, \theta)$. 
If $\cD \subeq \cE_+$ is a subspace invariant under the operators 
$U_t$, $t > 0$, for which 
\[ \cE_+^0 := \{ v \in \cE_+ \: (\exists T > 0)\ U_T v \in \cD \} \] 
is dense in $\cE_+$, then $\hat\cD \subeq \hat\cE$ is dense. 
\end{lem} 

\begin{prf} For $w \in \cE_+^0$ there exists a $T > 0$ with 
$U_T w \in \cD$, and this implies that 
$\hat U_t \hat w \in \hat\cD$ for $t \geq T$. Since the curve 
$\R_+ \to \hat\cE, t \mapsto \hat U_t w$, is analytic, 
$\hat U_t w \in \oline{\hat\cD}$ for every $t > 0$, and therefore 
$w \in \oline{\hat\cD}$ follows from the strong continuity of the semigroup 
$(\hat U_t)_{t \geq 0}$ (Proposition~\ref{prop:os-onepar}). As $\cE_+^0$ is dense in $\cE_+$ it
follows that $\widehat{\cD}$ is dense in $\widehat{\cE}$.
\end{prf}

\begin{rem} (Reduction to the $\cE_0$-cyclic case if 
$\hat{\cE_0}$ is cyclic in $\hat\cE$) 
Assume that $(U_t)_{t \in \R}$ is 
reflection positive on $(\cE,\cE_+,\theta)$ 
and that the image $q(\cE_0)$ of $\cE_0 = \cE_+^\theta$ in 
$\hat\cE$ is $\hat U$-cyclic. 

Let $\cE' \subeq \cE$ denote the 
closed $U$-invariant subspace generated by $\cE_0$ and 
$\cE'_+ := \cE' \cap \cE_+$. Then 
$\theta U_t \cE_0 = U_{-t} \theta \cE_0 = U_{-t} \cE_0$ implies that 
$\cE'$ is $\theta$-invariant. Therefore $U_t' := U_t\res_{\cE'}$ 
is a reflection positive unitary one-parameter group 
on $(\cE',\cE'_+, \theta\res_{\cE'})$. Since 
$q\res_{\cE_+'}$ has dense range, all the relevant data 
is contained in $\cE'$. It is therefore natural to assume that 
$\cE_0$ is $U$-cyclic in $\cE$ whenever $q(\cE_0) = \hat{\cE_0}$ is cyclic in $\hat\cE$. 
\end{rem} 

In quantization procedures it is of particular importance to which 
extent they are compatible with reduction. The following proposition 
is an instance of ``OS-quantization commutes with reduction''.

\begin{prop} \mlabel{prop:e.5} 
{\rm(OS-quantization commutes with reduction)}
Let $(U_t)_{t \in \R}$ be a reflection positive unitary one-parameter group 
on $(\cE,\cE_+, \theta)$. 
Suppose that $\cE_+$ is $U$-cyclic and 
write $(\hat U_t)_{t \geq 0}$ for the one-parameter semigroup of contractions 
induced by $U$ on $\hat\cE$. 

Let $\cE_{\rm fix}$ denote the subspace of elements fixed under 
all $U_t$ and $\hat\cE_{\rm fix}$ the subspace of fixed points 
for the semigroup $(\hat U_t)_{t > 0}$. Then the following assertions hold: 
\begin{description}
  \item[\rm(a)] $\cE_{\rm fix} \subeq \cE_0$, the space of $\theta$-fixed points 
in $\cE_+$. 
  \item[\rm(b)] The map $q\res_{\cE_{\rm fix}} 
\: \cE_{\rm fix} \to \hat\cE_{\rm fix}, v \mapsto \hat v$ 
is a unitary isomorphism. 
  \item[\rm(c)] $\cE_{\rm fix} = \cE_\infty := \bigcap_{t > 0} U_t\cE_+$.
\end{description}
\end{prop}

\begin{prf} (a) We write $P \: \cE \to \cE_{\rm fix}$ for the 
orthogonal projection onto the subspace of $U$-fixed points in $\cE$. 
Then 
\[ \lim_{N \to \infty} \frac{1}{N} \int_0^N U_t\, dt = P  \] 
holds in the strong operator topology (\cite[Cor.~V.4.6]{EN00}).
Let $\cD$ be as in Lemma~\ref{lem:eplus-cyc}(ii). 
For $v \in \cD$, there exists a $T > 0$ with 
$U_t v \in \cE_+$ for $t > T$. 
Since 
\[ \lim_{N \to \infty} \frac{1}{N} \int_0^T U_t\, dt = 0, \] 
we obtain 
$P v\in \cE_+$ for every $v \in \cD$. The density of $\cD$ in $\cE$ 
thus implies that 
$\cE_{\rm fix} = P\cE \subeq \cE_+.$ 
Since the fixed point spaces of $(U_t)_{t > 0}$ and $(U_{-t})_{t > 0}$ coincide, 
we also have 
\[ \theta P \theta = \lim_{N \to \infty} \frac{1}{N} \int_0^N U_{-t}\, dt = P,\] 
so that $\cE_{\rm fix}$ is $\theta$-invariant. 
Now the $\theta$-positivity of $\cE_+$ implies that 
$\theta\res_{\cE_{\rm fix}} \geq 0$, and thus $\cE_{\rm fix} \subeq \cE^\theta$. 

(b) Since $P$ commutes with $\theta$, 
Lemma~\ref{lem:d.1}(iv) shows that $P$ defines a hermitian contraction 
$\hat P \: \hat\cE \to \hat\cE$ with $\hat P \hat v = \hat{Pv}$ for $v \in \cE_+$. 
As $P^2 = P$ entails $\hat P^2 = \hat P$, $\hat P$ is a projection. 
For $v,w \in \cE_+$, we obtain 
\[ \lim_{N \to \infty} \frac{1}{N} \int_0^N \la \hat v, \hat U_t \hat w \ra\, dt 
=  \lim_{N \to \infty} \frac{1}{N} \int_0^N \la \theta v, U_t w \ra\, dt 
= \la \theta v, P w \ra = \la \hat v, \hat P \hat w \ra.\] 
Hence \cite[Cor.~V.4.6]{EN00} implies 
that $\hat P$ is the orthogonal projection onto $\hat\cE_{\rm fix}$.

Let $q \: \cE_+ \to \hat\cE, v \mapsto \hat v$, denote the canonical projection onto $\widehat\cE$. Then 
$q \circ P = \hat P \circ q$ implies that  
$q(\cE_{\rm fix}) = q(P\cE_+) = \hat P q(\cE_+)$, and hence that 
$q(\cE_{\rm fix}) \subeq \hat\cE_{\rm fix}$ is a dense subspace. 
On the other hand, $\cE_{\rm fix} \subeq \cE_+^\theta$ implies that 
$q\res_{\cE_{\rm fix}}$ is isometric, hence a unitary isomorphism onto 
$\hat\cE_{\rm fix}$. 

(c) The subspace $\cE_\infty$ is closed and it is easily seen to be invariant 
under $U$.  Therefore 
$\cF := \oline{\cE_\infty + \theta\cE_\infty}$ is invariant under 
$U$ and $\theta$, so that we obtain a reflection 
positive unitary one-parameter group 
$V_t := U_t\res_{\cF}$ on $(\cF,\cF_+ := \cE_\infty, \theta\res_{\cF})$ 
satisfying $V_t \cF_+ = \cF_+$ for every $t > 0$. 
Now Lemma~\ref{lem:2.5}(b) leads to $\hat V_t = \hat V_{t/2} \hat V_{t/2} = \1$ 
for every $t > 0$. Therefore $\hat\cF \subeq \hat \cE_{\rm fix}$, 
and (b) implies that $\hat\cE_{\rm fix} =q(\cE_{\rm fix})$, so that 
$\cE_\infty = \cF_+ \subeq \cE_{\rm fix} + \cN$. 

Since the elements of $\cE_{\rm fix}$ are $\theta$-fixed 
and $\cN = \cE_+ \cap \theta(\cE_+)^\bot$, we have $\cN \bot \cE_{\rm fix}$. 
>From $\cE_{\rm fix} \subeq \cE_\infty$ it thus follows that 
$\cE_\infty = \cE_{\rm fix} \oplus (\cN \cap \cE_\infty)$ 
is a $U$-invariant orthogonal decomposition. 
As $\cN \cap \cE_\infty$ is orthogonal to the $U$-cyclic subspace $\theta(\cE_+)$, 
it must be zero, and this shows that $\cE_\infty = \cE_{\rm fix}$. 
\end{prf}

\begin{rem} \mlabel{rem:3.8}
Let $\cE^1 := \cE_{\rm fix}^\bot$. 
Then all the structure of a reflection positive one-parameter group 
is adapted to the orthogonal decomposition $\cE = \cE_{\rm fix} \oplus \cE^1$: 
\[ \cE_+ = \cE_{\rm fix} \oplus \cE^1_+, \quad 
\theta = \1 \oplus \theta_1, \quad 
U_t = \1 \oplus U^1_t \] 
with respect to the obvious notation. The data corresponding to $\cE_{\rm fix}$ is 
trivial and the one-parameter group $(U^1_t)_{t \in \R}$ on $(\cE^1, \cE^1_+,\theta)$ 
has the additional property that $\cE^1_{\rm fix} = \{0\}$. We also have that 
$\hat\cE \cong \hat\cE_{\rm fix} \oplus \hat\cE_1$. 
\end{rem}

\begin{ex} Any unitary one-parameter group $(V_t)_{t \in \R}$ 
on a Hilbert space $\cH$ can be embedded into a $\theta$-positive one 
in a trivial manner. We simply put 
$\cE := \cH \oplus \cH$, $\cE_+ := \cH \oplus \{0\}$, 
$U_t := V_t \oplus V_{-t}$ and $\theta(v,w)  = (w,v)$. 
Then $\hat\cE = \{0\}$ and 
$\cE_\infty = \cE_+$ is not cyclic. 
\end{ex}

\begin{lem} \mlabel{lem:3.11} 
If $(U_t)_{t \in \R}$ is a reflection positive one-parameter group 
on $(\cE,\cE_+,\theta)$ with a cyclic $\theta$-invariant distribution vector, 
then $\cE_+$ is cyclic. 
\end{lem}

\begin{prf} Let $\alpha \in \cE^{-\infty}$ be a reflection positive 
cyclic distribution vector. 
Then $\cE_+$ is the closed subspace generated by 
$U^{-\infty}(\phi)\alpha$, $\phi \in C^\infty_c(\R_+)$. 
For $\phi \in C^\infty_c(\R)$, there exists a $t > 0$ with 
$\supp(\phi) + t \subeq \R_+$. Then 
\[ U^{-\infty}(\phi)\alpha 
= U^{-\infty}(\delta_{-t} * \delta_t * \phi)\alpha 
= U_{-t} U^{-\infty}(\delta_t * \phi)\alpha \in U_{-t} \cE_+,\]
 so that $\cE_+$ is cyclic. 
\end{prf}

\subsection{The connection to Sinai/Lax-Phillips scattering theory}
\label{se:LaxPhil}

One parameter groups and reflection positivity are closely  related
to the Sinai/Lax-Phillips scattering theory and translation invariant spaces \cite{LP64,LP67,LP81,Sin61}  as
was already noticed in \cite{JOl98}. In short, this theory says that every  
unitary representation
of $\R$ on a Hilbert space $\cE$ satisfying some 
simple conditions stated below can be realized as translation 
in $L^2(\R,\cM)$ for some Hilbert space $\cM$.  

Let $(U,\cE)$ be a unitary representation of $\R$. 
A closed subspace $\cE_+\subset \cE$ is called \textit{outgoing} if
\begin{description}
\item[\rm(LP1)] $\cE_+$ is invariant under $U_t, t > 0$,
\item[\rm(LP2)] $\cE_\infty := \displaystyle \bigcap_{t>0} U_t\cE_+ =\{0\}$,
\item[\rm(LP3)] $\displaystyle \bigcup_{t<0} U_t\cE_+$ is dense in $\cH$.
\end{description}

We have the classical Lax--Phillips Realization Theorem: 
\begin{thm}\label{thm:LP} {\rm(\cite[Thm.~1]{LP64})} 
If $\cE_+$ is outgoing for $(U,\cE )$, then there exists a Hilbert
space $\cM$ such that $\cE\simeq L^2 (\R ,\cM)$,  $\cE_+\simeq L^2( [0,\infty [,\cM)$, and $U$ is represented by translation
$(L_tf)(x)=f(x-t)$. This representation
is unique up to isomorphism of $\cM$.
\end{thm}

This realization of $(U, \cE)$ is called the \textit{outgoing realization of} $U$. 

\begin{rem} \mlabel{rem:outgo-real} It is 
instructive in this contest to recall the construction the Hilbert space $\cM$ and of the isomorphism
$\cE\simeq L^2(\R,\cM)$ as it is presented in \cite{LP81}. Let
$P_+ $ be the orthogonal projection onto $\cE_+$ and define
\[S_t=P_+ U_{-t} \quad \mbox{ for } \quad t\ge 0\, .\]

We claim that $(S_t)_{t\ge 0}$ is a strongly continuous contraction semigroup on 
$\cE_+$ satisfying 
\[ \lim_{t \to \infty} S_tv = 0 \quad \mbox{ for } \quad v \in \cE_+.\] 
This follows from the fact that the projections 
$U_t S_t = U_t P_+ U_{-t}$ onto $U_t \cE_+$ are decreasing and converge strongly 
to $0$ because $\cE_\infty = \{0\}$. 

Write $S_t=e^{-tH}$ in the sense of the Hille--Yoshida Theorem 
(\cite{EN00}). 
For $u\in \cD(H)\subeq \cE_+$, the domain for $H$, define
\[\|u\|_{ \cM}^2:= 2\Re \la H u ,u\ra\ge 0\, .\]
Let $\cK :=\{u\in \cD (H)\: \|u\|_{\cM}=0\}$ and let
$\cM$ be the completion of $\cD (H)/\cK$ in the norm $\|\cdot \|_{ \cM}$. Then $ \cM$ is
a Hilbert space.  We 
write  $[u]$ for the image of $u$ in $\cM$. For $u\in \cD (H)$ define
$f_u: [0,\infty[\to \cM $ by $f_u(t)=[S_t u]$. Then one can show that
\[\|u\|^2=\int_0^\infty \|f_u(t)\|^2_{\cM}\, dt\, .\]
Extend $f_u$ to be zero on $]-\infty ,0[$. Next, if $u\in \bigcup_{t\le 0 }U_t \cE_+$, then there exists
$r>0$ such that $v=U_r u\in \cE_+$. Define $f_u (s)=f_v (s+r)$, $s\ge -r$, and $f_u ( s )=0$ for
$s< -r$. Then $u\mapsto f_u$ is still an isometry and extends to an isometry $\cE \to L^2(\R ,\cM)$ with
the stated properties. 
\end{rem}

We now connect the Lax--Phillips construction to the dilation process.   
The following proposition is an obvious consequence of Theorem \ref{thm:LP} and
Proposition~\ref{prop:e.5}. 

\begin{prop} \mlabel{prop:4.11} Let $(U_t)_{t\in\R}$ be a 
reflection positive unitary one-parameter group on $(\cE,\cE_+,\theta)$ 
for which $\cE_+$ is cyclic and $\cE_{\rm fix} = \{0\}$. Then 
$\cE_+$ is outgoing, so that 
$(U,\cE)$ is unitarily equivalent to the translation representation 
on $L^2(\R,\cM)$ for some Hilbert space $\cM$. 
This realization is unique up to   isomorphism of $\cM$.
\end{prop}

The preceding proposition suggests to approach the structure of reflection 
positive one-parameter groups in an outgoing realization, but this turns out 
to be rather difficult because the involution $\theta$ is not well visible 
in this picture. 

\begin{rem} \mlabel{rem:3.15} Starting with the subspace 
$\cE_+ := L^2(\R_+,\cM)$ of $\cE = L^2(\R,\cM)$ and 
$(U_t f)(x) := f(x-t)$, we want to describe those unitary 
involutions $\theta$ for which $\cE_+$ is $\theta$-positive and 
$\theta U_t \theta = U_{-t}$ for $t \in \R$. 

If $\cF(f) = \hat f$ denotes the Fourier transform $\cE \to \cE$, then 
the commutant of $\cF U_\R \cF^{-1}$ is $L^\infty(\R,B(\cM))$, which implies 
that $\theta$ can be described as 
\[ \widehat{\theta f} (x) = m(x) \hat f(-x) \] 
for a measurable function $m \: \R \to \U(\cM)$ satisfying 
$m(-x) = m(x)^*$ for almost every $x \in \R$. 

To determine for which such functions $m$ the subspace $\cE_+$ is $\theta$-positive, 
we calculate for $f,g \in C^\infty_c(\R_+)$ and $v,w \in \cM$ as follows: 
\begin{align*}
 \la f v,\theta (gw) \ra 
&= \la \hat f v, \widehat{\theta(g w)} \ra 
= \int_{\R}  \hat f(x) \hat g^*(x) \la v, m(x) w\ra\, dx \\
&= \int_{\R}  \widehat{f * \oline g}(x) \la v, m(x) w\ra\, dx 
= \int_{\R}  (f * \oline g)(x) \la v, \cF^{-1}(m)(x) w\ra\, dx \\
&= \int_{\R_+} \int_{\R_+} f(x) \oline{g(y)}\la v, \cF^{-1}(m)(x+y) w\ra\, dx\, dy.  
\end{align*}
We conclude that $\cE_+$ is $\theta$-positive if and only if 
the restriction of $\cF^{-1}(m)$ to $\R_+$ is a 
positive definite $B(\cM)$-valued distribution with respect to the 
involution $s^\sharp = s$ on the semigroup~$\R_+$. 

For $\cM = \C$, it follows from 
\cite[Thm.~4.13]{NO12} that $\cF^{-1}(m)\res_{\R_+}$ is an analytic function 
which is the Laplace transform $\cL(\rho)(x) = \int_\R e^{-tx}\, d\rho(t)$ of a Radon measure $\rho$ on $[0,\infty[$. Actually the above calculation then implies that 
$\hat\cE \cong L^2([0,\infty[,\rho)$.

As we shall see in Subsection~\ref{ex:3.13} below, 
for every $\lambda > 0$, an interesting example with $\cM = \C$ is given by 
the function 
\begin{equation}\label{eq:m}
m(x) := \frac{\lambda - ix}{\lambda + ix} 
= -1 + \frac{2\lambda}{\lambda + ix}\, .
\end{equation} 
Then 
\[ \cF^{-1}(m) = - \delta_0 + 2 \lambda \sqrt{2\pi} e^\lambda \quad \mbox{ with} 
\quad e^\lambda(t) = e^{-\lambda t} \chi_{\R_+}.\] 
Restricting this distribution to $\R_+$, yields the Laplace transform 
of a multiple of the Dirac measure $\delta_\lambda$. 
This example already shows the difficulty of recovering $m$ from 
$\rho$ because this requires an extension of $\cL(\rho)$ to a 
distribution on $\R$ whose Fourier transform should be a $\T$-valued function.
\end{rem}

\subsection{The case where $\cE_0$ is cyclic}
\mlabel{subsec:3.3} 

In this section we explore reflection positive one-parameter groups 
$(U_t)_{t \in \R}$ on $(\cE,\cE_+,\theta)$ for which the subspace $\cE_0 = \cE_+^\theta $ 
(cf.~Definition~\ref{def:x.1}) is cyclic in $\cE$. Then 
the GNS construction permits us to reconstruct 
$U$ from the positive definite function 
$\phi(t) := P U_t P^*$, where $P \: \cE \to \cE_0$ is the 
orthogonal projection (Example~\ref{ex:vv-gns}). 
This function is reflection positive for $(\R,\R_+,-\id_\R)$ in the sense of 
\cite{NO12}, i.e., its restriction $\psi := \phi\res_{\R_+}$ is positive definite with respect to the 
trivial involution $s^\sharp = s$ and $\phi(t) = \psi(|t|)$. Since the latter formula 
produces for every positive definite function $\psi$ on $\R_+$ a positive definite 
function on $\R$ (cf.\ \cite[\S I.8.2, p.~29]{SzNBK10}, \cite[Rem.~3.2]{NO12}), every 
positive definite function on $\R_+$ extends to a reflection positive function on $\R$. 

Our starting point below is the observation that 
$\phi\res_{\R_+}$ is a representation if and only if $q|_{\cE_0}:\cE_0\to \hat\cE$ 
is unitary (Proposition~\ref{prop:3.9}), i.e., we can identify $\hat\cE$ with $\cE_0$. 
In this sense we then have $\hat U = \phi\res_{\R_+}$ and 
$U$ is called the {\it minimal unitary dilation} of the contraction semigroup $\hat U$. 
The following lemma provides a criterion for $q(\cE_0) = \hat\cE$. 

\begin{lem} \mlabel{lem:3.5} Let $(\cE,\cE_+,\theta)$ be a reflection positive 
Hilbert space. If $\cD \subeq \cE_0$ is a subspace whose image in $\hat\cE$ is dense, 
then the following assertions hold: 
\begin{description}
\item[\rm(i)] $\cD$ is dense in $\cE_0$. 
\item[\rm(ii)] $q(\cE_0) = \hat\cE$. 
\item[\rm(iii)] $\cN = \cE_+ \cap \cE_0^\bot = \cE_+ \cap \theta(\cE_+)^\bot$. 
\item[\rm(iv)] $\cN \oplus \cE_0 \oplus \theta(\cN)$ is an orthogonal direct sum, 
hence in particular, a closed subspace of $\cE$. 
\end{description}
\end{lem}

\begin{prf} (i),(ii) Since $q\res_{\cE_0}$ is isometric, the density of $q(\cD)$ 
in $\hat\cE$ implies that $q(\cE_0) = \hat\cE$ and that $\cD$ is dense in $\cE_0$. 

(iii) For $v_1 \in \cE_0^\bot \cap \cE_+$ and $v_0 \in \cE_0$ we have 
\[\la q(v_0), q(v_1) \ra = \la \theta v_0, v_1 \ra = \la v_0, v_1 \ra = 0\, ,\]
so that $q(v_1) = 0$ follows from $q(\cE_0) = \hat\cE$. 
On the other hand, $\cE_+ = \cE_0 \oplus (\cE_0^\bot \cap \cE_+)$, 
and $q\res_{\cE_0}$ is isometric. This implies $\cN = \cE_+ \cap \cE_0^\bot$. 
The relation $\cN = \cE_+ \cap \theta(\cE_+)^\bot$ follows from the fact that 
$q(v) = 0$ is equivalent to $\la \theta w, v\ra = 0$ for every $w \in \cE_+$. 

(iv) The orthogonality of the decomposition follows from (iii), and this in turn implies 
that the sum is a closed subspace. 
\end{prf}

For the notion of a positive definite function on an involutive semigroup 
and the corresponding GNS construction we refer to Example~\ref{ex:vv-gns}.
We now apply Lemma~\ref{lem:mult} 
to reflection positive one-parameter groups. 

\begin{prop} \mlabel{prop:3.9} 
Let $(U_t)_{t \in \R}$ be a reflection positive one-parameter 
group on $(\cE, \cE_+, \theta)$ for which 
$\cE_0 = \cE_+^\theta$ is cyclic and write 
$P \:  \cE \to \cE_0$ for the orthogonal projection. 
Then the following are equivalent: 
\begin{description}
\item[\rm(i)] $q(\cE_0) = \hat\cE$, i.e., $\Gamma := q\res_{\cE_0} \: \cE_0 \to \hat\cE$ is 
unitary. 
\item[\rm(ii)] The reflection positive function 
$\phi \: \R \to B(\cE_0)$, $\phi(t) := P U_t P^*$, is multiplicative 
on $\R_+$. 
\item[\rm(iii)] The orthogonal projections 
$P_+$ onto $\cE_+$, $P_0$ onto $\cE_0$, and $P_-$ onto $\theta(\cE_+)$ satisfy 
the Markov condition $P_+ P_0 P_- = P_+ P_-$.
\begin{footnote}
  {This is an abstraction of the Markov condition for Osterwalder--Schrader positive 
processes that one finds in \cite{Kl77, Kl78}.}
\end{footnote}
\end{description}
If this is the case, then 
\begin{description}
\item[\rm(a)] $\phi(t) = \Gamma^* \hat U_{|t|} \Gamma$ for $t \in \R$. In particular, 
$\Gamma$ intertwines $\phi\res_{\R_+}$ with $\hat U$. 
\item[\rm(b)] $\cE = \cN \oplus \cE_0 \oplus \theta(\cN)$ is an orthogonal decomposition 
and $\cE_+$ is maximal $\theta$-positive. 
\item[\rm(c)] $\cE_+$ is the closed subspace generated by $\bigcup_{t > 0} U_t\cE_0$. 
\end{description}
\end{prop}

\begin{prf} Let $S := [0, \infty[ = \{0\} \cup \R_+$. 
Since $\cE_0$ is $U$-cyclic, the map 
\[ \Phi \: \cE \to (\cE_0)^\R, \quad 
\Phi(v)(t) := P U_t v \] 
is an equivalence of the representation $U$ of $\R$ on $\cE$ to the 
GNS-representation defined by $\phi$ (cf.~Example~\ref{ex:vv-gns}) 
and the representation $\hat U$ of $S$ on $\hat\cE$ is equivalent to the 
GNS-representation  defined by $\phi\res_S$, where the map 
$q \: \cE_+ \to \hat\cE$ simply corresponds to the restriction 
$f \mapsto f\res_S$  (\cite[Prop.~1.11]{NO12}). 
The inclusion $\iota \: \cE_0 \into \cH_\phi$ is given by 
$\iota(v)(t) = P U_tv = \phi(t)v$ for $t \in \R$, and likewise the inclusion 
$\hat\iota \: \cE_0 \into \cH_{\phi\res_S}$ is given by 
$\iota(v)(t) = \phi(t)v$ for $t \geq 0$. 
Therefore (i) is equivalent to the surjectivity of the 
inclusion $\hat\iota$. In view of Lemma~\ref{lem:mult}, this is equivalent 
to the multiplicativity of $\phi\res_S$, which is equivalent to the 
multiplicativity on $\R_+$. 

To see that (i) is equivalent to (iii), we first observe  that the Markov 
condition is equivalent to 
\[\la P_0 \theta(v), w \ra = \la \theta(v), w \ra\quad \text{ for all }\quad v,w \in \cE_+\, .\] 
This in turn means that $\|q(v) \| = \|q(P_0 v)\|$ for $v \in \cE_+$, which 
implies that $q(\cE_+) = q(\cE_0)$, so that (i) follows from Lemma~\ref{lem:3.5}. 
Conversely, (i) implies that $\cN = \cE_+ \cap \cE_0^\bot$, which leads to 
$\|q(v)\| = \|q(P_0v)\|$ for $v \in \cE_+$. 

(a) Now we assume that (i) is satisfied, so that 
$\Gamma$ is unitary, 
and Lemma~\ref{lem:3.5}(ii) implies that $q = \Gamma \circ P\res_{\cE_+}$. 
For $t > 0$, the relation 
$\hat U_t \circ q = q \circ U_t\res_{\cE_+}$
leads to 
$\hat U_t \Gamma P\res_{\cE_+} = \Gamma P U_t\res_{\cE_+},$
so that $\Gamma^* \hat U_t \Gamma  = P U_t P^* = \phi(t),$
i.e., $\Gamma$ intertwines $\phi(t)$ with $\hat U_t$. 

(b), (c) Let $\cE^+ \subeq \cE$ be the closed subspace generated by the translates 
$U_t \cE_0$, $t > 0$, and note that $\cE^+ \subeq \cE_+$. 
For $t \geq 0$ we then have $U_t \cE_0 \subeq \cE^+$, 
and for $t < 0$ we have $U_t \cE_0 = \theta(U_{-t}\cE_0)\subeq \theta(\cE^+)$. 
Therefore $\cE^+ + \theta(\cE^+)$ is dense in $\cE$. 
For $\cN^+ := \cN \cap \cE^+$ we obtain the orthogonal decomposition 
\[ \cE^+ + \theta(\cE^+) 
= \cN^+ \oplus \cE_0 \oplus \theta(\cN^+) \] 
(Lemma~\ref{lem:3.5}(iv)), which implies that $\cE^+ + \theta(\cE^+)$ is closed, hence equal 
to $\cE$. 

To see that $\cE^+$ is a maximal $\theta$-positive subspace, let 
$\cF \subeq \cE$ be a $\theta$-positive subspace containing $\cE^+$. 
Then $\cF = \cN^+ \oplus \cE_0 \oplus (\cF \cap \theta(\cN^+))$, 
and for every $v \in \cF \cap \theta(\cN^+)$ we have $w := v - \theta v \in \cF$ with 
$\la \theta w, w \ra = - \la w,w \ra$. Therefore the $\theta$-positivity of $\cF$ implies that 
$w = 0$, so that $v = 0$ follows from $v \bot \theta v$. This shows that 
$\cE^+$ is maximal $\theta$-positive. In particular, we obtain $\cE^+ = \cE_+$. 
\end{prf}

\begin{ex}(OS-positive processes) \mlabel{ex:klein} 
(a) The following class of examples arises in A.~Klein's papers \cite{Kl77, Kl78} 
on path space models for Quantum Field Theory. 
Let $(X,\Sigma,\mu)$ be a probability space, 
$\Sigma_0$ be a $\sigma$-subalgebra of $\Sigma$, 
and assume that $\R$ acts by automorphisms $U_t$ on the probability 
space $(X,\Sigma, \mu)$. Further, $\tau$ is an involution of $(X,\Sigma,\mu)$ satisfying 
$\tau U_t \tau = U_{-t}$ and preserving $\Sigma_0$. It is further assumed that 
$\Sigma$ is generated by the subalgebras $U_t\Sigma_0$, $t \in \R$. 
We write  $\Sigma_\pm$ for the 
subalgebras generated by $\bigcup_{\pm t \geq 0} U_t \Sigma_0$ 
and note that $\Sigma_+$ is invariant under~$\R_+$. 

On the Hilbert $\cE := L^2(X,\Sigma,\mu)$, we then obtain a unitary one-parameter 
group $(U_t)_{t \in \R}$, a unitary involution $\theta(f) := f \circ \tau$, 
and a closed subspace $\tilde\cE_0 := L^2(X,\Sigma_0,\mu)$.
\begin{footnote}{The Hilbert space $L^2(X,\Sigma,\mu)$ has a 
natural realization as a space of functions on $\Sigma$ by the transform 
$\tilde f(E) := \la f, \chi_E \ra$. This is a realization as a reproducing 
kernel Hilbert space on $\Sigma$ with the kernel 
$K(E,F) = \mu(E \cap F)$. Here we only use that the characteristic functions 
of the elements of $\Sigma$ form a total subset of $L^2(X,\Sigma,\mu)$. 
In this picture, we obtain for each smaller $\sigma$-algebra 
$\Sigma_0 \subeq \Sigma$ a natural restriction map which corresponds to 
the conditional expectation. This restriction map 
restricts to a unitary isomorphism of the closed subspace 
generated by the characteristic functions 
$\chi_E, E \in \Sigma_0$, with $L^2(X,\Sigma_0,\mu)$. 
}\end{footnote}

Then $\cE_+ := L^2(X,\Sigma_+, \mu) \subeq \cE$ is the smallest closed subspace containing 
$U_t\tilde\cE_0$ for $t > 0$. Klein's condition of {\it Osterwalder--Schrader positivity}, 
{\it OS-positivity} for short, is equivalent to the 
$\theta$-positivity of $\cE_+$. 
Since $\tilde\cE_0$ is $\theta$-invariant and contained in $\cE_+$, this assumption 
implies in particular that $\theta f = f$ for $f \in \tilde\cE_0$, i.e., 
$\tilde\cE_0 \subeq \cE_0$. 
We thus obtain a reflection positive one-parameter group for which $\cE_0$ is cyclic. 

(b) Typical situations as described under (a) arise in QFT as follows. 
Let $\cH$ be a separable real Hilbert space. 
Then there exists a probability space $(X,\Sigma,\mu)$ 
and, for each $v \in \cH$, a random variable $\phi(v)$ on $X$ 
such that any tuple $(\phi(v_1), \ldots, \phi(v_n))$ is jointly 
Gaussian with covariance $(\la v_i, v_j\ra)_{1 \leq i,j \leq n}$ 
(\cite[Thm.~2.3.4]{Sim05}). Since $(X,\Sigma,\mu)$ is uniquely determined 
up to isomorphisms of measure spaces, we obtain a corresponding 
action of the orthogonal group $\OO(\cH)$ on $(X,\Sigma,\mu)$. 

If we start with a reflection positive one-parameter 
group $(U_t)_{t \in \R}$ on $(\cE,\cE_+, \theta)$ 
for which $\cE_0$ is $U$-cyclic and $\cE_+$ is generated by 
$U_t \cE_0$, $t > 0$, then all this structure is reflected 
in $(X,\Sigma,\mu)$. In particular, the $U_t$ and $\theta$ define 
automorphisms of measure spaces. 
We write $\Sigma_0 \subeq \Sigma$ for the $\sigma$-algebra generated by 
$\phi(\cE_0)$. Then $\Sigma_+$ is generated by the translates 
$U_t\Sigma_0$ and all assumptions from (a) are satisfied. 

Here $\cF := L^2(X,\Sigma,\mu)$ is the Fock space of $\cH$ 
and $\cF_+ := L^2(X,\Sigma_+,\mu)$ can be identified with the 
Fock space of $\cE_+$. We thus arrive at the situation from 
Remark~\ref{rem:1.7} which implies in particular that the passage 
to the Fock space, i.e., second quantization, leads to a reflection 
positive one-parameter group on $\cF(\cE)$. 

(c) Other situations as in (a) arise as follows. 
Let $M$ be a smooth manifold, $M_+ \subeq M$ open and 
$\tau \: M \to M$ a smooth involution whose fixed point set $M^\tau$  
is a hypersurface such that $M \setminus M^\tau$ is the disjoint 
union of $M_+$ and $\tau(M_+)$.  In addition, we 
assume that $(\sigma_t)_{t \in \R}$ is a one-parameter group of diffeomorphisms of 
$M$ preserving $D$ such that $\tau \sigma_t \tau = \sigma_{-t}$, 
$\sigma_t M_+ \subeq M_+$ for $t \geq 0$ and 
$\bigcup_{t >0} \sigma_{-t} M_+ = M$. 

On the space $X := C^{-\infty}(M)$ of distributions  
we consider the smallest $\sigma$-algebra $\Sigma$ for which 
all evaluation maps $\phi^*(D) := D(\phi)$, $\phi \in C^\infty_c(M)$, are measurable. 
Then $\alpha_t(D)(\phi) := D(\phi \circ \sigma_t)$ defines a one-parameter group 
of measurable isomorphisms on $(X,\Sigma)$. 

For every open subset $O \subeq M$ we obtain a sub-$\sigma$-algebra 
$\Sigma(O)$ which is minimal with the property that all function $\phi^*$, 
$\supp(\phi) \subeq O$, are measurable. Then 
$\Sigma_+ := \Sigma(M_+)$ is invariant under $(\alpha_t)_{t > 0}$, and 
$\bigcup_{t < 0} \alpha_t \Sigma_+$ generates $\Sigma$. 

If $\mu$ is a reflection positive probability measure on $(X,\Sigma)$ with respect to 
$(M,M_+,\tau)$ (cf.\ Example~\ref{ex:1.5}(a)) 
which is invariant under $(\alpha_t)_{t \in \R}$, then we obtain a 
reflection positive unitary one-parameter group 
$U_t F := F \circ \alpha_{-t}$ on $\cE := L^2(X,\mu)$ with respect to 
$(\theta F)(D) := F(\tau^*D)$ and 
the subspace $\cE_+$ generated by the functions $e^{i\phi^*}$, 
$\supp(\phi) \subeq M_+$. 

To match this with (a), we further need a $\sigma$-algebra $\Sigma_0 \subeq \cE_+$ 
such that $\alpha_t \Sigma_0$, $t \geq 0$, generates $\Sigma_+$. 
This comes usually from extending the natural map 
$C^\infty_c(M) \to \cE, \phi \mapsto e^{i\phi^*}$ to the space 
$C^\infty_c(M^\tau)$ of test functions on $M^\tau$. Then a natural 
choice for $\Sigma_0$ is the $\sigma$-algebra generated by 
$e^{i\phi^*}$, $\phi \in C^\infty_c(M^\tau)$. 
We shall see below that, for 
$(M,M_+,\tau) = (\R^d, \R^d_+,\theta)$, this situation 
arises for generalized free fields with a non-trivial ``time-zero'' subspace, 
i.e., for which the measure $\tilde\mu$ is tempered (see Remark~\ref{rem:d.10} and 
also \cite{GJ81}). 
\end{ex}

\subsection{The Hardy space of the real line} \mlabel{ex:3.13}
Some of the ideas and tools developed
here will show up again in Section \ref{ex:4.11}.
We consider the Hilbert space $\cE := L^2(\R)$ with the 
translation action $(U_t f)(x) := f(x-t)$ and the subspace 
$\cE_+ := L^2(\R_+)$ of all functions supported in $\R_+$. Note that 
$U_t \cE_+ \subeq \cE_+$ holds for every $t > 0$. 
The space $\cE_+$ can also be viewed as the Hardy space of the right complex half 
plane~$\C_+ = \R_+ + i \R$. 
To obtain this identification, one first observes that, for $\Re z > 0$, the functions 
$Q_z(x) := e^{-x\oline z} \chi_{\R_+}$ form a total subset of $\cE_+$ with 
\begin{equation}
  \label{eq:hardyker}
K(z,w) := \la Q_w, Q_z \ra = \int_0^\infty e^{-x(z + \oline w)}\, dx 
= \frac{1}{z + \oline w},
\end{equation}
that is the reproducing kernel of the Hardy space 
$\cH_K \subeq \cO(\C_+)$ that can be specified by 
\[ f(z) = \la f, K_z \ra \quad \mbox{ for } \quad z \in \C_+, \quad 
K_z(w) := K(w,z), \quad f \in \cH_K.\] 
Then the Laplace transform 
\[ \cL \: L^2(\R_+) \to \cH_K, \quad 
\cL(f)(z) := \la f, Q_z \ra = \int_0^\infty f(x) e^{-z x}\, dx \]
is unitary and satisfies 
\[ \cL(U_t f)(z) 
= \la f, U_{-t} Q_z \ra 
= \la f, e^{-t\oline z} Q_z \ra 
= e^{-t z} \la f, Q_z \ra = e^{-t z} \cL(f)(z).\] 
Note that multiplication with the function $e_{-t}(z) := e^{-tz}$ is isometric on the Hardy 
space because its boundary values on $i\R$ have absolute value~$1$ and 
\[ \|f\|_{\cH_K}^2 = \lim_{z \to 0} \int_{\R} |f(z + ix)|^2\, dx.\] 

For $\lambda > 0$, we consider on $L^2(\R)$ the unitary involution defined by 
\begin{equation}
  \label{eq:theta1}
\widehat{\theta f}(x) := m_\lambda(x) \hat f(-x), \quad \mbox{ where } \quad 
m_\lambda(x) := \frac{\lambda - ix}{\lambda + ix} 
= -1 + \frac{2\lambda}{\lambda + ix}\, ,
\end{equation}
compare to (\ref{eq:m}).
Note that $\theta^2 = \1$ follows from $m_\lambda(-x) = m_\lambda(x)^{-1}$. This involution satisfies 
$\theta U_t \theta = U_{-t}$ for $t \in \R$.

\begin{lem}\label{le:thPos} The following assertions hold: 
\begin{enumerate}
\item[\rm(i)] $L^2(\R_+)$ is $\theta$-positive. 
\item[\rm(ii)] If $e^\lambda (t)= e^{-\lambda t}\chi_{\R_+}$, 
then $\cE_0=\C e^\lambda$ is one dimensional.
\item[\rm(iii)] The quotient map $q: \cE_+\to \hat{\cE}\cong \C$ is equivalent to 
$q(f)=\sqrt{2\lambda}\, \cL (f)(\lambda)$. 
\end{enumerate}
\end{lem}

\begin{prf} (i) It suffices to show that, 
for $f \in \cS(\R)$ with $\supp(f) \subeq [0,\infty[$, the following integral is non-negative: 
\begin{align*}
 \la f, \theta f \ra 
&= \la \hat f, \widehat{\theta f} \ra 
=  \int_{\R} \oline{m_\lambda(x)} \hat f(x) \oline{\hat f(-x)}\, dx 
= \int_{\R} \frac{\lambda + ix}{\lambda - ix} \hat f(x) {\hat f}^*(x)\, dx 
= 2\lambda \int_\R \frac{1}{\lambda -ix} \hat f(x) {\hat f}^*(x)\, dx, 
\end{align*}
where we have used \eqref{eq:theta1} and 
\[ \int_\R \hat f(x) {\hat f}^*(x)\, dx  = \int_\R f(x) \oline{f(-x)}\, dx = \int_0^\infty f(x) \oline{f(-x)}\, dx = 0\]
as $f$ is supported on $\R_+$. 
>From \eqref{eq:hardyker} we further derive by taking boundary values on $i\R$ that 
the Fourier transform of $e^\lambda$ is given by 
\[ \gamma_\lambda(y) := \frac{1}{\sqrt{2\pi}} \frac{1}{\lambda + i y}.\] 
It now remains to verify the non-negativity of the integral: 
\begin{align*}
\la f, \theta f \ra 
&= 2\lambda \int_\R \frac{1}{\lambda -ix} \hat f(x) \hat f^*(x)\, dx 
= 2\lambda \sqrt{2\pi} \int_\R \oline{\gamma_\lambda (x)} \hat f(x) \hat f^*(x)\, dx \\
&=2\lambda \sqrt{2\pi} \la \widehat{f * \oline f}, \gamma_\lambda\ra_{L^2(\R)}
= 2\lambda \sqrt{2\pi}\la f * \oline f, e^\lambda\ra_{L^2(\R)}\\
&= 2\lambda \sqrt{2\pi}\cL(f * \oline f)(\lambda) 
= 2\lambda |\cL(f)(\lambda)|^2 \geq 0.
\end{align*}
This proves that $\cE_+ = L^2(\R_+)$ is $\theta$-positive, 
that $\hat\cE \cong \C$, and that the quotient map $q \: \cE_+ \to \cE$ can be identified 
with 
$q(f) = \sqrt{2\lambda} \cL(f)(\lambda).$ 
Hence the natural map $\cH_K \to \hat\cE \cong \C$ is given by evaluation in $\lambda \in \C_+$. 
>From 
\[ \widehat{\theta(e^\lambda)}(x) 
= \frac{\lambda - ix}{\lambda+ ix} \gamma_\lambda(-x) 
= \frac{1}{\sqrt{2\pi}}\frac{\lambda - ix}{\lambda+ ix} 
\frac{1}{\lambda - i x}
= \frac{1}{\sqrt{2\pi}}\frac{1}{\lambda+ ix}  = \widehat{\, e^\lambda\,}(x), \]
we get $\theta(e^\lambda) = e^\lambda$. Since $q\res_{\cE_0}$ is isometric and 
$\hat\cE$ is one-dimensional, 
it follows that $\cE_0 = \C e^\lambda$. 
\end{prf}

\begin{rem} We remark that the function $e^\lambda$ is cyclic in $\hat\cE$. This follows from the 
fact that its Fourier transform $\gamma_\lambda$ has no zeros. Hence $e^\lambda$ generates 
$L^2(\R)$ with respect to the multiplicative representation. We thus have a non-trivial 
example to which Proposition~\ref{prop:3.9} applies. 
\end{rem}

\section{An $L^2$-model  of the dilation Hilbert space} 
\mlabel{sec:dil}

Suppose that the one-parameter group 
$(U_t)_{t \in \R}$ is reflection positive with respect to 
$(\cE,\cE_+, \theta)$ and that the subspace $\cE_0 = \cE_+^\theta$ is $U$-cyclic. 
As we shall see in Section~\ref{ex:4.5}, this is not always the case. 
Then $(U,\cE)$ is equivalent to the GNS representation defined by the positive definite 
operator-valued function $\phi(t) := PU_t P^*$, where 
$P \: \cE \to \cE_0$ is the orthogonal projection. 
A particularly interesting case arises if $\cE_0$ is ``as big as possible'', 
i.e., if $q(\cE_0) =\hat\cE$, so that $q\res_{\cE_0} \: \cE_0 \to \hat\cE$ is unitary. 
In view of Proposition~\ref{prop:3.9}, this is 
equivalent to $\phi$ being the minimal unitary dilation $\phi(t) = \hat U_{|t|}$ 
of the hermitian contraction semigroup $(\hat U_t)_{t > 0}$ 
on $\hat\cE$. We consider the GNS representation 
associated to $\phi(t) = \hat U_{|t|}$ on $\cE \subeq C(\R,\hat\cE)$ as a 
``maximal'' euclidean realization of $(\hat\cE,\hat U)$. 
In this section we first discuss an $L^2$-model of the dilation space  $\cE$ 
in which we specify certain ``spectral subspaces'' $\cE_F$, $F \subeq \R$ a closed subset, 
by the requirement that the inverse Fourier transform of the elements of $\cE_F$ is supported 
in $F$. Here the most important piece of information is that 
$\cE_+ = \cE_{[0,\infty[}$ and $\cE_0 = \cE_{\{0\}}$, which can be viewed 
as a realization of $\cE_0$ by distributions supported in $\{0\}$. 
This is why, in Quantum Field Theory, the isomorphism 
$\hat\cE \cong \cE_0$ corresponds to the existence 
of a so-called ``time-zero realization'' of the Hilbert space 
$\hat\cE$ (\cite[Thm.~6.2.6]{GJ81}; see also 
 Remark~\ref{rem:d.10} and Example~\ref{ex:5.14}). 
We conclude this section with an explicit outgoing realization of $\cE$ 
in the sense of Lax--Phillips/Sinai scattering theory.

\subsection{The spectral subspaces $\cE_F$}\mlabel{sec:EF}
We start with a strongly continuous 
hermitian contraction semigroup $(C_t)_{t > 0}$ on the Hilbert space $\cH$. Here 
strong continuity means that $\lim_{t \to 0} C_t v = v$ for every 
$v \in \cH$ (\cite{EN00}), so that $C_0 := \id_\cH$ yields a strongly continuous extension 
to $[0,\infty[$. First we observe that $C$ always has a natural euclidean realization. 
We shall see below that this realization is rather large. 

\begin{prop} \mlabel{prop:exist-real} 
For every strongly continuous one-parameter 
semigroup $(C_t)_{t \geq 0}$ of hermitian contractions on a Hilbert space $\cH$, 
there exists a euclidean realization $(U_t)_{t \in \R}$ on 
$(\cE,\cE_+, \theta)$ for which 
$q \: \cE_0 \to \cH$ is  unitary and 
$\cE_0$ is $U$-cyclic in $\cE$. 
Any realization with these two properties 
is equivalent to the one obtained by dilation, i.e., 
from the $B(\cH)$-valued positive definite function $\phi(t) := C_{|t|}$ on $\R$. 
\end{prop}

\begin{prf} Let $P$ be the spectral measure of $C$, so that 
$C_t = \int_0^\infty e^{-t\lambda}\, dP(\lambda)$ 
in the sense of spectral integrals. 
Then \cite[Prop.~3.1]{NO12} implies that 
\[\phi(t) := C_{|t|} = \int_0^\infty e^{-\lambda|t|}\, dP(\lambda)\]
is a strongly continuous reflection positive $B(\cH)$-valued function for the triple 
$(\R,-\id,\R_+)$. The GNS construction yields a unitary one-parameter group 
$U_\phi$ on the corresponding reproducing kernel 
Hilbert space $\cE_\phi \subeq C(\R,\cH)$ in which 
$\cH$ is embedded via $\cH \to \cE_\phi, v \mapsto \phi \cdot v$ 
(cf.\ Example~\ref{ex:vv-gns}). 
>From \cite[Prop.~1.11]{NO12} it follows that 
$U_\phi$ is reflection positive with respect to the triple 
$(\cE_\phi, \cE_{\phi,+}, \theta)$, where 
$\cE_{\phi_+}$ is the closed subspace generated by the functions 
$\phi(\cdot + t)v$, $t > 0$, $v \in \cH$, $\theta(f)(t) := f(-t)$, 
and $\hat{U_\phi} = C$ on $\hat{\cE_{\phi}} \cong \cH$ with 
$q(f) = f\res_{\R_+}$. 
\end{prf}

Since the space $\cE_\phi$ is rather implicit, 
it is natural to develop a more explicit realization as an $L^2$-space. 
The tools for such a realization have been developed in \cite{Ne98} 
in terms of operator-valued measures. 
As in the proof of Proposition~\ref{prop:exist-real}, we write 
$P \: \fB([0,\infty[) \to B(\cH)$ for the corresponding spectral measure and 
$A = \int_0^\infty \lambda\, dP(\lambda)$ for its infinitesimal generator which satisfies 
$C_t = e^{-tA}$ for $t \geq 0$. Reformulating the integral representation 
of $\phi$ we now obtain 
\begin{align*}
\phi(t) 
&= e^{-|t|A} =  \int_{[0,\infty[} e^{-|t|\lambda} \, dP(\lambda) 
= \1_{\cH_{\rm fix}} \oplus \int_{\R_+}  e^{-|t|\lambda} \, dP(\lambda) \\
&= \1_{\cH_{\rm fix}} \oplus \frac{1}{\pi} \int_{\R_+} \int_{\R}  e^{-ixt} \frac{\lambda}{\lambda^2+x^2}\, dx
\, dP(\lambda) 
= \1_{\cH_{\rm fix}} \oplus \frac{1}{\pi} \int_\R \int_{\R_+}  e^{-ixt} \frac{\lambda}{\lambda^2+x^2}\, 
\, dP(\lambda)\, dx \\
&= \1_{\cH_{\rm fix}} \oplus \frac{1}{\pi} \int_\R e^{-ixt} A(A^2 + x^2)^{-1}\, dx.  
\end{align*}
Splitting off the subspace $\cH_{\rm fix}$ of fixed vectors, we assume in the following that 
$\cH_{\rm fix} = \{0\}$, which is equivalent to $P(\{0\}) = 0$, resp., $\ker A = \{0\}$. 
This is in particular justified by Remark~\ref{rem:3.8}.  The above formula expresses the operator-valued function $\phi$ as the 
Fourier transform of the bounded operator-valued measure 
\[ Q(x) dx  \quad \mbox{ for } \quad Q(x) := \frac{1}{\pi} A(A^2 + x^2)^{-1} \] 
with total mass $\int_\R Q(x)\, dx = \1$. Let $\cE := L^2(\R, Q;\cH)$ be the corresponding 
$\cH$-valued $L^2$-space, in which the scalar product is given by 
\begin{equation}
  \label{eq:scalprodincE}
\la f, g \ra _{\cE}
= \int_\R \la Q(x)f(x), g(x)\ra\, dx 
= \frac{1}{\pi}  \int_\R \la A(A^2 + x^2)^{-1}f(x), g(x)\ra\, dx
\end{equation}
(cf.\ \cite{Ne98}). 
For constant functions $v,w \in \cH$, we have $\la v,w\ra_\cE = \la v, w \ra_\cH$, 
so that the inclusion $\iota \:  \cH \into \cE$ is isometric. 
In the following we identify $\cH$ with the subspace $\iota(\cH)$ of constant functions in $\cE$. 
The orthogonal projection $P_\cH = \iota^*\: \cE \to \cH$ is given by 
\begin{equation}
  \label{eq:PH}
P_\cH = \frac{1}{\pi}  \int_\R Q(x)\, dx, \quad \mbox{ resp.} \quad 
P_\cH f = \frac{1}{\pi}  \int_\R A(A^2 + x^2)^{-1}f(x)\, dx.
\end{equation}

For the unitary one-parameter group 
\[ (U_t f)(x) = e^{-itx} f(x), \quad t,x\in \R,f \in \cE,  \] 
we then obtain for $v,w \in \cH$ 
\[ \la P_\cH U_t P^*_\cH v, w \ra 
 =   \la U_t v, w \ra 
= \frac{1}{\pi}  \int_\R e^{-itx} \la A(A^2 + x^2)^{-1}v,w\ra\, dx 
= \la \phi(t)v,w\ra.\] 
We conclude that 
\begin{equation}\label{eq:3.3.1}
\phi(t) = P_\cH U_t P^*_\cH \quad \mbox{ for } \quad t \in \R, \quad \mbox{ and hence }\quad 
C_t = P_\cH  U_t P^*_\cH \quad \mbox{ for } \quad t > 0.
\end{equation}
Since $\cH$ is $U$-cyclic in $\cE$ (cf.~\cite[Prop.~II.9]{Ne98}), 
the map 
\[ \Gamma  \: \cE \to C(\R,\cH), \quad 
\Gamma(f)(t) := P_\cH U_t f \] 
is a unitary equivalence from the 
representation of $\R$ on $\cE$ to the GNS representation 
$(\pi_\phi, \cH_\phi)$, defined by the positive definite 
function~$\phi$ (cf.~Example~\ref{ex:vv-gns}). 
On $\cE$ the unitary involution given by $(\theta f)(x) := f(-x)$ satisfies 
$\theta U_t \theta = U_{-t}$ and 
\[ \Gamma(\theta f)(t) = P_\cH U_{t} \theta f = P_\cH \theta U_{-t} f = P_\cH U_{-t} f 
= \Gamma(f)(-t),\] 
so that the proof of Proposition~\ref{prop:exist-real} implies that 
the closed subspace $\cE_+$ generated by 
$U_t \cH, t > 0$, is $\theta$-positive and the 
corresponding contraction semigroup on $\cH$ is $\hat U = C$. 
We conclude that the one-parameter group $(U_t)_{t \in \R}$ on $(\cE,\cE_+,\theta)$ 
is a euclidean realization of $(C,\cH)$. 

Lemma~\ref{lem:3.5}(iii) implies that the map $q \:  \cE_+ \to \hat \cE \cong \cH$ 
can be identified with the orthogonal projection $P_\cH\res_{\cE_+}$: 
\[ q(f) = P_\cH f = \Gamma(f)(0).\] 
>From $\cH \subeq \cE_0 = \cE_+^\theta$ and the fact that $q\res_{\cE_0}$ is isometric, 
it now follows that $\cE_0 =  \cH.$ Further $q(U_t f) = \hat U_t q(f)$ yields 
\begin{equation}\label{eq:gamma-equiv}
P_\cH U_t f =  C_t P_\cH f \quad \mbox{ for } \quad f \in \cE_+, t > 0.
\end{equation} 

The elements of the space $\cE$ are represented 
by measurable functions $f\: \R \to \cH$. To get a better description of the subspace 
$\cE_+$, we would like to define, for a closed subset $F \subeq \R$, the subspace 
$\cE_F$ of those functions $f \in \cE$ whose inverse ``Fourier transform'' is supported in~$F$. 
This requires to define integrals 
$\int_\R \phi(x) f(x)\, dx \in \cH$ for Schwartz functions $\phi \in \cS(\R)$. 
Unfortunately, we do not see how to do this directly, but it can be done 
by first cutting off with respect to the spectral measure of~$A$. 
For $0 < a < b$, we write $\cH_{a,b} := P([a,b])\cH$, where $P$ is the spectral 
measure of~$A$ and $P_{a,b} = P([a,b]) \: \cH \to \cH_{a,b}$ for the orthogonal 
projection. For $v \in \cH_{a,b}$, we have 
\[ a \|v\|^2 \leq \la Av,v \ra \leq  b \|v\|^2,\] 
so that 
\begin{equation}
  \label{eq:esti-3}
\frac{a}{b^2 + x^2} \|v\|^2 
\leq \la A(A^2 + x^2)^{-1}v, v\ra 
\leq \frac{b}{a^2 + x^2} \|v\|^2.
\end{equation}
We conclude that the corresponding closed subspace 
\[ \cE_{a,b} := \{ f \in \cE \: f(\R) \subeq \cH_{a,b} \} \] 
can be identified as a topological vector space 
with the Hilbert space $L^2\Big(\R, \frac{dx}{\pi (1 + x^2)}\Big) \hat\otimes \cH_{a,b}$. 
For  $\phi \in \cS(\R)$ and $f \in \cE_{a,b}$ we now have 
\begin{align*}
 \int_\R |\phi(x)|\, \|f(x)\|\, dx 
&= \int_\R (b^2 + x^2)^\shalf |\phi(x)| \, (b^2 + x^2)^{-\shalf}\|f(x)\|\, dx \\
&\leq  \Big(\int_\R |\phi(x)|^2 (b^2 + x^2)\, dx\Big)^{\shalf} 
\Big(\int_\R (b^2 + x^2)^{-1}\|f(x)\|^2\, dx\Big)^{\shalf}.
\end{align*}
This implies that we have a continuous sesquilinear map 
\[ \beta \: \cE_{a,b} \times \cS(\R) \to \cH_{a,b}, \quad 
(f,\phi) \mapsto \Xi(f)(\phi) 
:= \int_\R \hat{\oline\phi}(x) f(x)\, dx
= \int_\R \oline{\tilde\phi}(x) f(x)\, dx,\] 
where $\tilde\phi(x) = \hat\phi(-x)$ is the inverse Fourier transform. 
Any $f \in \cE_{a,b}$ can thus be viewed as an $\cH_{a,b}$-valued tempered distribution 
on $\R$, so that, in the sense of \eqref{eq:2}, we have 
$\Xi(f) = \hat f$. 
We now define, for a closed subset $F \subeq \R$, the {\it spectral subspace}: 
\begin{equation}
  \label{eq:spectral}
\cE_F := \{ f \in \cE \: (\forall 0 < a < b<\infty )\,\, \supp(\Xi(P_{a,b}f)) \subeq -F \}.
\end{equation}

The main point of the following proposition is a characterization of the closed 
subspace $\cE_+$ of $\cE$ as the spectral subspace $\cE_{[0,\infty[}$. 

\begin{prop} \mlabel{prop:3.15a} 
  The spectral subspaces have the following properties: 
  \begin{description}
  \item[\rm(i)] The subspaces $\cE_F$ of $\cE$ are closed.  
  \item[\rm(ii)] $U_t \cE_F = \cE_{F + t}$ for $t\in \R$. 
  \item[\rm(iii)] $\cE_0 = \cE_{\{0\}}$ and $\cE_+ = \cE_{[0,\infty[}$. 
  \end{description}
\end{prop}

\begin{prf} (i) follows from the continuity of $\beta$ in both arguments. 

(ii) For $t \in \R$ and $(U_t f)(x) = e^{-itx} f(x)$, we have 
$\widehat{U_tf} (y) = \hat f(y + t),$ 
so that $\supp(\widehat{ U_t f} ) = \supp(f) - t$. 

(iii) Clearly, $\cE_0 \subeq \cE_{\{0\}}$, so that (ii) further implies 
$\cE_+ \subeq \cE^+ := \cE_{[0,\infty[}$. To see that we
have equality, we show that $\cE^+$ is $\theta$-positive with 
\begin{equation}
  \label{eq:eplus-rel}
\la \theta f, f \ra = \|P_\cH f\|^2 \quad \mbox{ for } \quad f \in \cE^+.   
\end{equation}
Then $\cE^+ = \cE_+$ follows from Proposition~\ref{prop:3.9} because 
$\cE_0$ is cyclic. This in turn leads to 
\[ \cE_{\{0\}} = \cE_{[0,\infty[} \cap \cE_{]-\infty,0]} 
= \cE_+ \cap \theta(\cE_+) = \cE_+^\theta = \cE_0.\] 

It therefore remains to verify \eqref{eq:eplus-rel}. 
Since the union of the subspaces 
$\cE_{a,b}^+ := \cE_{a,b} \cap \cE^+$ is dense in $\cE^+$, it suffices to verify this relation 
for $f \in \cE_{a,b}^+$. Approximating every element in this space by linear combinations 
of elements of the form $f = f_0 \cdot v$, $v \in \cH_{a,b}$, it suffices to verify 
\begin{equation}
  \label{eq:eplus-rel2}
\la f, \theta g \ra = \la P_\cH f, g \ra \quad \mbox{ for } \quad 
f = f_0 v, \, g = g_0 w \in \cE_{a,b}.
\end{equation}

Since $(U_t)_{t \geq 0}$ is a strongly continuous 
one-parameter group of isometries of $\cE^+$, for any $\delta$-sequence 
$\delta_n \in C^\infty_c(\R_+)$ and $f \in \cE^+$, we have $\hat{\delta_n} f \to f$ 
in $\cE^+$. As $d\sigma(x) := \frac{dx}{1+x^2}$ 
is a finite measure, $L^2(\R,\sigma) \subeq L^1(\R,\sigma)$. 
>From 
\[ \int_\R |\hat\delta_n(x)| \|f(x)\|\, dx 
=  \int_\R |\hat\delta_n(x)| (1 + x^2) \frac{\|f(x)\|}{1 + x^2}\, dx < \infty \] 
and 
\[ \int_\R |\hat\delta_n(x)|^2 \|f(x)\|^2\ dx 
=  \int_\R |\hat\delta_n(x)|^2 (1 + x^2) \frac{\|f(x)\|^2}{1 + x^2}\, dx < \infty \] 
it thus follows that  $\hat{\delta_n} \|f\| \in L^1(\R)  \cap L^2(\R)$. 
We may therefore assume, in addition, that $f_0,g_0 \in L^1(\R) \cap L^2(\R)$. 

Using the fact that $\phi(x) := \sqrt{\frac{2}{\pi}} \frac{\lambda}{\lambda^2 + x^2}$ 
has the Fourier transform $\hat \phi(x) = e^{-\lambda|x|}$, we obtain 
\[ \frac{\lambda}{\lambda^2 + x^2} = \frac{1}{2} \int_\R e^{-ixy} e^{-\lambda|y|}\, dy,\] 
and this in  turn leads to the relation 
\[ A(A^2 + x^2)^{-1} = \frac{1}{2} \int_\R e^{-ixy} e^{-|y|A}\, dy\] 
in $B(\cH)$ with respect to the strong operator topology. We thus obtain 
for $f = f_0 v$: 
\begin{align*}
 P_\cH f 
&= \frac{1}{\pi} \int_\R A(A^2 + x^2)^{-1} f(x)\, dx 
= \frac{1}{\pi} \int_\R f_0(x) A(A^2 + x^2)^{-1} v\, dx \\
&= \frac{1}{2\pi} \int_\R \int_\R f_0(x) e^{-ixy} e^{-|y|A}v\, dy\, dx 
= \frac{1}{\sqrt{2\pi}} \int_\R \hat f_0(y) e^{-|y|A}v\, dy
= \frac{1}{\sqrt{2\pi}} \int_0^\infty \hat f_0(-y) e^{-yA}v\, dy.
\end{align*} 
This further leads to 
\begin{align*}
 \la P_\cH f, P_\cH g \ra 
&= \frac{1}{2\pi} \int_0^\infty \int_0^\infty \hat f_0(-y_1) \oline{\hat g_0}(-y_2) 
\la e^{-y_1 A}v, e^{-y_2 A}w \ra\, dy \\ 
&= \frac{1}{2\pi} \int_0^\infty \int_0^\infty \hat f_0(-y_1) \oline{\hat g_0}(-y_2) 
\la e^{-(y_1 + y_2)A}v, w \ra\, dy \\ 
&= \frac{1}{\sqrt{2\pi}} \int_0^\infty (\hat f_0 * \oline{\hat g_0})(-y) \la e^{-yA}v, w \ra\, dy 
= \frac{1}{\sqrt{2\pi}} \int_0^\infty \widehat{\, f_0 g_0^*\, }(-y)\la e^{-yA}v, w \ra\, dy \\ 
&= \la P_\cH f_0 g_0^* v, w \ra 
= \frac{1}{\pi} \int_\R \la A(A^2 + x^2)^{-1} v, w \ra (f_0 g_0^*)(x) \, dx  \\
&= \frac{1}{\pi} \int_\R \la A(A^2 + x^2)^{-1} f_0(x) v, g_0(-x) w \ra  \, dx  
= \la f, \theta g \ra.
\end{align*} 
This complete the proof of \eqref{eq:eplus-rel2} and hence of the lemma. 
\end{prf}

\begin{ex}
For $\dim \cE_0 = 1$ and $C_t = e^{-t\lambda}$, $\lambda \geq 0$, the minimal 
dilation $\phi(t) = e^{-\lambda|t|}$ leads to the Hilbert space 
$\cE \cong L^2(\R, \frac{\lambda}{\pi}\frac{dx}{\lambda^2 + x^2})$ (cf.\ \cite[Ex.3.5(b)]{NO12}). 
\end{ex}

We finish this section with the following proposition, but delay one part of the proof until
Example \ref{ex:4.11}:

\begin{prop} \mlabel{prop:4.9} Let $(C_t = e^{-tA})_{t > 0}$ be a strongly continuous 
hermitian contraction semigroup on $\cH$ with $\ker A = \{0\}$. 
For the corresponding reflection positive one-parameter group 
$(U_t f)(x) := e^{-itx} f(x)$ on $\cE := L^2(\R,Q;\cH)$ we obtain an outgoing 
realization by composition of the 
inverse Fourier transform $\cF^{-1} \: L^2(\R,\cH) \to L^2(\R,\cH)$ with the unitary isomorphism 
\[ T \: \cE \to L^2(\R, \cH), \quad 
(Tf)(x) := \frac{1}{\sqrt \pi} \sqrt{A}(A + i x)^{-1} f(x). \] 
\end{prop}

\begin{prf} First we decompose $\cH$ into cyclic subspaces. In view of the Spectral 
Theorem, these are isomorphic to 
$(C_t f)(\lambda)  = e^{-t\lambda}f(\lambda)$ on $L^2([0,\infty[,\rho)$, where 
$\rho$ is a finite measure. Therefore the discussion in Example~\ref{ex:4.11} 
below applies to the cyclic case, and the general case is obtained by combining all that in a 
direct sum. 
\end{prf}

For the space $\cE = L^2(\R,Q;\cH)$ from above, the outgoing 
realization leads to an isomorphism with $L^2(\R,\cH)$ and hence to a natural 
concept of a spectral subspace 
\[ L^2(\R,\cH)_F = \{ f \in L^2(\R,\cH) \: \supp(\hat f) \subeq - F\}.\] 
These spectral subspaces do not coincide with the ones defined in \eqref{eq:spectral}. 
A simple example 
is $L^2(\R,\cH)_{\{0\}} = \{0\}$ because there is no non-zero square integrable 
polynomial, whereas $\cE_{\{0\}} = \cH$ (the constant functions) in the sense 
of \eqref{eq:spectral}.

Since cyclic one-parameter semigroups of contractions are given by multiplication 
with functions on $L^2$-spaces, we take a closer look at this special situation. 
We give here one example, a second example is discussed in the following subsection.

\begin{ex} Let $D \in C^{-\infty}(\R)$ be a reflection positive distribution and write 
$D = \hat\nu$ for a tempered measure $\nu$ on $\R$, so that the corresponding 
reflection positive unitary one-parameter group is given on  
$\cF := L^2(\R,\nu)$ by $(V_tf)(p) = e^{-itp} f(p)$. Here 
the constant function $1$ is a distribution vector 
(Corollary~\ref{cor:3.17}) and 
\[ \cF_+ = \oline{\Spann\{ \hat \phi \: \supp(\phi)  \subeq \R_+ \}}\] 
(cf.\ \cite[Prop.~2.12, Thm.~4.11]{NO12}). 
The corresponding space $\hat\cF$ is equivalent to 
$L^2([0,\infty[,\rho)$, where the Radon measure $\rho$ on $[0,\infty[$ is determined by 
$D \res_{\R_+} = \cL(\rho)$ (the Laplace transform), and 
\[ (\hat V_t f)(\lambda) = e^{-t\lambda} f(\lambda)\] 
(\cite[Prop.~2.12]{NO12}). In view of Lemma~\ref{lem:2.3}, 
the measure $\rho$ is tempered because the Laplace transform 
$\cL(\rho)$ extends to a distribution on $\R$.  
If $\lim_{x \to 0+} D(x) = \nu(\R)$ is finite, then $1 \in \cF_0$ is a cyclic vector; 
in particular $\cF_0 \not=\{0\}$. 

Suppose that $\rho(\{0\}) = 0$, i.e., $\hat \cF_{\rm fix} = \{0\}$. 
In view of Proposition~\ref{prop:e.5}, this is equivalent to 
$\cF_{\rm fix} = \{0\}$, which means that $\nu(\{0\}) = 0$. In
this case $\cF_+$ is an outgoing subspace according to
Proposition~\ref{prop:4.11}.
Furthermore, the condition $\int_1^\infty \frac{d\rho(y)}{y} < \infty$ 
is equivalent to the local integrability of $D\res_{\R^\times}$ 
(Lemma~\ref{lem:2.3}), and in this case $\nu = \Theta \, dx$ as 
in Lemma~\ref{lem:delta-fin}. Therefore $\cF_0 \not=\{0\}$ 
is equivalent to the finiteness of the measures $\nu$ and $\rho$. 
\end{ex}  

\subsection{Semigroups of multiplication operators} 
\mlabel{ex:4.5}  Let $\mu$ be a non-zero $\sigma$-finite measure on 
$(Y,\fS)$ and $h: Y\to \R_+$ be a measurable function. 
On  $\cH := L^2(Y, \mu)$ consider the contraction semigroup 
given by
\[(C_tf)(y) = e^{-th(y)} f(y)\, .\] 
Then dilation leads to the reflection positive unitary one-parameter group $(U_t)_{t \in \R}$ on 
\begin{equation}
  \label{eq:zeta0}
\cE \cong   L^2(X,\zeta) \quad \mbox{ for } \quad  
X := \R \times Y \quad \mbox{ and } \quad d\zeta(x,y) = \frac{1}{\pi}\frac{h(y)}{h(y)^2 + x^2}
\, dx \, d\mu(y)
\end{equation}
given by $(U_t F)(x,y) = e^{-itx} F(x,y)$, the involution 
$(\theta F)(x,y) = F(-x,y)$, and $\cE_0 = \cE_+^\theta 
\cong L^2(Y,\mu)$ (the functions independent of $x$) 
(cf.~\eqref{eq:scalprodincE}). 

Let $\pr \: \R \times Y \to \R, \pr(x,y) := x$ denote the projection onto the first variable and 
$\nu := \pr_*\zeta$. 
Then $\cF := L^2(\R, \nu)$ can be identified with the $U$-invariant subspace of $\cE$ consisting 
of functions independent of $y$. For $\cF_+ := \cF \cap \cE_+$ we obtain on 
$(\cF,\cF_+,\theta)$ a reflection positive one-parameter group $(V_t)_{t \in \R}$ 
by restriction. Now $\cF_0 = \cF_+^\theta = \cE_0 \cap \cF$ is the space of those 
functions $f(x,y) = f(x)$ whose Fourier transform is supported in $\{0\}$ 
(Proposition~\ref{prop:3.15a}(iii)), so that $f$ is a polynomial. 
To determine the space $\cF_0$, we take a closer look at the measure~$\nu$.

By construction, $d\nu(x) = \Theta(x)\, dx$, where 
\[ \Theta(x) 
:= \frac{1}{\pi}  \int_{Y} \frac{h(y)}{h(y)^2 + x^2}\, d\mu(y)  
= \frac{1}{\pi}  \int_{\R_+} \frac{y}{y^2 + x^2}\, d\rho(y) \] 
and $\rho := h_*\mu$ is the image of $\mu$ under $h$ on $\R_+$. 
If $\Theta$ is identically $\infty$, then $\cF = \{0\}$. 
If this is not the case, i.e., if $\rho$ is tempered and 
$\int_1^\infty \frac{1}{y}\, d\rho(y) <  \infty$, 
then Lemma~\ref{lem:delta-fin}(b) 
shows that every polynomial function in $L^2(\R,\nu)$ is constant. 
This implies that 
\[ \cF_0 = \C 1 \cap L^2(\R,\nu) 
= \{0\} 
\quad \Leftrightarrow \quad \nu(\R) = \zeta(X) =\mu(Y) =\rho(\R_+) =  \infty.\] 

If, in addition, $\rho(]0,1])$ is finite, 
then $\nu$ is tempered by Proposition~\ref{prop:3.15}, 
so that Lemma~\ref{lem:dist-vec} implies that the 
constant function $1$ is a distribution vector in $\cF^{-\infty}$. 
Remark~\ref{rem:d.2} further implies that 
the corresponding distribution cyclic representation 
is equivalent to the representation $(V_t)_{t \in \R}$ on $\cF$ and that 
$\hat \cF \cong L^2(\R,\rho)$ with $(\hat V_t f)(y) = e^{-ty} f(y)$ 
(see also \cite[Prop.~2.12]{NO12}). 

For any case where $\mu$, resp., $\rho$ is infinite, 
we thus obtain a distribution cyclic reflection positive one-parameter group for which 
$\cF_0$ is trivial. 

We also note that the constant function $1$ is a distribution 
vector of $\cE$ if and only if this holds for the representation 
in the subspace $\cF$, which is equivalent to the temperedness of 
$\nu$ (Corollary~\ref{cor:3.17}). 

If $1 \not\in \cE^{-\infty}$ and 
$F \: \R_+ \to \R$ is such that $\tilde\rho := |F|^2 \rho$ satisfies 
$\tilde\rho(]0,1]) < \infty$ and 
$\int_1^\infty \frac{d\tilde\rho(y)}{y} < \infty$, 
then the preceding discussion implies that $F(x,y) := F(y)$ 
defines a reflection positive distribution vector for which the 
corresponding distribution is given by 
the locally integrable function $\cL(\tilde\rho)(|x|)$. 

We have seen in Proposition~\ref{prop:4.11} that $\cE_+$ 
is outgoing for $U$, so that an outgoing realization of $(U,\cE)$ exists. 
Since, in general, these realizations are not so easy to obtain  explicitly, 
it is interesting to observe that this works in a natural way 
for the space $\cE = L^2(\R, Q;\cH)$. We complete the proof 
of Proposition~\ref{prop:4.9} in the following example. 

\begin{ex}[The case $Y=\R_+$ and $h(\lambda )=\lambda$]\mlabel{ex:4.11}
Let $\rho$ be a finite positive Radon measure on $\R_+$. 
Recall from Lemma~\ref{lem:delta-fin} the corresponding finite 
measure $d\nu(x) = \Theta(x)  \, dx$ on $\R$. 
We consider the Hilbert space 
\[\cE = L^2(\R\times \R_+,\zeta) \quad \mbox{ for } \quad 
d\zeta(x,\lambda) = \frac{\lambda}{\pi (\lambda^2+x^2)}\, dx d\rho(\lambda) \]  
from the discussion after Proposition~\ref{prop:exist-real} with 
$\cH = L^2(\R_+, \rho)$ and $AF(\lambda) = \lambda F(\lambda)$. 

The involution $\theta$ is given
by $(\theta f)(x,\lambda)=f(-x,\lambda)$. Let 
$\cE_0$ be the space of functions that are independent of $x$ and 
$\cE_+ = \cE_{[0,\infty[}$ the space of functions who's inverse Fourier transform is supported
on the positive half-line (cf.\ Proposition~\ref{prop:3.15}(iii)). 
Then $(U_tf)(\lambda, x)
=e^{-itx}f(x,\lambda)$ is a reflection positive one-parameter group 
on $\cE$ and $\hat\cE\simeq L^2(\R_+,\rho)$, where the isomorphism is induced by 
\[ q \: \cE_+ \to \cH, \quad 
q(f)(\lambda) :=  \frac{\lambda}{\pi} \int_\R \frac{f(x, \lambda)}{\lambda^2+x^2}\, dx \] 
(cf.\ \eqref{eq:PH}). 
%In the following we identify $\hat\cE$ with $L^2(\R,\nu)$. 
%\edz{Why? I do not see that $\nu$ plays a role in the following.}
Since $\cE_{\rm fix} = \{0\}$ and $\cE_0 \subeq \cE_+$ is $U$-cyclic, 
Proposition~\ref{prop:e.5}(c) implies that $\cE_+$ is outgoing. 

>From the factorization 
\[\frac{1}{\pi} \frac{\lambda}{\lambda^2+x^2}=\left(\frac{1}{\sqrt{\pi}}\frac{\sqrt{\lambda}}{\lambda - ix}\right)
\, \left(\frac{1}{\sqrt{\pi}} \frac{\sqrt{\lambda}}{\lambda + ix}\right)\]
we derive the unitary isomorphism
\[ T  :\cE \to L^2(\R, \cH), \quad 
T(f)(x,\lambda )=\frac{1}{\sqrt{\pi}} \frac{\sqrt{\lambda }}{\lambda +ix}
\, f(x,\lambda ).\] 
The operator $T$ clearly commutes with the multiplicative $\R$-action $U$ 
on $\cE$ and $L^2 (\R,\cH)$ from above. 

To avoid confusion, we will use the notation $f_\lambda(x) := f(x,\lambda)$
 and  write $\widetilde{f}(\xi ,\lambda )
=(\cF f_\lambda)(\xi )$ for 
the Fourier transform of $f$ in the first variable. 
We want to show that 
\[ T(\cE_+) = L^2(\R,\cH)_+ := \{ f \in L^2(\R,\cH) \:  x > 0 \Rarrow \tilde f(x,\cdot) = 0\}.\] 
A simple use of Cauchy's Integral Formula shows that
\begin{equation}\label{eq:4.8}
\psi_\lambda(t) := \cF^{-1}\left(\frac{1}{\sqrt{\pi}}
\frac{\sqrt{\lambda}}{\lambda + ix}\right)(t)=\sqrt{ 2 \lambda }\,  e^{-\lambda t}\chi_{[0,\infty[}(t) 
\quad \mbox{ for } \quad \lambda > 0.
\end{equation}
Note that 
$\psi_\lambda\in L^p([0,\infty [)\subset L^p(\R )$ for all $p\ge 1$ and 
$\lambda > 0$. 

>From $\tilde{Tf}_\lambda = \psi_\lambda^\vee * \tilde f_\lambda$ we derive that 
$T(\cE_+) \subeq L^2(\R,\cH)_+$. A simple calculation shows that
the involution 
$\tilde\theta := T \circ \theta \circ T^{-1}$ on $L^2(\R,\cH)$ is given by 
\[ (\tilde \theta f)(x,\lambda) =\frac{\lambda - ix}{\lambda +ix}f(-x,\lambda) 
=m_\lambda (x) f(-x,\lambda) 
\quad \mbox{ for } \quad 
m_\lambda(x) := \frac{\lambda - ix}{\lambda + ix}. \]
Since $\cE_+ \subeq \cE$ is maximal $\theta$-positive 
by Proposition~\ref{prop:3.9}(b), the subspace 
$T(\cE_+)$ is maximal $\tilde\theta$-positive in $L^2(\R,\cH)$. 
Therefore $T(\cE_+) = L^2(\R,\cH)_+$ will follow if we show that 
$L^2(\R,\cH)_+$ is $\tilde\theta$-positive. 

For $f,g \in L^2(\R,\cH)_+$, we have 
\begin{align*}
\la f,\tilde\theta f\ra 
&= \int_{\R_+} \int_\R \frac{\lambda + ix}{\lambda -ix} f_\lambda(x) f_\lambda^*(x)\, dx\, d\rho(\lambda).
\end{align*}
Therefore it is enough to recall from Subsection~\ref{ex:3.13} on the Hardy space of $\C_+$ 
that, for every fixed $\lambda > 0$, we have 
\[ \int_\R \frac{\lambda + ix}{\lambda -ix} f(x) f^*(x)\, dx 
\geq 0 \quad \mbox{ for } 
\quad f \in L^2(\R)_+.\]
The inverse Fourier transform $\cF^{-1} \: L^2(\R,\cH) \to L^2(\R,\cH)$ 
maps $L^2(\R,\cH)_+$ onto $L^2(\R_+,\cH)$ and intertwines the multiplication 
action $U$ with the translation action $(L_t f)(y) = f(y-t)$. 
Therefore $S := \cF^{-1} \circ T \: \cE \to L^2(\R,\cH)$ 
is an outgoing realization of $(\cE,U)$. 
\end{ex}

\section{Examples of reflection positive representations} 
\mlabel{sec:6} 

So far we have only discussed reflection positive one-parameter groups. 
But that is only a tool to handle  reflection positivity for more general 
symmetric Lie groups. Here we only discuss a few examples. The more
general theory of integration of infinitesimal representation
will be developed in the forthcoming article \cite{MNO14}. 

In this section we first define a suitable concept 
of a reflection positive unitary representations and then we discuss 
various classes of examples where the dilation construction for 
reflection positive one-parameter groups automatically 
yields reflection positive representations of larger Lie groups. 
To organize our examples, we first study covariance properties 
of the dilation process in Subsection~\ref{subsec:5.2},  
then we study the affine group of the real line 
(Subsection~\ref{subsec:5.3}) and the representations of the euclidean 
motion group of $\R^d$ associated to generalized free fields 
(Subsection~\ref{subsec:5.4}). In 
Subsection~\ref{subsec:5.5} we show that the conformally invariant 
among these representations form the complementary series 
representations of the conformal group 
$\OO_{1,d+1}(\R)_+$ of $\R^d$, resp., its conformal completion 
$\bS^d$. This observation builds a bridge to the prequel \cite{NO12} 
where the reflection positivity of these representations was studied in some 
detail. We conclude this section with a dilation construction 
of reflection positive representations of the $3$-dimensional Heisenberg 
group. 

\subsection{Reflection positive representations}
\mlabel{subsec:5.1} 

Let $(G,H,\tau)$ be a {\it symmetric Lie group}, i.e., 
$\tau$ is an involutive automorphism of $G$ and $H$ an open subgroup 
of $G^\tau$. Write $\g = \fh \oplus \fq = \g^\tau \oplus \g^{-\tau}$ 
for the eigenspace decomposition of $\g$ with respect to $\tau$. Define
$\g^c=\fh\oplus i\fq$. As $[\fh,\fh]\subset \fh$, $[\fh,\fq]\subseteq \fq$, and
$[\fq,\fq]\subseteq \fh$, it follows that $\fg^c$ is a Lie algebra. It is called
the \textit{Cartan dual} or simply the $c$-dual of $(\fg,\fh)$. Extending
$\tau$ to a complex linear Lie algebra automorphism of $\fg_\C$, also denoted $\tau$, 
and then restricting
to $\g^c$ shows that $(\g^c,\tau)$ is a symmetric pair. Let 
$G^c$ be the simply connected Lie group with Lie algebra $\g^c$. Then
$\tau$ defines an involution on $G^c$ and $H^c:=(G^c)^\tau$ is
connected (\cite[Th.~3.4]{Lo69}). 
The symmetric Lie group $(G^c,H^c,\tau)$ is again called 
the \textit{Cartan dual}, or simply the {\it $c$-dual}, of
$(G,H,\tau)$. 

\begin{ex} Let $d\in \N$ and $p,q\in \N_0$ such that $p+q=d$. 
Let 
\[I_{p,q}=\begin{pmatrix} -I_p & 0 \\ 0 & I_q\end{pmatrix}\, .\]
On $G=\R^d\rtimes \OO_d(\R)$ we consider the involution
\[\tau (x,a)=(I_{p,q}x,I_{p,q}aI_{p,q})\, .\]
Then 
\[\fq=(\R^p\oplus 0_q)\oplus \left\{\left. \begin{pmatrix} 0 & X\\ -X^\top & 0\end{pmatrix}
\, \right|\, X\in M_{p,q}\right\}\, .\]
It is then easy to see that
\[\fg^c\simeq \R^{p,q}\rtimes \fo_{p,q}(\R)\, .\]
Hence $G^c$ is locally isomorphic to $\R^{p,q}\rtimes \OO_{p,q}(\R)$. 
Here $\R^{p,q} \cong  i \R^p \oplus \R^q$ stands for
the vector space $\R^d$ with the bilinear form $\beta_{p,q}(x,y)=
x_1y_1+\ldots + x_py_p - x_{p+1}y_{p+1}-\ldots - x_ny_n$.
\end{ex}
 
\begin{defn} Let  $(G,H,\tau)$ be a symmetric Lie group and 
$(\cE,\cE_+,\theta)$ be a reflection positive Hilbert space. 
A  unitary representation $\pi \: G \to \cE$ is said to be {\it reflection 
positive} on $(\cE,\cE_+,\theta)$ if the following 
three conditions hold:
\begin{description}
\item[\rm(RP1)] $\pi(\tau(g)) = \theta \pi(g)\theta$ for every $g \in G$. 
\item[\rm(RP2)] $\pi (h)\cE_+= \cE_+$ for every $h \in H$. 
\item[\rm(RP3)] There exists a subspace 
$\cD\subeq \cE_+\cap \cE^\infty$ (where $\cE^\infty$ is the subspace of 
smooth vectors for $\pi$), dense in $\cE_+$, 
such that  $\dd \pi (x)\cD\subset \cD$ for all $x\in \fq$.  
\end{description}
\end{defn}
 
\begin{rem} (a) If 
$\pi$ is a reflection positive representation on $(\cE,\cE_+,\theta)$, 
then it follows from Proposition \ref{prop:osquant-commute} that
$\pi^c_H(h)=\hat{\pi (h)}$ defines a  unitary representation 
$(\pi^c_H, \hat\cE)$ of $H$.

(b) Note that for the symmetric Lie group $(\R,\{0\}, -\id_\R)$, the above 
concept of a reflection positive representation is weaker than the 
concept of a reflection positive one-parameter group used before. This is 
due to the fact that there are also interesting examples of reflection 
positive representations $\pi$ of $\R$, where $\cE_+$ is not invariant 
under $\pi(\R_+)$ (cf.\ Example~\ref{ex:1.3}(c)). 
\end{rem}

Let $(\pi, \cE)$ be a reflection positive representation 
of $(G,H,\tau)$ on $(\cE,\cE_+,\theta)$. 
If $x\in\fq$ satisfies  $\pi (\exp \R_+ x )\cE_+\subseteq \cE_+$ then we
define $\pi_x(t) := \pi(\exp tx)$. We then get
$\theta \pi_x(t) \theta = \pi_x(-t)$ and $\widehat{\pi_x}(t)$ is 
a hermitian contraction semigroup on $\hat \cE$. 
Therefore $\hat\pi_x$ has a self-adjoint positive generator $A_x$ determined by 
$\hat\pi_x(t) = e^{-tA_x}$ for $t \geq 0$. This is a
spectral condition which is not satisfied in all cases of interest. 
It is therefore too much to ask for in general, which lead us to condition 
(RP3). This condition implies that that 
$\hat\cD := q(\cD)$ is a dense subspace of $\hat\cE$ on which 
\[ \beta(x + i y) := \hat{\dd\pi(x)} - i \hat\dd\pi(y), \quad 
x \in \fh, y \in \fq \] 
defines an infinitesimally unitary representation of the Lie algebra 
$\g^c$. Now it is a natural problem to determine when 
this representation integrates a unitary representation $\pi^c$ 
of the corresponding simply connected Lie group $G^c$, compatible with 
$\pi_H^c$. In this is the case we call $(\pi, \cE)$ a 
{\it euclidean realization of the representation $(\pi^c,\hat\cE)$}. 
Systematic tools to achieve this task will be developed in the 
sequel \cite{MNO14}. In this section we shall encounter various examples 
where $\pi^c$ exists by more direct arguments and where $\pi$ and $\pi^c$ 
are connected by the dilation process studied in Section~\ref{sec:dil}. 

\subsection{Covariant one-parameter groups} 
\mlabel{subsec:5.2}

We consider a group of the form 
$G := \R \rtimes_\alpha H$, where 
$\alpha \: H \to \Aut(\R)_0 \cong \R^\times_+$ is a continuous 
homomorphism. Then $\tau_G(t,h) := (-t,h)$ defines on $G$ the structure of a 
symmetric Lie group and $S = \R_+ \rtimes_\alpha \R$ is an open subsemigroup 
invariant under $s \mapsto s^\sharp = \tau(s)^{-1}$. Note that 
$\g \to \g^c, (t, x) \mapsto (it,x)$ is an isomorphism of Lie algebras, so that 
we may put $G^c := G$ (cf.\ Subsection~\ref{subsec:5.1}). 

\begin{lem}\mlabel{lem:7.1} Let $(U_t)_{t \in  \R}$ be a reflection positive one-parameter group 
on $(\cE,\cE_+,\theta)$ and $\rho \: H \to \U(\cE)$ be a representation satisfying 
\[ \rho(h) U_t \rho(h)^{-1} = U_{\alpha_h t} 
\quad \mbox{ and } \quad \rho(h)\cE_+ = \cE_+ 
\quad \mbox { for } \quad t \in \R, h \in H.\] 
Then the following assertions hold: 
\begin{description}
\item[\rm(i)]  $\pi \: G  \to \U(\cE),  \pi(t,h) := U_t \rho(h)$ 
defines a reflection positive representation of $(G,S,\tau_G)$. 
\item[\rm(ii)] $\hat\pi(t,h) := \hat U_t \hat{\rho(h)}$ is the corresponding 
contraction representation of $(S,\sharp)$ on $\hat\cE$.
\item[\rm(iii)] If $\hat U_t = e^{-tA}$, then 
$\pi^c(t,h) := e^{itA}\hat{\rho(h)}$ defines a unitary representation of 
$G$ on $\hat \cE$. 
\end{description}
\end{lem}

\begin{prf} (i) follows immediately from our assumptions. 

(ii) The contraction representation $\hat\pi$ of $S$ on $\hat\cE$ 
is obtained by combining Proposition~\ref{prop:osquant-commute} with 
Proposition~\ref{prop:os-onepar}. 

(iii) From the fact that $\hat\pi$ is a representation of $S$, we derive that 
$\hat{\rho(h)} A \hat{\rho(h)}^{-1} = \alpha_h A$ for $h \in H$. 
Further, the semigroup $C_t = e^{-tA}$ 
is strongly continuous and extends by 
$C_z := e^{-zA}$ to a strongly continuous representation 
on the right half space 
$\{z\in \C\: \Re z\ge 0\}$ which is holomorphic in the interior. 
This in turn implies that $\pi^c$ defines a representation of $G^c = G$. 
\end{prf}

Next we verify the compatibility of the dilation 
construction in  Section \ref{sec:dil} with the $H$-action on~$\R$. 

\begin{lem} \mlabel{lem:covar} Let $(C_t)_{t > 0}$ be a strongly continuous 
contraction semigroup on $\cH$ and $(\rho, \cH)$ a unitary representation of $H$ satisfying 
  \begin{equation}
    \label{eq:equiv}
\rho_h C_t \rho_h^{-1} = C_{\alpha_ht} \quad \mbox{ for } \quad 
t > 0, h \in H.
  \end{equation}
Then the following assertions hold: 
\begin{description}
\item[\rm(i)] $(\pi(t,h)f)(x) := \rho(h) f(\alpha_h^{-1}(x+t))$ 
defines a reflection positive unitary representation of $G$ on the GNS Hilbert space 
$\cH_\phi\subeq \cH^\R$ for $\phi(t) := C_{|t|}$. 
\item[\rm(ii)] $(\pi(t,h)f)(x) := e^{-itx}  f(\alpha_hx)$ 
defines a reflection positive unitary representation of $G$ on $\cE 
= L^2(\R, \cH, Q)$ 
(cf.~Section~\ref{sec:dil}). 
\item[\rm(iii)] The unitary isomorphism 
$\Gamma \: \cE = L^2(\R, \cH, Q) \to \cH_\phi, \Gamma(f)(t) := P_\cH U_t f$, intertwines $\pi$ and $\tilde \pi$. 
\end{description}
\end{lem}

\begin{prf} (i) For $f \in \cH_\phi$, $h \in H$ and $t \in \R$  put 
$(\tilde V_h f)(x) := \rho(h) f(\alpha_h^{-1}x)$ and 
$(\tilde U_t f)(x) := f(x+t)$. 
We recall from Example~\ref{ex:vv-gns} that $\cH_\phi$ contains the 
total subset of functions of the form 
$\tilde U_t \phi v$, $t \in \R, v \in \cH$, satisfying 
\[ \la \tilde U_t \phi v, \tilde U_s \phi w\ra 
= \la \phi(t-s)v,w\ra.\] 
We have
\[ (\tilde V_h \tilde U_t \phi v)(x) 
= \rho(h) \phi\big(\alpha_h^{-1}x + t\big)v 
=  \phi\big(x + \alpha_ht\big)\rho(h) v  
=  (\tilde U_{\alpha_ht} \phi \rho(h) v)(x)
=  (\tilde U_{\alpha_ht} \tilde V_h \phi v)(x),  \] 
which implies already that $\tilde V$ and $\tilde U$ combine to a representation 
of $G$ on the space $C(\R,\cH)$. Unitarity of $\tilde V$ follows from 
\cite[Prop.~II.4.3]{Ne00} via $\phi(\alpha_h t) = \rho(h) \phi(t) \rho(h)^{-1}$, 
or directly from 
\begin{align*}
\la \tilde V_h \tilde U_t \phi v, \tilde U_s \phi w \ra 
&= \la \tilde U_{\alpha_ht} \phi \rho(h)v, \tilde U_s \phi w \ra
= \la \phi\big(\alpha_ht - s\big) \rho(h)v, w \ra\\
&= \la \rho(h) \phi\big(t - \alpha_h^{-1}s\big) v, w \ra
= \la \phi\big(t - \alpha_h^{-1}s\big) v, \rho(h)^{-1} w \ra\\
&= \la \tilde U_t \phi v,  \tilde U_{\alpha_h^{-1}s} \phi \rho(h)^{-1} w \ra
= \la \tilde U_t \phi v, \tilde V_h^{-1}\tilde U_s \phi w\ra.
\end{align*} 

To verify reflection positivity, we recall from Section~\ref{sec:dil} 
that $\cE_0 = \Spann \{ \phi v \: v \in \cH\}$ and 
$\cE_+$ is generated by $\tilde U_t \cE_0$, $t > 0$. 
Since $\tilde V$ preserves $\cE_0$ and normalizes $\tilde U$, 
it also preserves $\cE_+$. Therefore the representation 
$\pi$ of $G$ is reflection positive. 

(ii) Put $(V_h f)(x) := f(\alpha_hx)$. 
First we observe that $V_h U_t = U_{\alpha_ht} V_h$ for $h \in H, t \in \R$ and that 
$V_h$ preserves the subspace $\cE_0 \cong \cH$ of constant functions. This implies that 
$V_h \cE_+ = \cE_+$ for every $h \in H$. Hence (ii) follows from Lemma~\ref{lem:7.1}. 

(iii) From \eqref{eq:equiv} we obtain 
$\rho(h) A  \rho(h)^{-1} = \alpha_h A$ 
and therefore 
\[
\rho(h)^{-1} Q(\alpha_h^{-1}x)\rho(h)
= \alpha_h^{-1} A \big(\alpha_h^{-2} A^2 + \alpha_h^{-2} x^2\big)^{-1} = \alpha_h Q(x)\, .
\]
This implies that the projection 
\[ P_\cH \: \cE \to \cH, \quad 
P_\cH f =  \int_\R Q(x) f(x)\, dx \] 
satisfies 
\begin{align*} 
P_\cH V_h f 
&=  \int_\R Q(x) \rho(h)f(\alpha_hx)\, dx
= \alpha_h^{-1}  \int_\R Q(\alpha_h^{-1}x) \rho(h) f(x)\, dx
= \rho(h) \int_\R Q(x) f(x)\, dx
= \rho(h) P_\cH f.
\end{align*}
Therefore the unitary isomorphism $\Gamma$ 
satisfies 
\begin{align*}
\Gamma(V_h f)(t) 
&= P_\cH U_t V_h f 
= P_\cH V_h U_{\alpha_h^{-1}t} f 
=  \rho(h) P_\cH U_{\alpha_h^{-1}t} f 
= \rho(h) \Gamma(f)(\alpha_h^{-1}t).
\end{align*}
This means $\Gamma$ intertwines the representations on $\cE$ and 
$\cH_\phi$, respectively. 
\end{prf} 

\subsection{The $ax+b$-group} 
\mlabel{subsec:5.3}

A special example of the setting discuss above is the $ax+b$-group 
$G = \R \rtimes_\alpha \R^\times_+$, where 
$\alpha_a x = ax$ and $H = (\R^\times_+, \cdot) \cong \R$. 
The subset $S := \R_+ \rtimes \R_+^\times \subeq G$ is an open 
$\sharp$-invariant subsemigroup w.r.t.\ $\tau_G(x,a) = (-x,a)$.

For every $s > 0$, the measure 
$d\mu_s(y) := y^{s-1}\, dy$ on $[0,\infty[$ is tempered and has the Laplace transform 
\[ \cL(\mu_s)(x) = \frac{\Gamma(s)}{x^s} \quad \mbox{ for } \quad x > 0.\] 
 This function is locally integrable for $0 < s < 1$, so that, 
in this case, $D(x) := |x|^{-s}$ defines a reflection positive 
distribution (Theorem~\ref{thm:2.2}). 
Proposition~\ref{prop:a.9} further implies that, for every $s > 0$, 
$\cL(\mu_s)$ possesses an extension to a reflection positive distribution on~$\R$. 

On $\cH := L^2(\R_+, \mu_s)$ we consider the contraction semigroup 
\begin{equation}
  \label{eq:c-def}
(C_t f)(x) = e^{-tx}f(x)
\end{equation}
and the unitary representation of $\R^\times_+$, given by 
\[ (\rho_a f)(x) = a^{s/2} f(ax).\] 
It satisfies $\rho_a C_t \rho_a^{-1} = C_{at}$.
The corresponding unitary representation of $G^c \cong G$ on $\cH$ is given by  
\begin{equation}
  \label{eq:pics}
 (\pi_s^c(t,a)f)(x) = e^{itx} a^{s/2} f(ax)
\end{equation}
(cf.\ Lemma~\ref{lem:covar}). 

It is easy to see that the representations 
$\pi^c_s$ are irreducible and pairwise equivalent. Up to equivalence, 
it is the unique infinite-dimensional irreducible 
unitary representations of $G$ satisfying the positivity condition 
$-i\dd\pi^c_s(1,0) \geq 0$. This condition is necessary 
for the existence of a euclidean realization $(\pi, \cE)$ 
for which the one-parameter group $U_t := \pi(t,1)$ is reflection 
positive. 

We now discuss some natural euclidean realizations of the representations 
$(\pi_s^c,\cH)$ of $G^c$. 

\begin{ex}
The most natural source for reflection positive representations 
is the representation 
\[ (\pi(t,a)f)(x) := e^{-itx} \sqrt{a} f(ax)\]
of $G$ on $\cE := L^2(\R)$. It decomposes into the two irreducible 
subrepresentations on $L^2(\R_+)$ and $L^2(\R_-)$. 
For $\theta(f)(x) := f(-x)$ we now describe reflection positive 
distribution vectors which are semi-invariant for $H \cong \R^\times_+$. 
In view of Lemma~\ref{lem:dist-vec}, the distribution 
vectors for the one-parameter group $U_t := \pi(t,1)$ can be identified with 
measurable functions $h \: \R \to \C$ for which 
$(1 + x^2)^{-N} h$ is square integrable for some $N \in \N$. 
The requirement of $\theta$-invariance and $H$-semi-invariance implies that 
such a distribution vector is of the form 
\[ h_s(x) := |x|^{s-\shalf} \quad \mbox{ for some } \quad s \in \C.\] 
Local integrability of $|h_s(s)|^2 = |x|^{2 \Re s -1}$ is equivalent to $\Re s > 0$, and 
then $h_s$ represents a distribution vector. The corresponding distribution on 
$\R$ is given by 
\[ D_s(\phi) := \la h_s, \pi^{-\infty}(\phi)h_s \ra  
= \int_{\R} |h_s(x)|^2 \oline{\hat\phi(x)}\, dx 
= \int_{\R} |x|^{2s-1}\oline{\hat\phi(x)}\, dx.\] 
It is the Fourier transform of the function $|x|^{2s-1}$. 
For, $0 < s < \shalf$, it follows from \cite[Ex.~VII.7.13]{Schw73}
that it is given by
\[ \cF(|x|^{2s-1}) = \pi^{\frac{1-4s}{2}} 
\frac{\Gamma\big(\frac{1}{4} + s\big)}{\Gamma\big(\frac{1}{2}-s\big)} |x|^{-2s}.\]
Since the map $\C_+ \to \cS'(\R), s \mapsto |x|^{2s-1}$ is weakly holomorphic, 
it follows by analytic extension that the above formula describes 
the restriction to $\R_+$ for every $s > 0$. 
The holomorphic function $\Gamma\big(\frac{1}{2}-s\big)^{-1}$ has simple zeros 
in $\shalf + \N_0$, hence changes sign in these points. 
Therefore $\cF(|x|^{2s-1})$ is reflection positive if and only if 
\[ s \in \Big]0,\frac{1}{2}\Big] 
\cup \Big[\frac{3}{2}, \frac{5}{2} \Big] \cup \Big[\frac{7}{2}, \frac{9}{2}\Big] 
\cup \ldots\] 
Its restriction to $\R_+$ vanishes for $s \in \shalf + \N_0$. 
We thus obtain a complete description of the $H$-semi-invariant 
reflection positive distribution vectors in $L^2(\R)^{-\infty}$. 
\end{ex}

\begin{ex} \label{ex:5.7} For the dilation space of $(C,\cH = L^2(\R_+,\mu_s))$ 
(see \eqref{eq:c-def}), we obtain 
\[ \cE = L^2(\R \times \R_+, \zeta) \quad \mbox{ with } \quad 
 d\zeta(x,y) 
= \frac{1}{\pi}\frac{y}{x^2 + y^2}\, dx \frac{dy}{y^{1-s}} 
= \frac{1}{\pi}\frac{y^s}{x^2 + y^2}\, dx\, dy\]  
and the involution $\theta$ given by $(\theta f)(x,y) = f(-x,y)$ and 
$\cE_0 \cong \hat\cE \cong L^2(\R_+, \mu_s)$. 
By Lemma~\ref{lem:covar}, we have a canonical representation 
of $G$ on this space given by 
\[ (\pi(t,a) f)(x,y) = e^{-itx} a^{s/2} f(ax, ay).\] 
It is reflection positive with respect to 
$\cE_+=\cE_{[0,\infty[}$ and $(G,S,\tau_G)$. Hence $\pi$ is 
reflection positive and Lemma \ref{lem:7.1} shows that 
the corresponding representation of $G^c \cong G$ is equivalent to 
$\pi^c_s$. 

Let $\nu := \pr_*\zeta$ be the image of $\zeta$ under the projection $\pr(x,y) = x$ 
onto the first factor. It has the density 
\[ \Theta(x) = \frac{1}{\pi} \int_0^\infty \frac{y^s}{x^2 + y^2}\, dy \] 
which has finite values if and only if $s < 1$. If this is the case, 
then 
\begin{align*}
\Theta(\lambda x) 
&= \frac{1}{\pi} \int_0^\infty \frac{y^s}{\lambda^2 x^2 + y^2}\, dy 
= \frac{1}{\pi \lambda^2} \int_0^\infty \frac{y^s}{x^2 + (\lambda^{-1}y)^2}\, dy 
= \frac{1}{\pi \lambda} \int_0^\infty \frac{\lambda^s y^s}{x^2 + y^2}\, dy \\
&= \frac{\lambda^{s-1}}{\pi} \int_0^\infty \frac{y^s}{x^2 + y^2}\, dy 
= \lambda^{s-1}\Theta(x).
\end{align*}
Hence $\Theta(x) = c |x|^{s-1}$ for some $c> 0$, and this measure is tempered. 
Therefore $1$ is a distribution vector for the additive subgroup 
and the corresponding distribution is the Fourier transform of $\Theta$, 
which leads to the distribution $D_{s/2}$ from above. We conclude that, 
for $0 < s < 1$, the distribution vector $1$ generates a reflection 
positive subrepresentation $\cF \cong L^2(\R)$. 
\end{ex}

\begin{prop} \mlabel{prop:ereal-axb} Every unitary representation $(\pi^c,\cH)$ 
of $G$ with $-i\dd\pi^c(1,0) \geq 0$ has a euclidean realization. 
\end{prop}

\begin{prf} Since every unitary representation of $G$ is of type I  
and, up to equivalence, $\pi^c_1$ is the only irreducible representation  
which is non-trivial on $N := \R \rtimes \{1\}$ 
and satisfying $-i\dd\pi^c_1(1,0) \geq 0$, 
it follows that $\cH = \cH_0 \oplus \cH_1$, where 
$N$ acts trivially on $\cH_0$ and the representation on $\cH_1$ 
is a multiple $L^2(\R_+,\mu_1) \hat\otimes \cM$ of 
$\pi^c_1$ from \eqref{eq:pics} above. 
The representation on $\cH_0$ has the obvious euclidean realization 
with $\cE = \hat\cE = \cH_0$. 
Therefore it only remains to recall from 
Example~\ref{ex:5.7} that the representation $\pi^c_1$ on $L^2(\R_+,\mu_1)$ 
has a euclidean realization. 
\end{prf}

\begin{rem} In algebraic QFT, reflection positivity for the 
$ax+b$-group shows up in a very natural way. 
Here 
the action of the multiplicative group $\R^\times_+$ corresponds to the 
modular automorphism group $(\Delta^{it})_{t \in \R}$ of a von Neumann 
algebra $\cM \subeq B(\cH)$ and the action of the translation group 
is given by a unitary one-parameter group $(U_t)_{t \in \R}$ 
satisfying $U_t \cM U_t^* \subeq \cM$ for $t \geq 0$ 
(cf.\ \cite[Thm.II.9]{Bo92}). These situations are used in 
\cite{Bo92} to show that every theory of local observables in two 
dimensions which is covariant under translations can be embedded into a 
theory of local observables covariant under the whole Poincar\'e group. 
\end{rem}

\subsection{From Poincar\'e group to the euclidean group} 
\mlabel{subsec:5.4}

In this subsection we demonstrate that the dilation construction 
from Section~\ref{sec:dil} can be used to obtain euclidean realizations 
of the unitary representations of the Poincar\'e group corresponding 
to so-called generalized free fields, resp., to invariant measures on the 
future light cone. 

We start with a description of the Lorentz invariant measures 
on $\oline{V_+}$. 

\begin{defn} \mlabel{def:mum}
For $m \geq 0$ or $d > 1$, we define a Borel measure 
$\mu_m$ on 
\[ H_m := \{ p \in \R^d \: p_0^2 - \bp^2 = m^2, p_0 > 0 \} \subeq \oline{V_+} 
= \{ p= (p_0,\bp) \in \R^d \: p_0 \geq 0, [p,p] = p_0^2 - \bp^2\geq  0\} \] 
by 
\[ \int_{\R^d} f(p)\, d\mu_m(p) = 
\int_{\R^{d-1}} f(\sqrt{m^2+ \bp^2}, \bp)\frac{d\bp}{\sqrt{m^2 + \bp^2}}\]  
(cf.\ \cite[Ch.~IX]{RS75}, \cite[Lemma~9.1.2/3]{vD09}). 
These measures are invariant under the Lorentz group $L^\uparrow$ and every 
Lorentz invariant measure $\mu$  on $\oline{V_+}$ 
is of the form 
\begin{equation}
  \label{eq:nu-rho-a}
\mu = c \delta_0  + \int_0^\infty \mu_m\, d\rho(m),  
\end{equation}
where $c \geq 0$ and $\rho$ is a Borel measure on $[0,\infty[$ 
(with $\rho(\{0\}) = 0$ for $d = 1$) 
whose restriction to $\R_+$ is a Radon measure 
(Theorem~\ref{thm:5.11b}). 
\end{defn}

\begin{rem} For $d = 1$, we have $H_m = \{m\}$ for $m > 0$ and $H_0 = \eset$. 
Therefore $\mu_0$ does not make sense. For $m > 0$, we have 
$\mu_m = \frac{1}{m} \delta_m$, where $\delta_m$ is the Dirac measure in $m$. 

For $d =2$, the measure $\mu_0$ 
defined by 
\[ \int_{\R^2} f(p)\, d\mu_0(p) 
= \int_\R f(|\bp|, \bp) \frac{d\bp}{|\bp|}, \]
is singular in~$0$, but every 
$f \in \cS(\R^2)$ with $f(0)=0$ is integrable 
(cf.\ \cite[p.~103]{GJ81}). 
In particular, this measure does not define a distribution. 
\end{rem}

\begin{ex} \mlabel{ex:3.15} (Free fields) 
For the free scalar field of mass $m$ and spin $s = 0$ on $\R^d$ 
(with $m > 0$ or $d > 1$), the corresponding one-particle  Hilbert space 
is $\cH := L^2(\R^d, \mu_m)$ (cf.\ \cite[p.~103]{GJ81}). 
Here the time translation semigroup $C_t$ acts by the contractions 
\[ (C_tf)(p) = e^{-t p_0} f(p).\] 
In this case the dilation construction from Section~\ref{sec:dil} leads to the space 
$\cE = L^2(\R^{d+1},\zeta)$ and the projection 
$\nu = \pr_*\zeta$ to $\R^d$ under $\pr(p_0,x, \bp) = (p_0,\bp)$ 
is given by 
\[ d\nu_m(p_0, \bp) 
= \frac{1}{\pi} \frac{\sqrt{m^2 + \bp^2}}
{m^2 + \bp^2 + p_0^2}\, dp_0 \frac{1}{\sqrt{m^2 + \bp^2}}\, d\bp 
= \frac{1}{\pi} \frac{1} {m^2 + p^2}\, dp.\] 
Since all elements of $L^2(\R^d,\mu_m)$ can be represented by functions 
not depending on $p_0$, we obtain a unitary isomorphism 
\[ \pr^* \: L^2(\R^d,\nu_m) \to \cE.\] 
The measure $\nu_m$ is finite if and only if $d = 1$ and $m > 0$. 
It is tempered if and only if $d > 2$ or $m > 0$. 

Here the remarkable point is that the 
measure $\nu_m$ on $\R^d$ is rotation invariant, so that 
dilation with respect to the semigroup $C$ 
leads us directly from the irreducible representation $\pi^c$  
of the Poincar\'e group on $L^2(\R^d,\mu_m)$ 
to a representation of the euclidean motion group $\Mot(\R^d)$ on~$\cE 
\cong L^2(\R^d,\nu_m)$. 
\end{ex}

The preceding examples are the building blocks in a more general 
picture: 

\begin{ex} \mlabel{ex:3.16} (Generalized free fields) 
Let $\mu$ be a Lorentz invariant Radon measure as in 
\eqref{eq:nu-rho-a} on the forward light cone 
$\oline{V_+}\subeq \R^d$ with $c = \mu(\{0\}) = 0$. 
Then we have a natural unitary representation 
of the Poincar\'e group $G^c := \R^d \rtimes L^\uparrow$ on 
\[ \cH := L^2(\oline{V_+},\mu) \quad \mbox{ by } \quad (\pi^c(x,g)f)(p) := e^{-i xp} f(g^{-1}p).\] 
Analytic continuation of the time-translation group leads to the contraction semigroup 
\[ (C_t f)(p) =e^{-t p_0} f(p), \] 
and the dilation construction produces the Hilbert space 
\[ \cE = L^2(X,\zeta), \quad \mbox{ where } \quad 
 X := \R \times \oline{V_+} \quad \mbox{ and } \quad 
d\zeta(x,p) = \frac{1}{\pi}\frac{p_0}{p_0^2 + x^2}\, dx \, d\mu(p).\] 

We consider the unitary representation $\pi$ of $\R^d$ on $\cE$, given by 
\[ \big(\pi(x_0, \bx)f\big)(t,p_0, \bp) 
= e^{-i (x_0 t + \bx \bp)} f(t, p_0, \bp).\] 
%A test function $\psi \in C^\infty_c(\R^d)$ then acts on $\cE$ by 
%\[ (\pi(\psi) F)(t,p_0, \bp) 
%= \hat\psi(t, \bp) F(t,p_0,\bp).\] 
The constant function $1$ on $X$ is a distribution vector for $\pi$ 
if and only if the projected measure $\nu = \pr_* \zeta$, 
where $\pr(x,p_0, \bp) = (x,\bp)$, 
it tempered. Then the corresponding distribution is 
$D = \hat\nu$ (Lemma~\ref{lem:dist-vec}). 
To make this requirement more explicit, we use \eqref{eq:nu-rho-a} to obtain 
\[ \nu = \Theta \cdot dp \quad \mbox{ with } \quad 
\Theta(p) = \frac{1}{\pi} \int_0^\infty\frac{1}{m^2 + p^2}\, d\rho(m),\] 
where we have used Example~\ref{ex:3.15} to see that $\mu_m$ contributes the 
function $\Theta_m(p) = \frac{1}{\pi}\frac{1}{m^2 + p^2}$ to the density~$\Theta$.
The measure $\nu$ has a non-zero $L^2$-space if and only if 
$\Theta(p)$ is finite for non-zero $p$ 
(Lemma~\ref{lem:delta-fin}). 
Proposition~\ref{prop:3.15} now provides a characterization of the 
measures $\rho$ for which the measure $\nu$ is tempered. 
For $d > 2$, this is always the case if $\Theta$ is finite, and 
only for $d = 1,2$, additional conditions on the behavior of $\rho$ 
near $0$ are required. 

Let us assume that these conditions are satisfied, so that $\nu$ is tempered. 
Then the corresponding distribution $D = \hat\nu$ is reflection positive 
by Proposition~\ref{prop:repo-cm}. 
Let $\cF \subeq \cE$ denote the closure of the subspace 
$\pi(C^\infty_c(\R^d))1$ and $\cF_+\subeq \cF$ the closure of $\pi(C^\infty_c(\R^d_+))1$. 
Proposition~\ref{prop:3.15a}  now yields $\cF_+ \subeq \cE_+$, 
in particular, $\cF_+$ is $\theta$-positive, which provides a second 
argument for the reflection positivity of $\hat\nu$. 
We also obtain that 
\[ \cF_0 = \cF_+^\theta \subeq \cE_+^\theta = \cE_0.\] 
Therefore $\cF_0$ consists of function in $\cE_0$ that are independent of~$p_0$. 
This is precisely the $L^2$-space of the 
projected measure $\tilde\mu := \pr_*\mu$ on $\R^{d-1}$ for 
$\pr(x,\bp) = \bp$. Since $\mu$ coincides with the image of 
$\zeta$ under the projection $(x,p_0,\bp) \mapsto (p_0,\bp)$, we see that the 
image $\tilde\nu$ of $\nu$ under the projection $(x, \bp) \mapsto \bp$ 
coincides with $\tilde\mu$ (cf.\ Remark~\ref{rem:b.2}). 
According to Theorem~\ref{thm:5.11b}, 
$L^2(\R^{d-1}, \tilde\mu) 
= L^2(\R^{d-1}, \tilde\nu)$ is non-zero 
if and only if the measure $\tilde\mu$ is tempered. 
This happens if and only if $\mu$ is tempered, and, in addition, 
\[ \int_1^\infty \frac{1}{m}\, d\rho(m) < \infty.\] 
Then $1$ is a distribution vector for the representation of 
$\R^{d-1}$ on $\cF \subeq \cE$, and the corresponding cyclic subspace 
coincides with $\cF_0 \cong L^2(\R^{d-1}, \tilde\mu)$. 
\begin{center}
{{Table 1: Regularity properties of the measures} \\[3mm]
\begin{tabular}{|l |l |l|l|}
       \hline
 $\Theta < \infty$ &  $\nu$ \text{temp.}  & $\mu$ \text{temp.} & $\tilde\mu  
= \tilde\nu$ \text{temp.} \\
   \hline   \hline
  $\int_0^\infty \frac{d\rho(m)}{1 + m^2} < \infty$ 
& $\int_0^\infty \frac{d\rho(m)}{1 + m^2} < \infty$ 
& $\rho$ \text{temp.} 
& $\int_0^\infty \frac{d\rho(m)}{1 + m} < \infty$ 
 \\[1mm] 
\hline 
& \phantom{\Big(}$\int_0^1 \frac{d\rho(m)}{m} < \infty; d = 1$ 
& $\int_0^1 \frac{d\rho(m)}{m} < \infty; d = 1$ 
& 
 \\[1mm] 
\hline 
& \phantom{\Big(}$\int_0^1 \ln\left(\frac{1}{m}\right)\, d\rho(m)< \infty; d = 2$ 
& $\int_0^1 \ln\left(\frac{1}{m}\right)\, d\rho(m)< \infty; d = 2$ 
& 
 \\[1mm] 
\hline 
\end{tabular}}\\[3mm] 
\end{center}

Assume that $\tilde\mu$ is tempered. 
Then $\hat\cE \cong L^2(\R^d,\mu)$ contains the subspace 
$\hat\cE_0 \cong L^2(\R^{d-1}, \tilde\mu)$ of functions not depending on 
$p_0$ (where $p = (p_0,\bp) \in \R^d$ are as before) and the canonical map 
$\cE_0 \to \hat\cE_0$ is unitary. 
In the corresponding dilation space $L^2(\R^{d+1}, \zeta)$, 
$\cE_0\cong \hat\cE_0$ can be identified with the subspace of functions not depending 
on the first two coordinates. 
Accordingly, the ``time zero-subspace'' is the same on the euclidean and the Minkowski side. 

Since $1$ is a distribution vector for $\R^{d-1}$ generating 
$\cE_0$ and it is a distribution vector for $\R^d$ generating $\cE$, it follows that 
$\cE_0$ is a cyclic subspace of $\cE$. Accordingly, $\hat\cE_0$ is $\hat U$-cyclic in $\hat\cE$. 

We want to determine the corresponding positive definite operator-valued function 
\[ \phi \: \R \to B(\cE_0), \quad \phi(t) = P_0 U_t P_0^*,\] 
where $P_0 \: \cE \to \cE_0$ is the orthogonal projection. This function is determined by 
the relation 
\[ \la \phi(t)f, g \ra = \la U_t f, g \ra  \quad \mbox{ for } \quad f,g \in \cE_0.\]  
This leads to 
\begin{align*}
\la \phi(t)f, g \ra 
&= \int_{\R^d} e^{-it x_0} f(\bx) \oline{g(\bx)}\, \Theta(x_0,\bx)\, dx_0\, d\bx 
= \int_{\R^{d-1}} f(\bx) \oline{g(\bx)} \int_\R e^{-it x_0} \Theta(x_0,\bx)\, dx_0\, d\bx \\ 
&= \int_{\R^{d-1}} f(\bx) \oline{g(\bx)} \Theta_t(\bx)\, d\bx 
= \int_{\R^{d-1}} f(\bx) \oline{g(\bx)} \frac{\Theta_t(\bx)}{\Theta_0(\bx)}\, d\tilde\nu(\bx),
\end{align*}
where 
\begin{align*}
\Theta_t(\bx) 
&= \int_\R e^{-it x_0} \Theta(x_0,\bx)\, dx_0 
= \frac{1}{\pi} \int_0^\infty \int_\R e^{-it x_0} \frac{\lambda}{\lambda^2 + \bx^2 + x_0^2}\, 
\, dx_0\, d\rho(\lambda) \\
&= \int_0^\infty \frac{\lambda}{\sqrt{\lambda^2 + \bx^2}} e^{-|t|\sqrt{\lambda^2 + \bx^2}}\, d\rho(\lambda).
\end{align*}
Note also that 
\[ \Theta_t(\bx) \leq \Theta_0(\bx) 
= \int_0^\infty \frac{\lambda}{\sqrt{\lambda^2 + \bx^2}} \, d\rho(\lambda).\]
We conclude that $\phi(t) \in B(\cE_0)$ is given by multiplication with the function 
$\Theta_t/\Theta_0$. These functions are bounded with $0\le \Theta_t/\Theta_0 \le 1$.

For the subspace $\hat\cE_0 \subeq \hat\cE$ and $f,g \in \cE_0$, the relation 
\[ \la \hat U_t \hat f, \hat g\ra 
= \la \theta U_t f, g \ra 
= \la U_t f, g \ra = \la \phi(t)f,g \ra = \la \phi(t)\hat f, \hat g \ra \] 
implies that $\phi\res_{\R_+}$ is the positive  definite function on $\R_+$ 
corresponding to the cyclic subspace $\hat\cE_0 \subeq \hat\cE$. 
\end{ex}

\begin{ex} \mlabel{ex:5.14}
For the special case where $\rho =\delta_m$ and $m > 0$ or $d > 2$, we have
\[ \Theta_t(\bx) = \frac{m}{\sqrt{m^2 + \bx^2}} e^{-|t|\sqrt{m^2 + \bx^2}} 
\quad \mbox{ and } \quad 
\frac{\Theta_t(\bx)}{\Theta_0(\bx)} =e^{-|t|\sqrt{m^2 + \bx^2}} \] 
is multiplicative for $t \geq 0$. This corresponds to the fact that 
$q(\cE_0) = \hat\cE$ in this case, which in turn is due to the fact that 
the inclusion $L^2(\R^{d-1}, \tilde\mu_m) \into L^2(\R^d,\mu_m)$ is surjective. 

This has the interesting consequence that, if we consider elements of $\hat\cE$ as 
functions 
\[ f \:\R_+ \times \R^{d-1} \to \C \] 
as in the preceding example, we have 
\begin{equation}
  \label{eq:time0sem}
f(t,\bp) = (\hat U_t f)(0,\bp)  =e^{-t\sqrt{m^2 + \bp^2}} f(0,\bp).
\end{equation}
This in turn leads by analytic continuation to 
\begin{equation}
  \label{eq:time0grp}
 f(it,\bp) = (U^c_t f)(0,\bp)  =e^{it\sqrt{m^2 + \bp^2}} f(0,\bp).
\end{equation}
These formulas provide rather conceptual direct arguments for 
formulas like \cite[Prop.~6.2.5]{GJ81}. 
\end{ex}

\begin{rem} A unitary representation $(\pi, \cH)$ of the Poincar\'e group 
is said to be of {\it positive energy} if the spectrum of the time translation 
group is non-negative. In view of the covariance with respect to $L^\uparrow$, 
this is equivalent to the spectral measure of 
$\pi\res_{\R^d}$ to be supported in the closed forward light cone 
$\oline{V_+}$ because this is the set of all orbits of $L^\uparrow$ on which 
the function $p_0$ is non-negative. 

If such a representation is multiplicity free on $\R^d$, then
$\cH = L^2(\oline{V_+},\mu)$ for a measure $\mu$ on $\oline{V_+}$ 
which is quasi-invariant under $L^\uparrow$. Since the action of $L^\uparrow$ 
on $\oline{V_+}$ has a measurable cross section and every orbit 
carries an invariant measure, the measure $\mu$ can be chosen 
$L^\uparrow$-invariant. The representation $\pi$ is irreducible 
if and only if the measure $\mu$ is ergodic, i.e., 
$\mu = \mu_m$ for some $m \geq 0$ (with $m > 0$ for $d = 1$) or 
$\mu = \delta_0$ (the Dirac measure in $0$). 
\end{rem}

For all the multiplicity free representations 
$(\pi^c, L^2(\oline{V_+}), \mu)$, Example~\ref{ex:3.16} 
provides a euclidean realization in the dilation space
$\cE = L^2(X,\zeta)$, as far as the representation of the 
subgroup \break $\R^d \rtimes \OO_{d-1}(\R)$ is concerned. 
Note that the subspace $\cE_0 \subeq \cE$ is invariant under the subgroup 
$(G^c)^\tau \cong \R^{d-1} \rtimes \OO_{d-1}(\R)$, which also implies the 
invariance of $\cE_+$ under this group. 

A euclidean realization for the full group is obtained 
in Example~\ref{ex:3.15} for irreducible representations, 
i.e., $\mu = \mu_m$. In the general case we assume that 
$\nu$ is tempered. Then the following theorem 
is the bridge between  the reflection  positive 
representation of $\Mot(\R^d)$ on $\cF \cong L^2(\R^d,\nu)$ 
and the representation 
$\pi^c$ of the Poincar\'e group on $\hat\cF \cong L^2(\R^d, \mu)$. 

\begin{thm} \mlabel{thm:6.11} 
If $\nu$ is tempered, then $1 \in \cE^{-\infty}$ is a 
reflection positive distribution vector for the action $\pi$ of $\R^d$.
Accordingly, we obtain a reflection positive representation 
of $\R^d$ on the subspace $\cF \subeq \cE$ generated by $\pi^{-\infty}(C^\infty_c(\R^d))1$. 
The corresponding reflection positive distribution $\hat\nu$ on $\R^d$ is 
rotation invariant, so that $\cF$ 
carries a reflection positive representation of $\Mot(\R^d)$ for which 
$\cF_0$ and $\cF_+$ are invariant under $H := \Mot(\R^d)^\tau \cong 
\R^{d-1} \rtimes \OO_{d-1}(\R)$. 

Moreover, $\hat\cF \cong L^2(\oline{V_+},\mu)$, 
$q \: \cF_+ \to \hat\cF$ is $H$-equivariant and 
$\hat{\pi(t,0)} = \pi^c(t,{\bf 0},\1)$ for the natural representation 
$\pi^c$ of the Poincar\'e group $\R^d \rtimes L^\uparrow$ on $\hat\cF$. 
\end{thm}

\begin{prf} We have already seen that $1 \in \cE^{-\infty}$ is equivalent to 
$\nu$ being tempered. To determine the corresponding space 
$\hat \cF$, we have to take a closer look at the corresponding 
reflection positive distribution $D = \hat\nu$ for 
$(\R^d, \R^d_+, \theta)$. In view of \cite[Prop.~2.12]{NO12}, 
we have to write $D\res_{\R^d_+}$ as a Fourier--Laplace 
transform $\cF\cL(\gamma)$ of a measure $\gamma$ on 
$[0,\infty[ \times \R^{d-1}$ to obtain $\hat\cF \cong L^2(\R^d, \gamma)$. 

This is done by verifying $\cF\cL(\mu) = \hat\nu\res_{\R^d_+}$. 
First we observe that the temperedness of $\mu$ implies that 
\[ \cF\cL(\mu)(x) := \int_{\oline{V_+}} 
e^{-x_0 p_0} e^{i \bx \bp}\, d\mu(p) \] 
exists pointwise and defines an analytic function on $\R^d_+$. 
Here the main point is that, on $\oline{V_+}$ we have 
$\|p\|^2 = p_0^2 + \bp^2 \leq 2 p_0^2$ (cf.\ \cite[Ex.~4.12]{NO12}). 
We have  
\begin{align*}
 \cF\cL(\mu)(x) 
&= \int_{\oline{V_+}} e^{-x_0 p_0} e^{i \bx \bp}\, d\mu(p) 
= \int_0^\infty \int_{\R^{d-1}} e^{-x_0 p_0} e^{i \bx \bp}\, d\mu(p_0, \bp) \\
&= \int_0^\infty \int_{\R^{d-1}} \frac{1}{\pi} \int_\R 
e^{itx_0}\frac{p_0}{p_0^2 + t^2}\, dt\, e^{i \bx \bp}\, d\mu(p_0, \bp) 
= \int_{\R \times \R^d}  e^{i(tx_0 + \bx \bp)}\, d\zeta(t, p_0, \bp) \\
&= \int_{\R^d}  e^{i(t x_0 + \bx \bp)}\, d\nu(t, \bp)  = \hat\nu(x).
\end{align*}
If $\mu$ is infinite, then the triple integral only exists as 
an iterated integral in the correct order, not in the sense that 
the integrand is Lebesgue integrable. One can repair this problem 
by integrating against a test function on $\R^d_+$, 
and then the above calculation show that $\cF\cL(\mu)$ coincides with 
$\hat\nu$ on $\R^d_+$ as a distribution. 
\end{prf}

\subsection{The conformally invariant case} 
\mlabel{subsec:5.5}

In this section we study the special case where the measure
$\mu$ on $\oline{V_+}$ is semi-invariant under homotheties. This turns out to provide 
a bridge to the complementary series representations of 
$\OO_{1,d+1}(\R)_+$ studied in \cite{NO12} because the representation of 
$\Mot(\R^d)$ on $L^2(\R^d,\nu)$ extends to the conformal group 
(cf.\ Theorem~\ref{thm:6.14}). 

\begin{lem} \mlabel{lem:7.4.1} An $L^\uparrow$-invariant measure 
$\mu = \int_0^\infty \mu_m\, d\rho(m)$ on $\oline{V_+}$ is semi-invariant 
under homotheties if and only if 
$\rho(m) = m^{s-1}\, dm$ for some $s \in \R$. 

If this is the case, then 
$\rho$ is tempered if and only if $s > 0$ and 
$\mu$ is tempered if, for $d = 1$, we have $s > 1$. 

For $d > 1$, the measure $\tilde\mu$ on $\R^{d-1}$ is tempered if and 
only if $s > 1$. For $d = 1$ the measure $\mu$ is never finite. 
\end{lem}

\begin{prf} Let $h_a(p) := ap$ for $a \in \R_+$. We first claim that 
\[ (h_a)_* \mu_m = a^{2-d}\mu_{am} \quad \mbox{ for } \quad a > 0, m \geq 0.\] 
This follows  from 
\begin{align*}
& \int_{\R^d} f( p )  \,  d((h_a)_* \mu_m)(p)
= \int_{\R^d} f(ap) \, d\mu_m(p)
= \int_{\R^{d-1}} f(a \sqrt{m^2 + \bp^2}, a \bp)\frac{ d\bp}{\sqrt{m^2 + \bp^2}} \\
&= \int_{\R^{d-1}} f(\sqrt{a^2m^2 + (a \bp)^2}, a\bp)\frac{ d\bp}{\sqrt{m^2 + \bp^2}} 
=a^{1-d} \int_{\R^{d-1}} f(\sqrt{a^2m^2 + \bp^2}, \bp)\frac{ d\bp}{\sqrt{m^2 + \bp^2/a^2}} \\
&= a^{2-d}\int_{\R^{d-1}} f(\sqrt{a^2m^2 + \bp^2}, \bp)\frac{ d\bp}{\sqrt{(am)^2 + 
\bp^2}} 
= a^{2-d} \int_{\R^d} f(p)\, d\mu_{am}(p).
\end{align*}
This leads to 
\[ (h_a)_* \mu = \int_0^\infty \mu_{am}a^{2-d}\, d\rho(m) 
= \int_0^\infty \mu_{m}a^{2-d}\, d((h_a)_*\rho)(m) 
= a^{2-d} \int_0^\infty \mu_{m}\, d((h_a)_*\rho)(m).\] 
Since $\rho$ is determined uniquely by $\mu$, it follows that $\mu$ 
is semi-invariant under the maps $h_a$ if and only if $\rho$ has the corresponding 
property, and this means that $d\rho(m) = m^{s-1}\, dm$ for some $s \in \R$. 
Then $(h_a)_*\rho = a^{-s} \rho$, so that 
$(h_a)_*\mu = a^{2-d-s}\mu$.

According to Theorem~\ref{thm:5.11b}(a), $\mu$ is tempered if and only if 
$\rho$ is tempered and satisfies, for $d = 1,2$, the additional conditions (C2/3). 
Clearly, $\rho$ is tempered if and only if $s > 0$. 
For $d = 1$, we find the additional condition
$\int_0^1 m^{s-2}\, dm < \infty,$ which is equivalent to $s > 1$. 
For $d = 2$, condition (C3) reads 
\[ \infty > \int_0^1 \ln(m^{-1}) m^{s-1}\, dm 
= \int_0^\infty \ln(e^x) e^{(1-s)x} e^{-x}\, dx 
= \int_0^\infty x e^{-sx}\, dx, \] 
which is satisfied for $s > 0$. 

For $d > 1$, the corresponding measure $\tilde\mu$ on $\R^{d-1}$ 
is tempered if and only if $\mu$ is tempered and  
$\int_1^\infty \frac{1}{m} m^{s-1}\, dm = \int_1^\infty m^{s-2}\, dm$ is finite 
(Theorem~\ref{thm:5.11b}(b)), which is equivalent to $s < 1$.

For $d = 1$, we have $\mu_m = m^{-1} \delta_m$, so that 
$d\mu(p) = \int_0^\infty m^{s-2}\, dm$, and this measure is never finite. 
\end{prf}

>From now on we write $d\rho_s(m) = m^{s-1}\, dm$. 
The measure $\mu$ is also semi-invariant 
under homotheties, so that we can expect the corresponding
representation of the Poincar\'e group to extend to the conformal group 
$\SO_{2,d}(\R)$ of Minkowski space. 

\begin{lem} \mlabel{lem:6.13} 
The measure $\nu = \Theta \cdot dp$ is tempered if and only if 
$0 < s < 2$ for $d > 1$ and, if $0 < s < 1$ for $d = 1$. In this case 
$\Theta$ is a multiple of $\|p\|^{s-2}$ and the Fourier transform 
$\hat\nu$ is a positive multiple of $\|x\|^{-d+2-s}$. 
\end{lem}

\begin{prf} 
The integral $\int_0^\infty \frac{m^{s-1}}{1 + m^2}\, dm$ is finite if and only if 
$0 < s < 2$. In this case 
$\Theta_s(p) := \int_0^\infty \frac{m^{s-1}}{m^2 + p^2}\, dm$ is finite and 
a direct calculation shows that it is a positive multiple of $\|p\|^{s-2}$. 

The function $\|p\|^{s-2}$ on $\R^d$ is locally integrable if and only if 
$s > 2 - d$ and then $\nu$ is a tempered measure. 
Only for $d = 1$, we obtain the additional restriction $s >1$. 

For $0 < s < d$, it follows from \cite[Ex.~VII.7.13]{Schw73} that  
\[ \cF(\|x\|^{-s}) = \pi^{s-d/2}
\frac{\Gamma\big(\frac{d-s}{2}\big)}{\Gamma\big(\frac{s}{2}\big)} \|x\|^{s-d}.\]
We conclude that, for $0 < s < 2$, the Fourier transform $D_s$ of the measure
$d\nu(p) \sim \frac{dp}{\|p\|^{2-s}}$ is a positive multiple 
of $\frac{1}{\|x\|^{d-2+s}}$. 
\end{prf}

The preceding lemma implies in particular that the distribution 
$\frac{1}{\|x\|^a}$ is reflection  positive for $d - 2 < a < d$, 
which has been obtained in \cite[Prop.~6.1]{NO12} by other means. 
This connection is made more precise in the following theorem:

\begin{thm} \mlabel{thm:6.14} For $0 < s < 2$, resp., $1 < s < 2$ for $d =1$, 
the following assertions hold: 
  \begin{description}
  \item[\rm(i)] The canonical representation 
of the conformal motion group 
$\R^d \rtimes (\OO_d(\R) \times \R^\times_+)$ 
on $\cE := L^2(\R^d, \nu) \cong \cH_D$ for 
$D(x) = \|x\|^{-d+2-s}$ 
extends to a complementary series representation of the orthochronous 
euclidean conformal $\OO_{1,d+1}(\R)_+$. 
  \item[\rm(ii)] The corresponding representation 
of the conformal Poincar\'e group 
$\R^d \rtimes (L^{\up} \times \R^\times_+)$ is irreducible and 
extends to a representation of a covering of 
the relativistic conformal group $\SO_{2,d}(\R)_0$.   
  \end{description}
\end{thm}

\begin{prf} (i) From Lemma~\ref{lem:6.13} 
we know that $\cE := L^2(\R^d,\nu)$ can be identified 
naturally with the Hilbert space $\cH_D$ obtained by completion of 
$C^\infty_c(\R^d)$ with respect to the scalar product 
\[ \la \phi, \psi \ra_s := \int_{\R^d} \int_{\R^d}
\frac{\phi(x) \oline{\psi(y)}}{\|x-y\|^{d-2+s}}\, dx\ dy.\] 
Now \cite[Prop.~6.1]{NO12} (and the proof of Lemma~5.5 in \cite{NO12}) implies 
that the representation 
of $\Mot(\R^d)$ in this space extends to an irreducible complementary 
series representation of the conformal group 
$\OO_{1,d+1}(\R)_+$. 

(ii) The irreducibility of the representation $\pi^c$ follows from the 
transitivity of the action of $\R^\times_+ L^\uparrow$ on the open light cone $V_+$. 
To see that this representation extends to $\SO_{2,d}(\R)_0$, 
we can use the fact that the representation 
$\pi$ of the conformal group $G$ is reflection positive 
with respect to the open subsemigroup 
of strict compressions of the open half space $\R^d_+$ in the conformal 
compactification $\bS^d$. 
As explained in \cite[\S\S 6,10]{JOl00}, see also \cite{HN93}, \cite{JOl98}, 
the reflection positivity and the L\"uscher--Mack 
Theorem now provides an irreducible  representation of the 
simply connected $c$-dual group $G^c$ on~$\hat\cE$. 
\end{prf}

\subsection{The Heisenberg group} 
\mlabel{subsec:heisenberg}

In this subsection we apply the dilation process to the Heisenberg group.
Let 
$G = \Heis(\R^2,\omega)$ 
be the $3$-dimensional Heisenberg group. We fix basis vectors 
$P,Q,Z$ of its Lie algebra $\g$ satisfying the relations 
\[ [P,Q] = Z, \quad [P,Z] = [Q,Z] = 0.\]  
We  consider the involution $\tau$ on $\g$ defined by 
\[ \tau(Q) = Q, \quad \tau(P) = -P \quad \mbox{ and } \quad 
\tau(Z)= - Z.\] 
The map $Q\mapsto Q$, $P\mapsto iP$ and $Z\mapsto iZ$ defines a Lie algebra
isomorphism $\g \to \g^c$. We therefore identify 
$G$ and $G^c$ in the following.

Denote by $\pi= \pi_1$ the Schr\"odinger representation of $G$ acting
on $\cH := L^2(\R)$. It is determined by
\[ (\dd\pi_1 (Q)f)(x) = i x f(x), \quad 
(\dd\pi_1(P) f)(x) =  f'(x), \quad 
(\dd\pi_1 (Z)f)(x) = i  f(x).\] 

For $\lambda \in \R^\times$, we define the automorphisms 
$\gamma_\lambda \in \Aut(G) \cong \Aut(\g)$ of $G$ specified on $\g$ by 
\[\gamma_\lambda(Q) := Q, \quad 
\gamma_\lambda(P) := \lambda P, \quad 
\gamma_\lambda(Z) := \lambda Z\] 
and obtain a twisted 
Schr\"odinger representation $\pi_\lambda := \pi_1\circ \gamma_\lambda$ 
on $L^2(\R)$, where 
\[ (\dd\pi_\lambda(Q)f)(x) = i x f(x), \quad 
(\dd\pi_\lambda(P)f)(x) = \lambda f'(x), \quad 
(\dd\pi_\lambda(Z)f)(x) = i \lambda f(x).\] 
We also obtain for 
$\lambda = 0$ a representation 
$\pi_0$ of $G$ on $L^2(\R)$ for which 
$\dd\pi_0(P) = \dd\pi_0(Z) = 0$, so that it is not irreducible. 

Consider the contraction semigroup $C_t f= e^{-t}f$ on $\cH = L^2(\R)$. 
We thus obtain the dilation Hilbert space 
\[ \cE = L^2(\R, \mu;\cH) \cong L^2(\R,\mu) \hat\otimes \cH, \quad \mbox{ where } \quad 
d\mu(x) =\frac{1}{\pi} \frac{dx}{1 + x^2}.\] 
Then $\cE_0\simeq\cH$, the space of constant
functions, and $\cE_+=\cE_{[0,\infty[}$ (Proposition~\ref{prop:3.15a}). 
Define a representation $\pi$ of $G$ on $\cE$ by
\[ (\pi(g)f)(x) = \pi_{-x}(g) f(x), \] 
so that 
\[ (\pi(\exp tZ)f)(x) = e^{-itx} f(x) \quad \mbox{ and } \quad 
(\pi(\exp t Q)f)(x,y) = e^{ i t y} f(x,y). \]
The canonical subspace $\cE_+$ is $\theta$-positive and invariant 
under the semigroup $\pi(\exp \R_+ Z)$, actually it is generated by 
$\pi(\exp tZ) \cE_0$, $t > 0$. As $Q$ commutes with $Z$, 
the invariance of $\cE_0$ under $G^\tau = \exp(\R Q)$ further implies the 
invariance of $\cE_+$ under $G^\tau$. 
%Together these two families form a unitary representation 
%of the abelian subgroup $\exp (\R Q + \R Z)\subset G $ 
%which is isomorphic to a multiple of the restriction of $\pi_1 $ to this group.
% KH: Not true because the action of $Z$ is different!

On the space $\hat\cE \cong L^2(\R)$ we have the unitary representation 
$\pi^c := \pi_1$ of $G^c \cong G$. This representation is compatible 
with $\pi(\exp tZ)\,\hat{} = C_t = e^{-t}$ for $t \geq 0$, and the map 
$q \: \cE_+ \to \hat\cE$ is also equivariant with respect to $G^\tau$ 
which commutes with $\theta$. 

For $X := pP + z Z \in \g^{-\tau}$ with $p \not=0$, we have 
$\Spec(-i\dd\pi^c(X)) = \R$, so that we cannot expect that 
$\pi(\exp(t X))$ is reflection positive with respect to $(\cE,\cE_+,\theta)$. 
Therefore $\pi$ is not reflection positive with respect to $(G,S,\tau)$ for an 
open subsemigroup $S$ of $G$. 
However, we can exhibit a dense subspace $\cS_+ \subeq \cE_+$, invariant under  
$\dd\pi(\g)$, so that 
\begin{equation}
  \label{eq:heis-inter}
 q \circ (\dd\pi(X) - i \dd\pi(Y))\res_{\cS_+} 
= \dd\pi^c(X + i Y)\res_{q(\cS_+)} \quad \mbox{ for } 
\quad \tau(X) = X, \tau(Y) = -Y.
\end{equation}

Before considering the the action of the one-parameter group $\exp (\R P)$, 
let us prove the following lemma that will again be used in the next section.

\begin{lem}\label{le:8.5.1} 
Let 
\[ \cS(\R)_+ := \{ f \in \cS(\R) \: \hat f\res_{\R_+} = 0\} \quad \mbox{ and } \quad 
\cS_+ := \cS(\R)_+ \otimes \cS(\R) \subeq \cE_+.\] 
Then $\cS_+$ is dense in $\cE_+$. Let $(Tf)(x,y) := ixf (x,y)$ for $f\in\cS_+$.
Then $P_0(Tf)=P_0(f)$ for all $f\in\cS_+$ and the orthogonal 
projection $P_0 \: \cE \to \cE_0 \cong \cH$.
\end{lem}

\begin{prf} It is clear that $\cS_+$ is dense in $\cE_+$ and that $T$ 
preserves $\cS_+$. From \eqref{eq:PH} in Section~\ref{sec:dil}, we recall that 
\[ (P_0 f)(y)=\frac{1}{\pi}\int_{\R } \frac{1}{1+x^2}\, f(x,y)\, dx \quad 
\mbox{ for } \quad f\in \cE.\]
Hence the following holds for $f\in \cS_+$:
\begin{eqnarray*}
P_0( T f)(y) & =&
\frac{1} {\pi } \int_{\R}\frac{i x}{1+x^2}   f (x,y)\, dx
 =
\frac{1} {\pi } \int_{\R}\frac{1}{1+x^2}\, f(x,y)\, dx- \frac{1}{\pi}\int_{\R} \frac{1 - i x}{1+x^2}   f (x,y)\, dx\\
&=&   (P_0 f) (y ) -\frac{1}{\pi }  \int_{\R} \frac{1}{1+i x}  f  (x,y)\, dx \, .
\end{eqnarray*}
For any $\phi \in \cS(\R)_+$, the Fourier transform of 
$\Phi(x) := \frac{\phi(x)}{1 + ix}$ is the convolution of $\hat\phi$ with 
the Fourier transform of 
\[ (1 + ix)^{-1} 
= \int_0^\infty e^{-(1 + ix)t}\, dt 
= \int_0^\infty e^{-itxt} e^{-t}\, dt \] 
which vanishes on $\R_+$. Therefore $\hat\Phi$ vanishes on $\R_+$, which in 
particular leads to $0 = \hat\Phi(0) = \int_\R \Phi(x)\, dx$. 
This proves our assertion. 
\end{prf}

The one-parameter group $\exp (\R P)$ acts on $\cE$ by
\[ (\pi(\exp t P)f)(x,y)=f(x,y- tx)\, .\]
Hence $\cE_+$ is not invariant under $\exp(\R_+ P)$ or $\exp(\R_- P)$, 
but for $f\in\cS_+$, it follows
from Lemma \ref{le:8.5.1} that the infinitesimal generator 
\[ (\dd \pi ( P)f) (x,y)= - x \frac{\partial f }{\partial y}(x,y) \]
maps $\cS_+$ into itself and by Lemma \ref{le:8.5.1}  we get
\[ - i P_0 (\dd \pi (P ) f)  =\dd \pi_1 (P)f \quad \mbox{ for } \quad 
f \in \cS_+\, .\]
This proves that the restriction of $\dd\pi(\g)$ to $\cS_+$ satisfies
\eqref{eq:heis-inter}, which means that the reflection positive 
representation $\pi$ of $G$ on $(\cE,\cE_+,\theta)$ is a euclidean 
realization of the Schr\"odinger representation $(\pi_1, \cH)$. 
This is a minimal euclidean realization for which 
$\pi_Z(t) = \pi(\exp(tZ))$ is a reflection positive one-parameter group 
because the Lax--Phillips Theorem leads to the requirement that 
$\Spec(i\dd\pi(Z))= \R$. This can only happen for direct integrals of the 
representations $\pi_\lambda$, $\lambda \in \R^\times$, with full support. 

Composing with $\gamma_\lambda$, $\lambda > 0$, we likewise obtain euclidean 
realizations of $\pi_\lambda$, $\lambda > 0$. We thus obtain the following 
analog of Proposition~\ref{prop:ereal-axb} for the Heisenberg group: 

\begin{prop} \mlabel{ereal-heis} Every unitary representation $(\pi^c,\cH)$ 
of $G = \Heis(\R^2,\omega)$ with $-i\dd\pi^c(Z) \geq 0$ has a euclidean realization. 
\end{prop}

\begin{rem} In this remark we show that the representations of the 
Heisenberg group obtained by an outgoing realization of a reflection 
positive one-parameter group do not lead to non-trivial reflection 
positive representations. 

Let $(U_t)_{t \in \R}$ be a reflection positive one-parameter group 
on $(\cE,\cE_+,\theta)$ for which 
$\cE_{\rm fix} = \{0\}$, so that, by the Lax--Phillips Theorem, 
we obtain an equivalence to 
$\cE \cong L^2(\R,\cM)$, $\cE_+ \cong L^2(\R_+,\cM)$ with 
$(U_t f)(x) = f(x-t)$.  

We then have another unitary one-parameter group 
\[ (V_s f)(x) := e^{-isx} f(x) \] 
on $\cE$, combining with $U$ to a unitary representation 
of the $3$-dimensional Heisenberg group $G := \Heis(\R^2)$ because we have
\begin{equation}
  \label{eq:ccr}
 U_t V_s = e^{ist} V_s U_t.
\end{equation}

Although $V$ leaves the subspace $\cE_+$ invariant, we cannot expect it 
to be compatible with the involution $\theta$ which satisfies 
$\theta U_t \theta = U_{-t}$. More concretely, we have 
\[ \widehat{\theta f} (y) = m(y) \hat f(-y),\] 
where $m \: \R \to \U(\cM)$ satisfies 
$m(-y) = m(y)^{-1} = m(y)^*$. 

Let $P$ and $Q$ denote the infinitesimal generators of $U$ and $V$ in 
$\g = \L(G)$ and $Z = [P,Q]$. 
Suppose that $\tau$ is an involution on $G$ with 
$\tau(P) = -P$ and $\tau(Z) = Z$. 
Then $\tau$ induces $-\id$ modulo $\R Z$, so that 
\[ \tau(Q) = - Q + c Z \quad \mbox{ for some } \quad c \in \R,\] 
and this leads to 
\[ \theta V_s \theta = e^{isc} V_{-s}.\] 
More concretely, we obtain 
\[ \widehat{\theta V_s \theta f}(y)
= m(y) \widehat{V_s \theta f} (-y)
= m(y) \widehat{\theta f}(s-y)
= m(y) m(s-y) \hat f(y-s) 
= m(y) m(s-y) \widehat{V_{-s} f} (y),\] 
and thus 
$m(y) m(s-y) = e^{isc}$, which for $y = 0$ leads to 
\[ m(s) = m(0)^{-1} e^{isc} = e^{isc} m(0).\] 
Therefore the inverse Fourier transform of $m$ is the Dirac measure 
\[ \cF^{-1}(m) = m(0) \delta_{-c},\] 
hence positive definite on $\R_+$ if and only if it vanishes on $\R_+$, 
and this is the case for $c \geq 0$. In any case we obtain 
$\hat\cE = \{0\}$ if $\cE_+$ is $\theta$-positive. 
\end{rem}

\appendix 

\section{Rotation invariant reflection positive measures on $\R^d$} 

In this section we take a closer look at rotation invariant tempered measures 
$\nu$ on $\R^d$ which are {\it reflection positive} in the sense that their 
Fourier transform $\hat\nu$ is a reflection positive distribution w.r.t.\ 
($\R^d, \R^d_+,\theta)$, where $\theta(x_0,\bx) = (-x_0,\bx)$. 
In Euclidean Quantum Field Theories, 
the Fourier transforms of these measures are the Schwinger 
distributions $S_2(x,y) = \hat\nu(x-y)$ describing the euclidean 
$2$-point ``functions'' (cf.\ \cite[p.~91]{GJ81}).

The following lemma will be useful in our discussion of examples. 

\begin{lem} \mlabel{lem:delta-fin} 
For a non-zero positive Borel measure $\rho$ on $[0,\infty[$, we consider the 
rotation invariant measure $d\nu(x) := \Theta(x)dx$ on $\R^d$ with the density 
$\Theta(x) := \int_0^\infty \frac{1}{y^2 + x^2}\, d\rho(y)$. 
Then the following are equivalent: 
  \begin{description}
  \item[\rm(i)] There exists a non-zero $x \in \R^d$ with $\Theta(x) < \infty$. 
  \item[\rm(ii)] $\Theta(x) < \infty$ for every non-zero $x \in \R^d$. 
  \item[\rm(iii)] $L^2(\R^d,\nu) \not=\{0\}$. 
  \item[\rm(iv)] $\int_0^\infty \frac{d\rho(y)}{1+y^2} < \infty$. This implies that 
$\rho$ is tempered. 
  \end{description}
If these conditions are satisfied, then the following assertions hold: 
\begin{itemize}
\item[\rm(a)] $\lim_{\|x\| \to \infty} \Theta(x)x^2 = \rho([0,\infty[)> 0.$ 
\item[\rm(b)] $L^2(\R^d,\nu)$ contains a non-zero polynomial if and only 
if $d = 1$ and $\nu$ is finite, and then  every  polynomial in $L^2(\R^d,\nu)$ 
is constant. We have $\nu(\R) = \pi \int_0^\infty \frac{d\rho(y)}{y}$. 
\end{itemize}
\end{lem}

\begin{prf} The equivalence of (i)-(iv) is easy to see. 
Suppose that these conditions are satisfied. 

(a) follows from 
\[ \lim_{x \to \infty} \Theta(x)x^2 
= \lim_{x \to \infty} \int_0^\infty \frac{x^2}{y^2 + x^2}\, d\rho(y)
=  \int_0^\infty \, d\rho(y) > 0.\] 

(b) Let $f = \sum_{j = 0}^N f_j \: \R^d \to \C$ be a polynomial of degree $N$, where 
the $f_j$ are homogeneous of degree $j$. For a suitably normalized measure 
$\sigma$ on the sphere $\bS^{d-1}$, we have the integral formula 
\[ \int_{\R^d} F(x)\, dx = \int_0^\infty \int_{\bS^{d-1}} F(r x)\, d\sigma(x)\, r^{d-1}\, dr,\]
so that 
\[ \int_{\R^d} F(x)\, d\nu(x) = \int_0^\infty \int_{\bS^{d-1}} F(r x)
\, d\sigma(x)\, \Theta(r)r^{d-1}\, dr.\]
The function $|f|^2$ is a polynomial of degree $2N$, and this implies that 
$h(r) := \int_{\bS^{d-1}} |f(r x)|^2\, d\sigma(x)$ also is a polynomial of degree $2N$. 
If 
\[ \int_{\R^d} |f(x)|^2\, d\nu(x) 
= \int_0^\infty h(r)\Theta(r)r^{d-1}\, dr\] 
is finite, (a) implies that 
$\int_1^\infty h(r)r^{d-3}\, dr < \infty$, so that we obtain 
$2N + d- 3 < -1$, i.e., $2N + d < 2$. This leaves only the possibilities 
$d = 1$ and $N= 0$, so that $f$ must be constant. But $1 \in L^2(\R^d,\nu)$ is equivalent 
to the finiteness of the measure $\nu$. In view of 
\[ \nu(\R) 
=  \int_0^\infty \int_{\R} \frac{dx}{y^2 + x^2}\, d\rho(y) 
=  \int_0^\infty \frac{\pi}{y}\, d\rho(y),\] 
$\nu$ is finite if and only if  
$\int_0^\infty \frac{1}{y}\, d\rho(y) < \infty$. 
\end{prf}

\begin{defn} We call a positive Borel measure $\rho$ on $[0,\infty[$ {\it tame} if 
the conditions (i)-(iv) from the preceding lemma are satisfied. 
\end{defn}

\begin{prop} \mlabel{prop:3.15} 
If $\rho$ is tame, then $\nu = \Theta \cdot dx$ with 
$\Theta(x) := \int_0^\infty \frac{1}{y^2 + x^2}\, d\rho(y)$ is 
a tempered measure on $\R^d$ if and only if $d > 2$ or 
\begin{description}
\item[\rm(a)] $d = 1$ and $\int_0^1 \frac{1}{y}\,  d\rho(y) < \infty$. 
\item[\rm(b)] $d= 2$ and $\int_0^1 \ln(y^{-1})\,  d\rho(y) < \infty$. 
\end{description}
In particular, $\rho(\{0\}) = 0$ for $d = 1,2$. 
\end{prop}

\begin{prf} That  $\nu$ is a tempered measure means that, 
for some $N \in \N$, the integral 
\[\int_{\R^d} (1 + p^2)^{-N}\, d\nu(p) 
= \int_0^\infty \int_{\R^d} (1 + p^2)^{-N}\frac{1}{y^2 + p^2}\, dp\, d\rho(y) \]
is finite. This is equivalent to the finiteness of the integral 
\begin{equation}
  \label{eq:doubint}
\int_0^\infty \int_0^\infty (1 + r^2)^{-N}\frac{r^{d-1}}{y^2 + r^2}\, dr\, d\rho(y).
\end{equation}
If \eqref{eq:doubint} is finite, we obtain in particular that 
\[ \int_0^\infty \int_1^2(1 + r^2)^{-N}\frac{r^{d-1}}{y^2 + r^2}\, dr\, d\rho(y) < \infty,\]
and this is equivalent to the tameness of $\rho$. 

Now we split the double integral \eqref{eq:doubint}  according to $r < 1$ and $r > 1$ to find 
\begin{align*}
&\int_0^\infty \int_1^\infty (1 + r^2)^{-N}\frac{r^{d-1}}{y^2 + r^2}\, dr\, d\rho(y) 
= \int_0^\infty I(y)\, d\rho(y) 
\end{align*}
with 
\[ I(y) = \int_1^\infty (1 + r^2)^{-N}\frac{r^{d-1}}{y^2 + r^2}\, dr \leq I(0).\] 
This integral is finite if and only if 
$2N + 2 > d$. 
To see if $I$ is $\rho$-integrable, we 
have to estimate the asymptotics of $I(y)$ for $y \to \infty$. The 
Monotone Convergence Theorem implies that $\lim_{y \to \infty} I(y)= 0$. 
The asymptotics of $I$  is the same as the asymptotics of the function 
\[ J(y) := \int_1^\infty \frac{r^{d-1-2N}}{y^2 + r^2}\, dr 
= \frac{1}{y^2} \int_1^\infty \frac{r^{d-1-2N}}{1+ \frac{r^2}{y^2}}\, dr 
\leq \frac{1}{y^2} \int_1^\infty r^{d-1-2N}\, dr,\] 
provided $2N > d$. 
Now the tameness of $\rho$ implies that $I$ is $\rho$-integrable for $2N > d$. 

Next we turn to the other part of the integral: 
\[ \int_0^\infty \int_0^1 (1 + r^2)^{-N}\frac{r^{d-1}}{y^2 + r^2}\, dr
\, d\rho(y). \] 
Since $1 \leq 1 + r^2 \leq 2$ for $0 \leq r \leq 1$, it is finite if and only 
if the integral 
\[ \int_0^\infty \int_0^1 \frac{r^{d-1}}{y^2 + r^2}\, dr\, d\rho(y) 
= \int_0^\infty K(y)\, d\rho(y) \quad \mbox{ for }  \quad 
K(y) := \int_0^1 \frac{r^{d-1}}{y^2 + r^2}\, dr,\]  
is finite. The integral $K(y)$ is finite for $y > 0$, and 
$K(0) =  \int_0^1 r^{d-3} \, dr$ 
is finite if and only if $d > 2$. 
We thus have to estimate 
the asymptotics of $K$ for $y \to \infty$, and, 
for $d \leq 2$, also for $y \to 0$. 
In view of 
\[ K(y) \leq \int_0^1 \frac{dr}{y^2 + r^2} \leq \frac{1}{y^2},\] 
the finiteness of the integral 
$\int_1^\infty K(y)\, d\rho(y)$ follows for every~$d$ 
from the tameness of $\rho$. 
For $d > 2$, the function $K$ is continuous on $[0,\infty[$, 
hence $\rho$-integrable. This completes the proof for $d > 2$. 

For  $d = 1$ we have 
\[ K(y) 
= \int_0^1 \frac{1}{y^2 + r^2}\, dr
= \frac{1}{y}(\arctan \frac{1}{y}- \arctan 0)  
= \frac{1}{y} \arctan\frac{1}{y} \sim \frac{\pi}{2y}\] 
for $y \to 0$. We thus find the necessary and sufficient 
condition (a) for the finiteness of the $\rho$-integral of~$I$. 

For $d = 2$ we have 
\[ K(y) 
= \int_0^1 \frac{r}{y^2 + r^2}\, dr
= \frac{1}{2} \big( \ln(y^2 + 1) - \ln(y^2)\big) 
\sim \ln(y^{-1}) \] 
for $y \to 0$. Therefore $I$ is $\rho$-integrable over $[0,1]$ 
if and only if (b) is satisfied. 
\end{prf}

The following proposition can also be derived from 
\cite[Prop.~6.2.5]{GJ81}. We include it for the sake of easier reference.

\begin{prop} \mlabel{prop:repo-cm} 
If $\rho$ is a tame measure on $[0,\infty[$ for which 
the measure $\nu = \Theta \cdot dx$ is tempered, then 
the distribution $\hat{\nu} \in C^{-\infty}(\R^d)$ is reflection positive 
for $(\R^d, \R^d_+, \theta)$ and $\theta(x_0, \bx) = (-x_0, \bx)$, i.e., 
\[ \int_{\R^d} \theta \hat  \psi \cdot \oline{\hat \psi}\, d\nu \geq 0 
\quad \mbox{ for } \quad f \in C^\infty_c(\R^d_+).\] 
\end{prop}

\begin{prf} Writing 
\[ \nu = \int_0^\infty \nu_y\, d\rho(y) \quad \mbox{ with } \quad 
d\nu_y(p)= \frac{dp} {y^2 + p^2},\] 
we see that it suffices to show that the distributions 
$\hat{\nu_y}$ are reflection positive. 
First we observe that 
\[ \int_{\R^d} \theta \hat  \psi \cdot \oline{\hat \psi}\, d\nu_y 
= \int_{\R^d} \oline{\hat \psi(-p_0, \bp)} \frac{\hat \psi(p_0, \bp)}{y^2 + p^2}\, dp 
\quad \mbox{ for }\quad  p = (p_0, \bp)\in \R \times \R^{d-1}.\] 
For each $\bp \in \R^{d-1}$, the function 
$h_\bp(p_0) := \hat \psi(-p_0, \bp)$ is a Schwartz function with 
$\supp(\hat h_{\bp}) \subeq ]0, \infty[$, and 
\[ \int_{\R^d} \theta \hat  \psi \cdot \oline{\hat \psi}\, d\nu_y 
= \int_{\R^{d-1}} \int_\R 
\frac{\oline{h_{\bp}(p_0)} h_{\bp}(-p_0)}{p_0^2 + y^2 + \bp^2}\, dp_0 \, d\bp.\] 
If $y = 0$, then $d > 2$, so that $\{0\}$ 
is a Lebesgue zero set in $\R^{d-1}$. Therefore it suffices to show that, for 
$f \in \cS(\R)$ with $\supp(\hat f) \subeq \R_+$ and $y > 0$, the measure 
$d\nu_y^1(x) := \frac{dx}{x^2 + y^2}$ satisfies 
\[ \int_\R f(x) \oline{f(-x)}\, \frac{dx}{x^2 + y^2} 
=  \int_\R f(x) \oline{f(-x)}\, \nu_y^1(x) \, dx \geq 0. \] 
In view of $\hat{\nu_y^1}(p)= \frac{\pi}{y}\ e^{-y|p|}$, 
this follows from 
\begin{align*}
 &\int_\R (\hat f * \oline{\hat f})(p) e^{-y |p|}\, dp 
= \int_{\R_+} (\hat f * \oline{\hat f})(p) e^{-y p} \, dp 
= \cL(\hat f * \oline{\hat f})(y)
= |\cL(\hat f)(y)|^2 \geq 0.
\end{align*} 
More conceptually, the preceding calculation means that $\hat{\nu_y^1}$ is a reflection positive 
function on $\R$, resp., its restriction to the semigroup 
$\R_+$ is also positive definite with respect to the trivial 
involution (cf.\ \cite{NO12}). 
\end{prf}

Note that Proposition~\ref{prop:a.9} below implies in particular that, for $d = 1$, 
the preceding proposition does not cover all reflection positive measures. 
Example that are not covered arise for $\hat\nu(x) = |x|^{-s}$ (on $\R^\times$) 
for $s > 1$.  

\begin{lem} \mlabel{lem:2.3} 
Suppose that $\rho$ is a positive Borel measure on $[0,\infty[$ whose Laplace transform 
$\cL(\rho)$ exists on $\R_+$. Then the function $\cL(\rho)$ on $\R_+$ 
extends to a symmetric distribution 
on $\R$ if and only if $\rho$ is tempered, i.e., 
$\int_1^\infty \frac{1}{y^k}\, d\rho(y) < \infty$ for some $k \in \N$. 
More precisely, for $\ell \in \N_0$, 
\[ \int_0^1 x^\ell \cL(\rho)(x)\, dx < \infty \quad \Leftrightarrow \quad 
\int_1^\infty \frac{1}{y^{\ell + 1}}\, d\rho(y) < \infty.\] 
In particular, $\cL(\rho)$ is locally integrable on $[0, \infty[$ if and only if 
$\int_1^\infty \frac{1}{y}\, d\rho(y) < \infty$. 
\end{lem}

\begin{prf} Since the Laplace transform of the finite measure 
$\rho\res_{[0,1]}$ extends to a continuous function on 
$[0,\infty[$, we may w.l.o.g.\ assume that $\rho([0,1]) = 0$. 
We put $D^\times(x) = \cL(\rho)(|x|)$ for $x \in \R^\times$. In view of 
\cite[Thm.~VIII, \S VII.4]{Schw73}, $D^\times$ extends to a distribution on $\R$ 
if and only if there exists an $\ell \in \N_0$ with 
\begin{equation}
  \label{eq:locint}
\int_0^1 x^\ell D^\times(x)\, dx < \infty.
\end{equation}
We rewrite this integral as follows 
\begin{align*}
\int_0^1 x^\ell D^\times(x)\, dx 
&= \int_0^1 x^\ell \int_1^\infty  e^{-y x}\, d\rho(y) dx 
= \int_1^\infty \int_0^1 x^\ell e^{-y x}\, dx\, d\rho(y) \\ 
&= \int_1^\infty \int_0^y \frac{u^\ell}{y^\ell} e^{-u}\, \frac{du}{y}\, d\rho(y) 
= \int_1^\infty \int_0^y u^\ell e^{-u}\, du\, \frac{1}{y^{\ell+1}}\, d\rho(y).
\end{align*}
In view of $0 < \int_0^\infty u^\ell e^{-u}\,du < \infty$, it follows that 
$\int_0^1 x^\ell D^\times(x)\, dx$ is finite if and only if 
$\int_1^\infty \frac{1}{y^{\ell+1}}\, d\rho(y)$ is finite. 
\end{prf}

\begin{thm} \mlabel{thm:2.2} The reflection positive distributions on $\R$ which are represented 
by a locally integrable function are the Fourier transforms $D = \hat\nu$ 
of measures of the form 
\begin{equation}
  \label{eq:nu}
 \nu = c \delta_0 + \Theta \cdot dx, \qquad 
\Theta(x) = \frac{1}{\pi} \int_0^\infty \frac{y}{y^2 + y^2}\, d\rho(y), 
\end{equation}
where $c \geq 0$ and $\rho$ is a positive Radon measure on $\R_+$ satisfying 
\begin{equation}
  \label{eq:r1}
\rho(]0,1]) < \infty \quad \mbox{ and } \quad 
\int_1^\infty \frac{1}{y}\, d\rho(y) < \infty.
\end{equation}
\end{thm}

\begin{prf} Let $D \in C^{-\infty}(\R)$ be reflection positive. 
Then its restriction to $\R_+$ is positive definite w.r.t.\ the involution 
$x^* = x$, so that \cite[Thm.~4.13]{NO12} 
implies $D\res_{\R_+}$ is represented by an analytic function which is the 
Laplace transform $\cL(\rho)$ of a positive Radon measure $\rho$ on $[0,\infty[$. 
If $D$ is represented by a locally integrable function, $\cL(\rho)$ is locally integrable, 
so that Lemma~\ref{lem:2.3} leads to 
\[ \int_1^\infty \frac{1}{y}\, d\rho(y) < \infty.\] 
We now have almost everywhere on $\R$  
\[ D(x) = \cL(\rho)(|x|) 
= \int_0^\infty e^{-y|x|}\, d\rho(y) =  \hat\nu(x) 
\quad \mbox{ for } \quad \nu = \rho(\{0\})\delta_0 + \Theta \cdot dx,\] 
provided $\nu$ is a tempered measure. Here we used that $e^{-y|x|}$ is the Fourier transform 
of the measure $\frac{1}{\pi}\frac{y \cdot dx}{y^2 + x^2}$. 
In view of Proposition~\ref{prop:3.15}, the temperedness of $\nu$ is equivalent to 
\eqref{eq:r1}. 

Suppose, conversely, that \eqref{eq:r1} is satisfied. 
 Then Proposition~\ref{prop:3.15} implies that $\nu$ 
is tempered, so that $D := \hat\nu$ is a positive definite distribution. 
That it is reflection positive follows from the reflection positivity 
of the functions $\hat\nu_y(x) = e^{-y|x|}$. 
\end{prf}

\begin{rem}
  \mlabel{rem:4.4}  
The integral representation from Theorem~\ref{thm:2.2} does not extend 
to all reflection positive distributions on $\R$. 
Any distribution of the form 
$E := P({1\over i}{d\over dx})\delta_0$, where $P$ is a non-negative even polynomial, 
is reflection positive. In this case $E \res_{\R_+} = 0$. 

We also conclude that, for any reflection positive distribution 
$D$ on $\R$, the distribution $D + E$ is another extension of 
$D\res_{\R_+}$, so that reflection positive extensions are not unique. 
\end{rem} 

\begin{prop} \mlabel{prop:a.9} For every tempered measure $\rho$ on $[0,\infty[$, there exists a 
reflection positive distribution $D$ on $\R$ with $D\res_{\R_+} = \cL(\rho)$. 
\end{prop}

\begin{prf} In view of Lemma~\ref{lem:2.3}, we may assume that 
$\int_1^\infty \frac{1}{\lambda^{4N}}\, d\rho(\lambda) < \infty.$
Then the measure 
\[ d\tilde\rho(\lambda) :=  \frac{1}{1 + \lambda^{4N}}\, d\rho(\lambda) \] 
on $[0,\infty[$ is finite. Its Laplace transform 
therefore defines a continuous reflection positive function 
$\tilde D(x) := \cL(\tilde\rho)(|x|)$ 
(\cite[Cor.~3.3]{NO12}). For the polynomial $P(x) := 1 + x^{4N}$, we 
now obtain the positive definite distribution 
$D := P({1\over i}{d\over dx}) \tilde D,$ 
and its restriction to $\R_+$ is given by 
\[ D(x) = P\Big({1\over i}{d\over dx}\Big) \cL(\tilde\rho)(x)  
= \int_0^\infty (1 + \lambda^4)e^{-\lambda x}\, d\tilde\rho(\lambda) 
= \int_0^\infty e^{-\lambda x}\, d\rho(\lambda) = \cL(\rho)x).\] 
Therefore $D$ is a reflection positive extension of $\cL(\rho)$. 
\end{prf}

\section{Lorentzian invariant tempered measures on 
the light cone} 

For the following theorem, we recall the $L^{\up}$-invariant 
measures $\mu_m$ from Definition~\ref{def:mum}.  
It is our version of the 
Lehmann Spectral Formula from Quantum Field Theory 
(\cite[Thm.~6.2.4]{GJ81}) describing the $2$-point functions of 
Poincar\'e invariant field theories. The proof given below provides many 
details skipped in \cite{GJ81}.

\begin{thm} \mlabel{thm:5.11b} 
For every $L^{\up}$-invariant Borel measure $\mu$ 
on the closed  forward light cone $\oline{V_+}$, there exists a $\sigma$-finite 
Borel measure $\rho$ on $[0,\infty[$ and a constant $c \geq 0$ so that 
\begin{equation}
  \label{eq:nu-rho}
\mu = c \delta_0 + \int_0^\infty \mu_m\, d\rho(m),
\end{equation}
where $\rho(\{0\}) = 0$ for $d = 1$. 
\begin{description}
\item[\rm(a)] The measure 
$\mu$ is tempered if and only if the following conditions are satisfied: 
\begin{description}
\item[\rm(C1)] $\rho$ is a tempered measure on $[0,\infty[$. 
\item[\rm(C2)] $\int_0^1 \frac{1}{m}\,  d\rho(m) < \infty$ for $d = 1$. 
\item[\rm(C3)] $\int_0^1 \ln(m^{-1})\,  d\rho(m) < \infty$ for $d = 2$. 
In particular, $\rho(\{0\}) = 0$.  
\end{description}
\item[\rm(b)] Let $\pr \: \R^d \to \R^{d-1}, (p_0, \bp) \mapsto \bp$ be the projection map. 
Then the measure $\tilde\mu = \pr_*\mu$ on $\R^{d-1}$ is tempered if and only if 
$\mu$ is tempered and, in addition, 
\[ \int_1^\infty\frac{1}{m}\, d\rho(m) < \infty.\] 
If this condition is not satisfied, then every Borel subset $E \subeq \R^{d-1}$ 
satisfies $\tilde\mu(E) \in \{ 0, \infty\}$, so that 
$L^2(\R^{d-1}, \tilde\mu) = \{0\}$. 
\end{description}
\end{thm}

\begin{prf} In \cite[Thm.~IX.33]{RS75}, the decomposition \eqref{eq:nu-rho} is stated only 
in the case $d = 4$, but the proof works in the general case, where it leads to an 
integral representation of $\mu$ in terms of a measure $\rho$ on $[0,\infty[$, 
whose restriction to $]0,\infty[$ is a Radon measure with possibly 
infinite mass for the interval $[0,1]$ 
(see also \cite[Lemma~9.1.2/3]{vD09} for a description of the ergodic 
$\OO_{1,d-1}(\R)$-invariant measures on $\R^d$). 

(a) Now the problem consists in characterizing the temperedness of $\mu$ in terms 
of properties of the measure $\rho$. To this end, we may w.l.o.g.\ assume that 
$c = 0$. Temperedness of $\mu$ is equivalent to the 
existence of an $N \in \N$ such that the following integral is finite 
\begin{align*}
 \int_{\R^d} \frac{1}{(1 + p^2)^N}\, d\mu(p) 
&= \int_0^\infty \int_{\R^d} \frac{1}{(1 + p^2)^N}\, d\mu_m(p)\, d\rho(m) \\
&= \int_0^\infty \int_{\R^{d-1}} \frac{1}{(1 + (m^2 + \bp^2) + \bp^2)^N}
\frac{1}{\sqrt{m^2 + \bp^2}}\, d\bp\, d\rho(m).
\end{align*}

For $d = 1$, the preceding formula simplifies to 
\begin{align*}
 \int_{\R} \frac{1}{(1 + p^2)^N}\, d\mu(p) 
&= \int_0^\infty \frac{1}{(1 + m^2)^N}\frac{1}{m}\, d\rho(m).
\end{align*}
There exists an $N \in \N_0$ for which such an integral is finite 
if and only if (a) is satisfied and $\rho$ is a tempered. 

>From now on we assume that $d > 1$. 
Evaluating the integral over $\R^{d-1}$ in polar coordinates, 
we see that this is equivalent to the finiteness of the double integral 
\begin{equation}
  \label{eq:rho-int}
\int_0^\infty \int_0^\infty \frac{1}{(1 + m^2 + 2r^2)^N}
\frac{1}{\sqrt{m^2 + r^2}} r^{d-2}\, dr\, d\rho(m).
\end{equation}
The finiteness of \eqref{eq:rho-int} implies in particular 
the finiteness of the integral 
\[ \int_0^\infty \int_1^2\frac{1}{(1 + m^2 + 2r^2)^N}
\frac{1}{\sqrt{m^2 + r^2}} r^{d-2}\, dr\, d\rho(m),\] 
which is equivalent to the finiteness of the integral 
\[ \int_0^\infty \frac{1}{(1 + m^2)^N}\frac{1}{\sqrt{m^2 + 1}}\, d\rho(m) 
=  \int_0^\infty \frac{1}{(1 + m^2)^{N+ \frac{1}{2}}}\, d\rho(m).\] 
This implies that $\rho([0,M]) < \infty$ for every $M > 0$ and 
that $\rho$ is tempered. 

Suppose, conversely, that $\rho$ is tempered. We discuss the integral \eqref{eq:rho-int} 
by splitting it into the two pieces corresponding to $r < 1$ and $r > 1$. 
For the $r >1$-part we obtain for $N \geq N_0$ and $N_0 > \frac{d-1}{2}$: 
\begin{align*}
& \int_1^\infty \frac{1}{(1 + m^2 + 2r^2)^N}
\frac{1}{\sqrt{m^2 + r^2}} r^{d-2}\, dr \\
&\leq  \frac{1}{\sqrt{1 + m^2}} 
\int_1^\infty \frac{1}{(1 + m^2 + 2r^2)^{N_0}(1 + m^2 + 2r^2)^{N-N_0}} r^{d-2}\, dr \\
&\leq  \frac{1}{\sqrt{1 + m^2}} \frac{1}{ (3 + m^2)^{N-N_0}} 
\int_1^\infty \frac{1}{(1 + 2r^2)^{N_0}} r^{d-2}\, dr.  
\end{align*}
Our condition on $N_0$ ensures that this integral is finite, and now the temperedness 
of $\rho$ implies that 
\[ \int_0^\infty \int_1^\infty \frac{1}{(1 + m^2 + 2r^2)^N}
\frac{1}{\sqrt{m^2 + r^2}} r^{d-2}\, dr\, d\rho(m) < \infty\] 
if $N$ is sufficiently large.

Next we note that, 
\begin{align*}
&\int_1^\infty \int_0^1 \frac{1}{(1 + m^2 + 2r^2)^N}
\frac{1}{\sqrt{m^2 + r^2}} r^{d-2}\, dr\, d\rho(m) \\
&\leq 
\int_1^\infty \int_0^1 \frac{1}{(1 + m^2)^N}
\frac{1}{m} \, dr\, d\rho(m) 
= \int_1^\infty \frac{1}{(1 + m^2)^N}\frac{1}{m}\, d\rho(m) < \infty.
\end{align*}
This already shows that $\int_1^\infty \mu_m\, d\rho(m)$ is tempered 
if $\rho$ is tempered. 

It remains to consider the integral 
\[ \int_0^1 \int_0^1 \frac{1}{(1 + m^2 + 2r^2)^N}
\frac{1}{\sqrt{m^2 + r^2}} r^{d-2}\, dr\, d\rho(m), \] 
which is finite if and only if 
\[ \int_0^1 \int_0^1 \frac{r^{d-2}}{\sqrt{m^2 + r^2}}\, dr\, d\rho(m) < \infty. \] 
Let 
\[ I(m) := \int_0^1 \frac{r^{d-2}}{\sqrt{m^2 + r^2}} \, dr,\] 
and observe that this defines a continuous function for $m > 0$ with 
$I(0)= \int_0^1 r^{d-3}\, dr < \infty$ if and only if $d > 2$. 
Since $\rho([0,1]) < \infty$, we conclude that, for $d\geq 3$, 
$\mu$ is a tempered measure if and only if $\rho$ has this property. 

This leaves the case $d = 2$. Then 
\[ I(m) = \int_0^1 \frac{1}{\sqrt{m^2 + r^2}}\, dr 
= \ln(1 + \sqrt{1 + m^2}) - \ln(m).\] 
We thus find condition (b), 
which is needed to ensure that $\mu$ is a tempered measure. 

(b) The measure $\tilde\mu$ is of the form $\tilde\Theta(\bp)\, d\bp$ for 
\begin{equation}
  \label{eq:tildemu}
\tilde \Theta(\bp) = \int_0^\infty \frac{1}{\sqrt{m^2 + \bp^2}}\, d\rho(m).
\end{equation}

For $d = 1$, the measure 
$\tilde\mu$ is a point measure which is tempered if and only if it is 
finite. This is equivalent to 
\[ \tilde\Theta(0) =  \int_0^\infty \frac{1}{m}\, d\rho(m) < \infty.\] 

Suppose that $d > 1$. If $\tilde\mu$ is finite on a set of positive 
Lebesgue measure, then there exists a non-zero $\bp \in \R^{d-1}$ with 
$\tilde\Theta(\bp) < \infty$. This implies that compact subsets of $[0,\infty[$ have finite 
$\rho$-measure ($\rho$ is a Radon measure), and that 
\begin{equation}
  \label{eq:1m-int}
\int_1^\infty \frac{1}{m}\, d\rho(m) < \infty.
\end{equation}
In particular, $\rho$ is tempered. 
If \eqref{eq:1m-int} is not satisfied, then $\tilde\Theta(\bp) = \infty$ 
implies that $\tilde\mu(E) = \infty$ for 
every Borel subset $E \subeq \R^{d-1}$ of positive Lebesgue measure and 
$\tilde\mu(E) = 0$ if $E$ is a Lebesgue zero set. 

Let us assume that $\rho$ is a Radon measure on 
$[0,\infty[$ satisfying \eqref{eq:1m-int}. 
Temperedness of $\tilde\mu$ is equivalent to the 
existence of an $N \in \N$ such that 
\begin{align*}
 \int_{\R^{d-1}} \frac{1}{(1 + \bp^2)^N}\, d\tilde\mu(p) 
&= \int_0^\infty \int_{\R^{d-1}} \frac{1}{(1 + \bp^2)^N}
\frac{1}{\sqrt{m^2 + \bp^2}}\, d\bp\, d\rho(m) < \infty.
\end{align*}
Evaluating the integral over $\R^{d-1}$ in polar coordinates, 
we see that this is equivalent to the finiteness of the double integral 
\begin{equation}
  \label{eq:rho-intb}
\int_0^\infty \int_0^\infty \frac{1}{(1 + r^2)^N}
\frac{1}{\sqrt{m^2 + r^2}} r^{d-2}\, dr\, d\rho(m).
\end{equation}

We discuss the integral \eqref{eq:rho-intb} 
by splitting it into the two pieces corresponding to $r < 1$ and $r > 1$. 
For the $r >1$-part we obtain for $N  > \frac{d-1}{2}$: 
\[ \int_1^\infty \frac{1}{(1 + r^2)^N}
\frac{1}{\sqrt{m^2 + r^2}} r^{d-2}\, dr 
\leq  \frac{1}{\sqrt{1 + m^2}} 
\int_1^\infty \frac{r^{d-2}}{(1 +r^2)^{N}}\, dr.\] 
Our condition on $N$ ensures that this integral is finite, 
and now the finiteness of $\int_0^\infty\, \frac{1}{\sqrt{1+m^2}}\, d\rho(m)$ 
implies the finiteness of the double integral 
\[ \int_0^\infty \int_1^\infty \frac{1}{(1 + r^2)^N}
\frac{1}{\sqrt{m^2 + r^2}} r^{d-2}\, dr\, d\rho(m).\] 

Next we note that, 
\[ \int_1^\infty \int_0^1 \frac{1}{(1 + r^2)^N}
\frac{ r^{d-2}}{\sqrt{m^2 + r^2}}\, dr\, d\rho(m) 
\leq \int_1^\infty \frac{1}{m} \, d\rho(m) \cdot 
\int_0^1 \frac{r^{d-2}}{(1 + r^2)^N}\, dr < \infty.\] 
It remains to consider the integral 
\[ \int_0^1 \int_0^1 \frac{1}{(1 + r^2)^N}
\frac{r^{d-2}}{\sqrt{m^2 + r^2}} \, dr\, d\rho(m), \] 
which is finite if and only if 
\[ \int_0^1 \int_0^1 \frac{r^{d-2}}{\sqrt{m^2 + r^2}}\, dr\, d\rho(m) < \infty. \] 
In (a) we have seen that this integral is finite 
if $d > 2$, and that, for $d = 2$, its finiteness is equivalent to (C3). 
This proves (b). 
\end{prf}

\begin{rem} \mlabel{rem:b.2} 
For the measure $\nu = \Theta\cdot dx$ on $\R^d$, where 
$\Theta(x) = \frac{1}{\pi} \int_0^\infty\frac{1}{m^2 + x^2}\, d\rho(m)$, the projection 
$\tilde\nu$ to $\R^{d-1}$ under $\pr(x_0,\bx) := \bx$, is of the form 
$\tilde\Theta(\bx)\, d\bx$, where 
\begin{align*}
\tilde\Theta(\bx) 
&= \frac{1}{\pi} \int_0^\infty \int_\R \frac{1}{m^2 + \bx^2 + x_0^2}\, dx_0\, d\rho(m) 
= \int_0^\infty \frac{1}{\sqrt{m^2 + \bx^2}} 
\frac{1}{\pi} \int_\R \frac{\sqrt{m^2 + \bx^2}}{m^2 + \bx^2 + x_0^2}\, dx_0\, d\rho(m) \\
&= \int_0^\infty \frac{1}{\sqrt{m^2 + \bx^2}}\, d\rho(m).
\end{align*}
This implies that $\tilde\nu = \tilde\mu$ (cf.\ \eqref{eq:tildemu}).
\end{rem}

\section{Positive definite functions} 
\mlabel{app:d}

In this appendix we collect some definitions and results concerning 
positive definite functions and kernels. 

\begin{defn} \mlabel{def:1.5} Let $X$ be a set and $\cF$ be a complex Hilbert space.

\par (a)  A function $K \: X \times X \to B(\cF)$ is called a 
{\it $B(\cF)$-valued kernel}. A $B(\cF)$-valued kernel $K$ on $X$ is said to be 
{\it positive definite} if, 
for every finite sequence $(x_1, v_1), \ldots, (x_n,v_n)$ in $X \times \cF$, 
\[ \sum_{j,k = 1}^n \la K(x_j, x_k)v_k, v_j \ra \geq 0. \] 

\par (b) If $(S,*)$ is an involutive semigroup, then a 
function $\phi \: S \to B(\cF)$ is called {\it positive definite} 
if the kernel $K_\phi(s,t) := \phi(st^*)$ is positive definite. 
\end{defn}

Positive definite kernels can be characterized as those 
for which there exists a Hilbert space $\cH$ and a 
function $\gamma \: X \to B(\cH,\cF)$ such that 
\begin{equation}
  \label{eq:evalprod}
K(x,y) = \gamma(x)\gamma(y)^* \quad \mbox{ for } \quad x,y \in X 
\end{equation}
(cf.\ \cite[Thm.~I.1.4]{Ne00}). 
Here one may assume that the vectors 
$\gamma(x)^*v$, $x \in X, v \in \cF$, span a dense subspace of 
$\cH$. If this is the case, then the pair $(\gamma,\cH)$ is called a {\it realization 
of $K$}. 
The map $\Phi \: \cH \to \cF^X, \Phi(v)(x) := \gamma(x)v$, 
then realizes $\cH$ as a Hilbert subspace of $\cF^X$ 
with continuous point evaluations $\ev_x \: \cH \to \cF, f \mapsto f(x)$. 
Then $\Phi(\cH)$ is the unique Hilbert space in $\cF^X$ with continuous point evaluations 
$\ev_x$, for which $K(x,y) = \ev_x \ev_y^*$ for $x,y \in X$. 
We write $\cH_K \subeq \cF^X$ for this subspace and call it 
the {\it reproducing kernel Hilbert space with kernel~$K$}.

\begin{ex} \mlabel{ex:vv-gns} (Vector-valued GNS construction) 
(cf.\ \cite[Sect.~3.1]{Ne00}) Let $(\pi, \cH)$ be a representation of the 
unital involutive semigroup $(S,*)$, $\cF \subeq \cH$ be a closed subspace for which 
$\pi(S)\cF$ is total in $\cH$ and $P \: \cH \to \cF$ denote the orthogonal projection. 
Then $\phi(s) := P\pi(s)P^*$ is a $B(\cF)$-valued positive definite function 
on $S$  with $\phi(\1) = \1_\cF$ because $\gamma(s) := P\pi(s) \in B(\cH,\cF)$ 
satisfies 
\[ \gamma(s)\gamma(t)^* = P\pi(st^*)P^* = \phi(st^*).\] 
The map 
\[ \Phi \: \cH \to \cF^S, \quad \Phi(v)(s) = \gamma(s)v = P \pi(s)v \] 
is an $S$-equivariant realization of $\cH$ as the reproducing kernel space
$\cH_\phi \subeq \cF^S$, on which $S$ acts by right translation, i.e., 
$(\pi_\phi(s)f)(t) = f(ts)$. 

Conversely, let $S$ be a unital involutive semigroup 
and $\phi \: S \to B(\cF)$ be a positive definite function with 
$\phi(\1) = \1_\cF$. 
Write $\cH_\phi \subeq \cF^S$ for the corresponding reproducing kernel space and 
$\cH_\phi^0$ for the dense subspace spanned by 
$\ev_s^*v, s \in S, v \in \cF$. 
Then $(\pi_\phi(s)f)(t) := f(ts)$ defines a 
$*$-representation of $S$ on $\cH_\phi^0$. 
We say that $\phi$ is {\it exponentially bounded} if 
all  operators $\pi_\phi(s)$ are bounded, so that we actually 
obtain a representation of $S$ by bounded operators on $\cH_\phi$. 
As $\1_\cF = \phi(\1) = \ev_\1 \ev_\1^*$, the map 
$\ev_\1^* \: \cF \to \cH$ is an isometric inclusion, so that we may identify 
$\cF$ with a subspace of $\cH$. Then $\ev_\1 \: \cH \to \cF$ corresponds to the 
orthogonal projection onto $\cF$ and 
$\ev_\1 \circ \pi_\phi(s) = \ev_s$ leads to 
\begin{equation}
  \label{eq:factori}
\phi(s) = \ev_s \ev_\1^* =\ev_\1 \pi_\phi(s) \ev_\1^*.
\end{equation}

If $S = G$ is a group with $s^* = s^{-1}$, then $\phi$ is always exponentially bounded and the 
representation $(\pi_\phi, \cH_\phi)$ is unitary. 
\end{ex}

\begin{lem} \mlabel{lem:mult} Let $(S,*)$ be a unital involutive semigroup 
and $\phi \: S \to B(\cF)$ be a positive definite function  
with $\phi(\1) = \1$. We write $(\pi_\phi, \cH_\phi)$ for the representation 
on the corresponding reproducing kernel Hilbert space 
$\cH_\phi \subeq \cF^S$ by $(\pi_\phi(s)f)(t) := f(ts)$. 
Then the inclusion 
\[ \iota \: \cF \to \cH_\phi, \quad \iota(v)(s) := \phi(s)v \]
 is surjective if and only if 
$\phi$ is multiplicative, i.e., a representation. 
\end{lem}

\begin{prf} If $\phi$ is multiplicative, then 
$(\pi(s) \iota(v))(t) = \phi(ts)v = \phi(t)\phi(s)v \in \iota(\cF)$. 
Therefore the $S$-cyclic subspace $\iota(\cF)$ is invariant, which 
implies that $\iota$ is surjective. 

Suppose, conversely, that $\iota$ is surjective. 
This means that each $f \in \cH_\phi$ satisfies 
$f(s) = \phi(s)f(\1)$. For $v \in \cF$ and $t,s \in S$, this leads to 
\[ \phi(st)v = \pi(t)(\iota(v))(s) 
= \phi(s) (\pi(t)\iota(v))(\1) 
= \phi(s) \iota(v)(t) = \phi(s) \phi(t)v.\] 
Therefore $\phi$ is multiplicative. 
\end{prf}

\begin{rem} The preceding lemma can also be expressed without referring to positive definite 
functions and the corresponding reproducing kernel space. In this context it asserts the following.
Let $\pi \: S \to B(\cH)$ be a $*$-representation of a unital involutive semigroup 
$(S,*)$, $\cF \subeq \cH$ be a closed cyclic subspace and 
$P \: \cH \to \cF$ the orthogonal projection. Then the function 
\[ \phi \: S \to B(\cF), \quad \phi(s) := P \pi(s)P^* \] 
is multiplicative if and only if $\cF = \cH$. 
\end{rem}

The following lemma is an abstraction of \cite[Thm.~6.2.2]{GJ81}. 

\begin{lem} \mlabel{lem:exp-crit}
Let $V$ be a vector spaces over $\K \in \{\R,\C\}$ and $\beta \: V \times V\to \K$ 
be a hermitian form on $V$ (for $\K = \R$ this means that it is symmetric and bilinear). 
Then the kernel $e^\beta$ is positive definite if and only if 
$\beta$ is positive semidefinite.   
\end{lem}

\begin{prf} If $\beta$ is positive semidefinite, then  
$\beta$ defines a positive definite kernel on $V \times V$. 
Hence the kernels $\beta^n(x,y) := \beta(x,y)^n$ are also positive definite, 
see \cite[Rem. I.17(b)]{Ne00},
and therefore $e^\beta = \sum_{n = 0}^\infty \frac{\beta^n}{n!}$ is positive definite. 

If, conversely, $e^\beta$ is a positive definite kernel, 
let $x_1,\ldots,x_N\in V$ and $c_1,\ldots,c_N\in \K$. For $t\in \R$
define
\[ y_1=tx_1,\ldots, y_N=t x_N,\quad y_{N+1}=0,\ldots,y_{2N}=0\] 
and
\[ d_1=c_1,\ldots, d_N=c_N,\quad d_{N+1}=-c_1,\ldots,d_{2N}=-c_N. \] 
Then
\[ 0 \leq \sum_{i,j=1}^{2N}e^{\beta(y_i,y_j)}d_i\oline{d_j}=
\sum_{i,j=1}^{N}(e^{t^2\beta(x_i,x_j)}-1)c_i\oline{c_j}\] 
and letting $t$ tend to zero after dividing by $t^2$ yields the claim.
\end{prf}

\section{Distribution vectors and tempered measures} 

In this appendix we take a closer look at distribution vectors 
for representations of vector spaces by multiplication operators 
on $L^2$-spaces.

\subsection{Representations of cones} 

Let $V$ be a finite-dimensional real vector space, 
$\tau \: V \to V$ be a linear involution and 
 $\Omega \subeq V$ be an open convex cone invariant 
under the involution $v \mapsto v^\sharp = -\tau(v)$. The cone 
\[ \hat\Omega = \{ \alpha \in V_\C^* \: (\forall
x = x^\sharp \in \Omega)\ \alpha(x) \geq 0\}\]
parametrizes the bounded characters 
of $\Omega$ by assigning to $\alpha \in \hat\Omega$ the character 
$e_\alpha(v) := e^{-\alpha(v)}$.  

Any $\sigma$-finite measure $\mu$ on $\hat\Omega$ 
defines a contraction representation of $\Omega$ on $\cH := L^2(\hat\Omega, \mu)$ by 
\begin{equation}
  \label{eq:pimu}
(\pi_\mu(v)f)(\alpha) := e^{-\alpha(v)} f(\alpha) \quad \mbox{ for } \quad 
v \in \Omega, \alpha \in \hat\Omega.
\end{equation}
On the subspace $V^c := i V^{-\tau} \oplus V^{\tau} \subeq V_\C$ 
formula \eqref{eq:pimu} defines a unitary representation $\pi_\mu^c$ because 
every $\alpha \in \hat\Omega$ is purely imaginary on $V^c$.

The following lemma is a supplement to the 
generalized Bochner--Schwartz Theorem from \cite{NO12} because
it tells us which measures $\mu$ on $\hat\Omega$ actually have 
Fourier--Laplace transforms defining distributions on $\Omega$. 

\begin{lem} \mlabel{lem:auto} Let $V$ be a finite-dimensional vector space, 
$\tau \in \GL(V)$ an involution and 
$\Omega \subeq V$ be an open convex cone invariant under the involution 
$v \mapsto v^\sharp = - \tau(v)$. Further, let $\mu$ be a positive Borel 
measure on $\hat\Omega$ such that, for every 
$\phi \in C^\infty_c(\Omega)$ the Fourier--Laplace transform 
$\hat\phi(\alpha) := \int_{\Omega} \phi(x) e^{-\alpha(x)}\, dx$ is $\mu$-integrable 
on $\hat\Omega$. 
Then $\hat\mu(\phi) := \mu(\oline{\hat\phi})$ defines a distribution on $\Omega$. 
In the special case $V = \Omega$ and $\tau = \id_V$, the measure 
$\mu$ is tempered. 
\end{lem}

\begin{prf} We consider the linear functional 
\[ E \: C^\infty_c(\Omega) \to \C, \quad 
E(\phi) := \mu(\hat\phi).\] 
Then, for every test function $\psi$ and $\psi^\sharp := \psi^* \circ \tau$, 
the regularized functional 
$E_\psi(\phi) := E(\psi * \phi * \psi^\sharp) = \mu(|\hat\psi|^2 \hat\phi)$ is continuous 
because it is the Fourier transform of the bounded measure $|\hat\psi|^2 \mu$. 

Now we consider a sequence $0 \leq \psi_n \in C^\infty_c(\Omega)$ with 
$\hat\psi_n(0) = \int_{\Omega} \psi_n(x)\, d\mu_V(x) = 1$ 
and $\supp(\psi_n) \to \{0\}$ (a $\delta$-sequence). Then 
\[ \hat\psi_n(\alpha) =  \int_{V} \psi_n(x) e^{-\alpha(x)}\,  d\mu_V(x) \to e^{-\alpha(0)} = 1 \] 
holds pointwise. 
We also have $\|\hat\psi_n\|_\infty \leq \|\psi_n\|_1 = 1$ for the 
sup-norm of $\hat\Omega$, so that Dominated Convergence implies 
that $E_{\psi_n} \to E$ pointwise on $C^\infty_c(V)$. 

The Uniform Boundedness Theorem, applied to the 
restriction of $E_{\psi_n}$ to the Fr\'echet spaces 
$C^\infty_K(\Omega)$, $K \subeq \Omega$ compact, now implies that 
$E$ is continuous on every subspace $C^\infty_K(\Omega)$, hence continuous 
on $C^\infty_c(\Omega)$. This means that $E$ is a positive definite distribution. 

For $V = \Omega$ and $\tau = \id_V$, the 
Bochner--Schwartz Theorem (\cite[Thm.~XVIII, \S VII.9]{Schw73}) 
further implies that $\mu$ is tempered. 
\end{prf}

\begin{rem}  \mlabel{rem:d.2}
Suppose that $D \in C^{-\infty}(\Omega)$ is a positive definite distribution, 
so that we can use the  generalized Bochner--Schwartz Theorem (\cite[Thm.~4.11]{NO12}) 
to write it as the Fourier--Laplace transform of a measure $\mu$ in the sense of  
\[ D(\phi) = \int_{\hat\Omega} \oline{\hat\phi}\, d\mu \quad \mbox{ for } \quad 
\phi \in C^\infty_c(\Omega), \quad \mbox{ where }\quad 
\hat\phi(\alpha) := \int_{\Omega} \phi(x) e^{-\alpha(x)}\, dx\] 
is the Fourier--Laplace transform of $\phi$. 
Then we obtain an isomorphism 
\[ \Gamma \: L^2(\hat\Omega, \mu) \to \cH_D \subeq C^{-\infty}_c(\Omega), \quad 
\Gamma(f)(\phi) = \la f, \hat \phi \ra, \] 
under which $D$ corresponds to the 
constant function~$1$. In particular, 
$\hat \phi \in L^2(\hat\Omega, \mu)$ for $\phi \in C^\infty_c(\Omega)$. 
\end{rem}

\begin{rem} In the special case $\tau = -\id_V$ the integrability of 
the functions $\hat\phi$ is equivalent to the existence of the 
Laplace transform 
\[ \hat\cL(\mu)(x) = \int_{\hat\Omega} e^{-\alpha(x)}\, d\mu(\alpha) 
\quad \mbox{ for } \quad x \in \Omega.\] 
This follows from 
\[\int_{\hat\Omega} \hat\phi\, d\mu = \int_\Omega \cL(\mu)(x) \phi(x)\, dx\]
because $\cL(\mu)$ is continuous whenever it exists (by  the Dominated Convergence Theorem). 
\end{rem}

\begin{ex} The constant function $1$ on $\hat\Omega$ need 
not be a distribution vector for the representation 
$\pi_\mu^c$ of $V^c$, even if $\hat\mu$ defines a distribution on $\Omega$. 
A simple example is provided for $V = \R$, 
$\tau = -\id_V$ and $\Omega = \R_+$ by the measure 
$d\mu(t) := e^{\sqrt t} \, dt$ on $\hat\Omega = [0,\infty[$. 
Then the Laplace transform $\cL(\mu)$ is defined on $\Omega$, 
so that  each function $\hat\phi$, $\phi \in C^\infty_c(\Omega)$, is integrable, 
but $\cL(\mu)$ does not extend to a distribution on $\R$ because 
$\int_1^\infty x^{-\ell}\, d\mu(x) = \infty$ for every $\ell \in \N$ 
(cf.~Lemma~\ref{lem:2.3}). 
\end{ex}

\begin{ex} We now turn to the special case 
$V = \R^d$, $\Omega = \R^d_+$ and $\tau(x) = (-x_0, \bx)$ for 
$x = (x_0, \bx)$. 

(a) Identifying $\hat\Omega$ with $[0,\infty[  \times \R^{d-1}$, the contraction 
representation of $\Omega$ on $L^2(\hat\Omega,\mu)$ is given by 
\[ (\pi_\mu(x)f)(p) = e^{-x_0 p_0} e^{-i \bx \bp} f(p) \quad \mbox{ for } \quad 
p = (p_0, \bp).\] 
Here we use the embedding 
\[\iota \: [0,\infty[ \times \R^{d-1} \to \hat\Omega \subeq V_\C^* \cong \C^d, \quad 
\iota(p_0, \bp)(x) := p_0 x_0 + i \bx \bp.\] 
The corresponding unitary representation of 
$V^c = i\R \oplus \R^{d-1} \cong \R^d$ is given by 
\[ (\pi_\mu^c(x)f)(p) = e^{-i xp} f(p) = e^{-ix_0 p_0} e^{-i \bx \bp} f(p) 
\quad \mbox{ for } \quad 
x = (x_0, \bx), p = (p_0,\bp).\] 
In view of Corollary~\ref{cor:3.17}, the function $1$ is a distribution vector for 
$\pi^c_\mu$ if and only if $\mu$ is a tempered measure on $\R^d$. 

(b) Note that $\pi_\mu(x)1 \in \cH$ for every $x \in \Omega$ 
is equivalent to the existence of the Fourier--Laplace transform 
of $\mu$ as a function on the open half space $\Omega$. 
In view of $|e^{-\alpha(x)}| = e^{-p_0 x_0}$ for 
$\alpha = (p_0, i \bp)$, 
this in turn is equivalent to the square-integrability of 
all functions $e_x(\alpha) := e^{-\alpha(x)}$, $x \in \Omega$, 
on $\hat\Omega$. This in turn is equivalent to 
the integrability of all functions $e^{-t p_0}$, $t > 0$. 
If $\mu$ is supported by the forward light cone 
$\oline{V_+} \subeq \R^d$, then the temperedness of $\mu$ implies the 
existence of its Laplace transform on $\Omega$.   
\end{ex}

The following lemma makes the condition from 
Lemma~\ref{lem:auto} more explicit in terms of the 
product decomposition $\Omega = \R_+ \times \R^{d-1}$. 

\begin{lem} For $V = \R^d$, $\Omega = \R^d_+$ and $\tau(x) = (-x_0, \bx)$, 
a positive Borel measure $\mu$ on $\hat\Omega$ defines a distribution 
$\hat\mu$ on the open half space $\Omega = \R^d_+$ 
if and only if all functions $e_t(p_0,\bp) := e^{-tp_0}$, $t > 0$, on $\hat\Omega$  
are distribution vectors for the 
representation $\pi_\mu^c\res_{\R^{d-1}}$ of the subgroup 
$\{0\} \times \R^{d-1}$ on $L^2(\hat\Omega, \mu)$. 
This in turn is equivalent to the temperedness of the measures 
$\nu_t := \pr_*(e_t \mu)$ on $\R^{d-1}$, where $\pr(x_0,\bx) = \bx$. 
\end{lem}

\begin{prf} In view of Lemma~\ref{lem:auto}, $\hat\mu \in C^{-\infty}(\Omega)$ 
is equivalent to the integrability of all function $\hat\psi$, $\psi \in C^\infty_c(\Omega)$. 
Consider such a function of the form 
$\psi(x_0,\bx) = \psi_0(x) \psi_1(\bx)$ with $\psi_0 \in C^\infty_c(\R)$ 
and $\psi_1 \in C^\infty_c(\R^{d-1})$. Then 
$\hat\psi(p_0,\bp) = \hat\psi_0(p_0) \hat\psi_1(\bp)$  
is $\mu$-integrable on $\hat\Omega$. If $\psi_0 \geq 0$ is supported in 
the interval $[a,b] \subeq \R_+$, then 
\[ \hat\psi_0(p_0) 
= \int_0^\infty e^{-p_0 x_0}\psi_0(x_0)\, dx_0 
\geq e^{-p_0 b} \int_0^\infty \psi_0(x_0)\, dx_0 \] 
implies that $e^{-p_0 b} \hat \psi_1(\bp) \in L^2(\hat\Omega, \mu)$ for every 
$b > 0$ and $\psi_1 \in C^\infty_c(\R^{d-1})$. According to 
Corollary~\ref{cor:3.17}, this means that, for every 
$b > 0$, the function $e^{-b p_0}$ is a distribution vector for the 
representation of $\R^{d-1}$ on $L^2(\hat\Omega, \mu)$. In view of 
Lemma~\ref{lem:dist-vec}, this is equivalent to the temperedness of the measure 
$\pr_*(e^{-b p_0} \mu)$ on $\R^{d-1}$.

Suppose, conversely, that this condition is satisfied and let 
$\phi \in C^\infty_c(\Omega)$ and $0 < a < b$ with 
\[ \supp(\phi) \subeq [a,b] \times \R^{d-1}.\] 
For $\phi_{x_0}(\bx) := \phi(x_0,\bx)$ and 
\[ \hat\phi(p_0,\bp)
= \int_0^\infty \int_{\R^{d-1}} \phi(x_0, \bx) e^{-p_0 x_0} e^{-i \bp \bx}\, d\bx dx_0
= \int_0^\infty \hat\phi_{x_0}(\bp) e^{-p_0 x_0} dx_0,\] 
we obtain the estimate 
$|\hat\phi(p_0,\bp)| \leq e^{-p_0 a} \int_a^b |\hat\phi_{x_0}(\bp)|\, dx_0.$ 
Therefore 
\begin{align*}
\int_{\hat\Omega} |\hat\phi(p_0,\bp)| \, d\mu(p) 
&\leq \int_{\hat\Omega} e^{-p_0 a} \int_a^b |\hat\phi_{x_0}(\bp)|\, dx_0 \, d\mu(p) 
= \int_{\R^{d-1}} \int_a^b |\hat\phi_{x_0}(\bp)|\, dx_0 \, d p_*(e^{-p_0 a}\mu)(\bp) \\
&= \int_a^b \int_{\R^{d-1}} |\hat\phi_{x_0}(\bp)| \, d p_*(e^{-p_0 a}\mu)(\bp)\, dx_0.
\end{align*}

To see that this integral is finite, we first observe that the function 
$\R_+ \to C^\infty_c(\R^{d-1}), x_0 \mapsto  \phi_{x_0}$ is continuous. 
Since the measure $\nu := p_*(e^{-p_0 a}\mu)$ on $\R^{d-1}$ is tempered, 
the map $\gamma \: C^\infty_c(\R^{d-1}) \to L^1(\R^{d-1},\nu), \psi \mapsto \hat\psi$ is continuous, 
and hence 
\[ \int_{\hat\Omega} |\hat\phi(p_0,\bp)| \, d\mu(p) 
\leq \int_a^b \|\gamma(\phi_{x_0})\|_1\, dx_0 < \infty.\qedhere\] 
\end{prf}

\subsection{Representations of vector groups} 

We have seen above that the temperedness of a measure on $\R^d$ 
is closely connected with distribution vectors. 
The following lemma provides a useful criterion to check this 
condition. Its main ingredient is the automatic continuity result 
from Lemma~\ref{lem:auto}.

\begin{lem} \mlabel{lem:dist-vec} Let $(X,\fS,\mu)$ be a measure space. 
We write $M(X,\C)$ for the vector space of measurable 
functions $X \to \C$. For $H_j \in M(X,\R)$, $j = 1,\ldots, d$, we consider 
the continuous unitary representation of $\R^d$ on $L^2(X,\mu)$, given by 
\[ U_{\bt}(f) := e^{i \sum_{j = 1}^d t_j H_j} f\quad \mbox{ for } \quad 
\bt = (t_1,\ldots, t_d).\] 
 Put $R := \sqrt{\sum_{j = 1}^d H_j^2}.$ 
Then 
\[ \cH^{-\infty} \cong \Big\{ h \in M(X,\C) \: (\exists n \in \N) \ \|(1+ R^{2})^{-n} f \|_2 < \infty\Big\},\] 
where the pairing $\cH^{-\infty} \times \cH^\infty \to \C$ is given by 
$(h,f) \mapsto \int_X h \oline f \, d\mu$. 

The following assertions are equivalent: 
\begin{description}
\item[\rm(i)] The constant function $1$ is a distribution vector. 
\item[\rm(ii)] For the measurable map 
$\eta := (H_1, \ldots, H_d) \: X \to \R^d$, 
the measure $\eta_*\mu$ on $\R^d$ is tempered. 
\item[\rm(iii)] $\hat\phi \circ \eta \in L^2(X,\mu)$ for every $\phi \in C^\infty_c(\R^d)$. 
\end{description}
If these conditions are satisfied, then the corresponding distribution 
on $\R^d$ is given by the Fourier transform of 
$\eta_*\mu$. 
\end{lem}

\begin{prf} Put $\cH := L^2(X,\mu)$. 

(a) For $d = 1$ and the measurable function $H \: X \to \R$, 
the subspace of smooth vectors is 
\[ \cH^\infty = \{ f \in \cH \: (\forall n \in \N) \ \|H^n f \|_2 < \infty\}\] 
and the natural Fr\'echet topology on this space is defined by the seminorms 
$p_n(f) := \| (1 + H^2)^n f\|_2.$ 
>From 
\[ \|fh\|_1 = \| f(1+ H^{2})^n (1 + H^{2})^{-n}h\|_1 
\leq p_n(f) \|(1 + H^{2})^{-n}h\|_2 \] 
it follows that the space of distribution vectors can be identified with 
\[ \cH^{-\infty} =  \Big\{ h \in M(X,\C) \: (\exists n \in \N) \ \|(1+ H^{2})^{-n} f \|_2 
< \infty\Big\}.\] 

(b) In the general case it follows from (a) that the subspace of smooth vectors is 
\[ \cH^\infty = \{ f \in \cH \: (\forall \bn \in \N_0^d) 
\ \|H_1^{n_1} \cdots H_d^{n_d}f \|_2 < \infty\}\] 
and the natural Fr\'echet topology on this space is defined by the seminorms 
$p_n(f) := \| (1 + R^{2})^n f\|_2.$ 
Therefore 
\[ \cH^{-\infty} \cong \Big\{ h \in M(X,\C) \: (\exists n \in \N) \ \|(1+ R^{2})^{-n} f \|_2 < \infty\Big\}. \] 
In particular,  $1$ is a distribution vector if and only if some 
function $(1 + R^{2})^{-n}$ is integrable, which is equivalent to the 
integrability of the function 
$(1 + \|x\|^2)^{-n}$ with respect to the measure $\eta_*\mu$ on $\R^d$, 
and this means that $\eta_*\mu$ is tempered. This proves the equivalence of (i) and (ii). 

(ii) $\Rarrow$ (iii) The integrated representation of $L^1(\R^d)$ on $\cH$ is given by 
$U(\phi) f= (\hat \phi \circ (-\eta)) \cdot f$ and 
the integrated representation $U^{-\infty}$ 
on $\cH^{-\infty}$ is given by the same formula. 
Hence (iii) follows from $\hat\phi \circ (-\eta) = U^{-\infty}(\phi) 1 \in L^2(X,\mu)$ 
for $\phi \in C^\infty_c(\R^d)$. 

(iii) First we use the Dixmier--Malliavin Theorem (\cite[Thm.~3.1]{DM78}) to see that 
every test function on $\R^d$ is a finite sum of convolution products 
$\phi * \psi$. Therefore (iii) implies that, for every 
$\phi \in C^\infty_c(\R^d)$, the function $\hat\phi \circ \eta$ is $\mu$-integrable, 
resp., $\hat\phi$ is integrable w.r.t.\ $\eta_*\mu$. 
Now Lemma~\ref{lem:auto}, applied to $\Omega = \R^d$ and 
$\tau = \id_{\R^d}$, implies that $\eta_*\mu$ is tempered. 

For the last assertion, we simply calculate: 
\[ D(\phi) := \la 1, U^{-\infty}(\phi) 1 \ra 
= \int_X \oline{\tilde\phi(\eta(x))}\, d\mu(x) 
= \int_{\R^d} \oline{\tilde\phi(x)}\, d(\eta_*\mu)(x) 
= (\eta_*\mu)\,\hat{}(\phi).\qedhere\]
\end{prf}

Applying Lemma~\ref{lem:dist-vec} to $X = V^* \cong \R^d$ and $\eta(x) = -x$, 
we obtain: 

\begin{cor} \mlabel{cor:3.17} Let $V$ be a finite-dimensional real vector space. 
For a positive Borel measure on $V^*$, the following are equivalent: 
\begin{itemize}
\item[\rm(i)] $\mu$ is tempered. 
\item[\rm(ii)] $1 \in L^2(V^*,\mu)^{-\infty}$ 
for the representation $(U_v f)(\alpha) = e^{-i\alpha(v)} f(\alpha)$. 
\item[\rm(iii)] $\hat\phi \in L^2(V^*,\mu)$ for every $\phi \in C^\infty_c(V)$. 
\end{itemize}
\end{cor}

\begin{rem} The first part of Lemma~\ref{lem:dist-vec} remains true for vector-valued 
$L^2$-spaces $L^2(X,Q,\mu;\cH)$, where $Q\: X \to \Herm_+(\cH)$ 
is an operator-valued density, 
\[ \la f,g \ra = \int_X \la Q(x) f(x), g(x)\ra\, d\mu(x) \quad \mbox{ and }\quad 
U_{\bt}(f) := e^{i \sum_{j = 1}^d t_j H_j} f, \quad  
\bt = (t_1,\ldots, t_d).\]
\end{rem}

\begin{rem} \mlabel{rem:d.10} (Projectable measures and ``time-zero'' subspaces) 

(a) Let $\mu$ be a measure on $\R^d$ and $\pr \: \R^d \to \R^{d-1}, 
p = (p_0, \bp) \mapsto \bp$ be the projection. 
We call $\mu$ {\it projectable} if $\pr_*\mu$ is tempered. In view of 
Corollary~\ref{cor:3.17}, this is equivalent to the integrability 
of all Fourier transforms $\hat\phi \circ \pr$, $\phi \in C^\infty_c(\R^{d-1})$, 
which in turn means that $1$ is a distribution vector for the 
canonical action of $\R^{d-1}$ on $L^2(\R^d,\mu)$. From Lemma~\ref{lem:dist-vec} 
we now see that $1$ is also a distribution vector for $\R^d$, 
which means that $\mu$ is tempered. 

In particular, we obtain an isometric embedding 
\[ \pr^* \: L^2(\R^{d-1}, \pr_*\mu) \into L^2(\R^d,\mu), \quad 
f \mapsto f \circ \pr\] 
which maps onto the subspace of those functions not depending on the first 
argument $p_0$. Accordingly their Fourier transform is supported in 
the $x_0 = 0$ hyperplane. 
Therefore $L^2(\R^{d-1}, \pr_*\mu)$ is also called the {\it time-zero subspace} 
of $L^2(\R^d,\mu)$. 

(b) For the distribution $D := \hat\mu$, the projectability of $\mu$ means that 
$D$ can be restricted to the subspace $\R^{d-1}$ by 
\[ D_r(\phi) := \int_{\R^d} \oline{\hat\phi(\bp)}\, d\mu(p) 
= (\pr_*\mu)\,\hat{}(\phi),\] 
respectively extended to the space $C^\infty_c(\R^{d-1})$. 
The scalar product 
$\la \phi, \psi \ra := D_r(\phi^* * \psi)$ on $C^\infty_c(\R^{d-1})$ leads to 
$\cH_{D_r} \cong L^2(\R^{d-1}, \pr_*\mu)$, which is naturally realized as a subspace 
of $L^2(\R^d,\mu)$. 
Accordingly, we see that, as a space of distributions on $\R^d$, we may consider 
$\cH_{D_r}$ as a subspace of $\cH_D$. 

(c) The Lorentz invariant measures $\mu_m$ are interesting examples of measures 
$\mu$ on $\R^d$ for which the inclusion 
$L^2(\R^{d-1}, \pr_*\mu) \into L^2(\R^{d}, \mu)$ is actually surjective. 

(d) For $d = 1$, a measure $\mu$ is projectable if and only if it is finite, i.e., 
if its Fourier transform $\hat\mu$ is a continuous function. 
Then $L^2(\R^0, \pr_*\mu)$ is the subspace of constant functions in $L^2(\R,\mu)$. 

(e) A similar picture prevails for operator-valued measures, distributions and 
$L^2$-spaces such as those discussed in Section~\ref{sec:dil}. In this context 
the subspace $\cH \subeq L^2(\R,Q;\cH)$ of constant functions plays the role 
of the time-zero subspace. This matches well with Proposition~\ref{prop:3.15a}(iii), 
where $\cE_0 \cong \cH$ is identified as consisting of those functions whose 
Fourier transform is supported in $\{0\}$. 
\end{rem}

\end{document}